%% file: main.tex
\documentclass{article}
\usepackage{style}
\usepackage{caption}
\usepackage{subfig}
\usepackage{authblk}

\usepackage{tikz}
\usepackage{tikz-3dplot}
\usetikzlibrary{shapes,calc,positioning}
\tikzset{every picture/.style={line width=0.75pt}}

\title{\textbf{Exact recovery of the support of piecewise constant images via total variation regularization}}
\author[1]{Yohann De Castro}
\author[2,3]{Vincent Duval}
\author[4]{Romain Petit\thanks{Most of this work was conducted when Romain Petit was affiliated with $2,3$.}}
\affil[1]{Institut Camille Jordan, CNRS UMR 5208, \'Ecole Centrale de Lyon, F-69134 \'Ecully, France}
\affil[2]{CEREMADE, CNRS, UMR 7534, Universit\'e Paris-Dauphine, PSL University, 75016 Paris, France}
\affil[3]{INRIA-Paris, MOKAPLAN, 75012 Paris, France}
\affil[4]{MaLGa Center, Department of Mathematics, University of Genoa, 16146 Genoa, Italy}
\date{\today}

\crefname{equation}{}{}
\newcommand{\customref}[2]{\textup{(}\hyperref[#1]{\textup{#2}}\textup{)}}

\begin{document}
\maketitle
\begin{abstract}
This work is concerned with the recovery of piecewise constant images
from noisy linear measurements. We study the noise robustness of a variational
reconstruction method, which is based on total (gradient) variation
regularization. We show that, if the unknown image is the superposition of a
few simple shapes, and if a non-degenerate source condition holds,
then, in the low noise regime, the reconstructed images have the same
structure: they are the superposition of the same number of shapes, each a
smooth deformation of one of the unknown shapes. Moreover, the reconstructed
shapes and the associated intensities converge to the unknown ones as the noise
goes to zero.
\end{abstract}



\input{intro}
\input{prelim}
\input{faces_tv}
\input{curv_pb}

\input{main_res}

\section*{Conclusion}
We have showed that, in the low noise regime, the support of piecewise constant
images can be exactly recovered from noisy linear measurements, provided that
the measurement operator is smooth enough and some non-degenerate source
condition holds. We have also provided numerical evidence that this last
condition is satisfied for some radial images in the deconvolution setting. 
The investigation of its validity beyond the radial case, which we briefly discussed, 
is an interesting avenue for future research. It is also natural to wonder whether some 
quantitative version of our main result could be proved. This might be achieved
by studying the stability of solutions to the prescribed curvature problem for 
non-smooth perturbations, possibly by adapting the selection principle of 
\cite{cicaleseSelectionPrincipleSharp2012}.
Finally, another direction could be to study the denoising case, which is not
covered by our assumptions. In this setting, dual certificates are a priori non-smooth, 
which is a major source of difficulties.

\section*{Acknowledgments}

The authors warmly thank Jimmy Lamboley and Rapha\"el Prunier for fruitful
discussions about convergence of smooth sets and stability in geometric
variational problems. They are also deeply indebted to Claire Boyer, Antonin
Chambolle, Fr\'ed\'eric de Gournay and Pierre Weiss for insightful discussions
about the faces of the total variation unit ball.

This work was supported by a grant from R\'egion Ile-De-France, by the ANR
CIPRESSI project, grant ANR-19-CE48-0017-01 of the French Agence Nationale de la
Recherche and by the European Union (ERC, SAMPDE, 101041040). Views and opinions
expressed are however those of the authors only and do not necessarily reflect
those of the European Union or the European Research Council. Neither the
European Union nor the granting authority can be held responsible for them.

\bibliography{main}
\bibliographystyle{apalike}

\appendix
\input{appendix}

\end{document}

%% file: intro.tex

\section{Introduction}

\subsection{Reconstruction of images from noisy linear measurements}

In their seminal work \cite{rudinNonlinearTotalVariation1992}, Rudin, Osher and
Fatemi proposed a celebrated denoising method, which has the striking feature of
removing noise from images while preserving their edges. This is achieved by
minimizing a functional with a regularization term, the total~(gradient)
variation, which penalizes oscillations in the reconstructed image, while
allowing for discontinuities. This approach was later applied outside the
denoising setting, in order to solve general linear inverse problems (see e.g.
\cite{acarAnalysisBoundedVariation1994,chaventRegularizationLinearLeast1997}).
Although state of the art algorithms now have much better performance, this work
pioneered the use of the total variation in imaging, and is still an important
baseline for image reconstruction methods.

It is well known that using the total variation as a regularizer promotes
piecewise constant solutions. In denoising, for instance, Nikolova has explained in~\cite{nikolovaLocalStrongHomogeneity2000} that the non-differentiability of the regularizer tends create large flat zones instead of oscillating regions (see also~\cite{ringStructuralPropertiesSolutions2000,jalalzaiRemarksStaircasingPhenomenon2016}), which is known as the \emph{staircasing
effect}. Alternatively, in inverse problems with few linear measurements, it is possible to prove that some solutions to the variational problem are indeed piecewise constant, by appealing to a representation principle derived in \cite{boyerRepresenterTheoremsConvex2019,brediesSparsitySolutionsVariational2019} 
. Considering variational problems with a
convex regularization term, they pointed out the link between the structural
properties of the solutions on the one hand, and the structure of the unit
ball defined by the regularizer on the other hand. In the context of total variation
regularization, these results show that, under a few assumptions, some solutions are of the 
form $\sum_{i=1}^N a_i\mathbf{1}_{E_i}$. This suggests that such functions are the sparse objects naturally associated to this regularizer. In the present article, we follow this line of work and analyze total variation regularization from a
new perspective, by drawing connections with the field of sparse recovery.

\subsection{Problem formulation}
We consider an unknown function $u_0\in\LD$ which models the image to reconstruct. We assume that, in order to recover $u_0$, we have access to a set of linear observations $y_0=\Phi u_0$, where~$\Phi$ is defined by:
$$\! \begin{aligned}[t]
                \Phi : \LD &  \to \mathcal{H} \\
               u   & \mapsto \int_{\RR^2}\phi(x)\,u(x)\,dx\,,
    \end{aligned}$$
with $\phi\in \mathrm{L}^2(\mathbb{R}^2,\mathcal{H})$ and $\mathcal{H}$ a separable Hilbert space (typically $\RR^m$ or $\LD$). To account for the presence of noise in the observations, we also consider the recovery of $u_0$ from~$y_0+w$ where~$w\in\mathcal{H}$ is an additive noise. Following the above-mentioned works, we aim at recovering~$u_0$ from $y_0$ by solving
\begin{equation}
   \underset{u \in \LD}{\text{inf}}~\TV(u) ~\text{ s.t. }~ \Phi u=y_0\,,
  \tag{$\mathcal{P}_0(y_0)$}
  \label{primal_noreg}
\end{equation}
where $\TV(u)$, defined below, denotes the total (gradient) variation of $u$. To recover $u_0$ from~${y_0+w}$, we solve instead, for some $\lambda>0$, the following problem:
\begin{equation}
    \underset{u \in\LD}{\text{inf}} ~ \frac{1}{2}\|\Phi u-
    y\|_{\mathcal{H}}^2+\lambda\,\TV(u)\,,
    \tag{$\mathcal{P}_{\lambda}(y)$}
    \label{primal_reg}
\end{equation}
with $y=y_0+w$.

The question this work is concerned with is the following: if $w$ is small and $\lambda$ well chosen,  are the solutions of \customref{primal_reg}{$\mathcal{P}_{\lambda}(y_0+w)$} close to some solutions of \cref{primal_noreg}? If $u_0$ is the unique solution to~\cref{primal_noreg}~(in this case, we say that $u_0$ is \emph{identifiable}), answering positively amounts to proving that the considered variational method enjoys some noise robustness, i.e. that solving~\customref{primal_reg}{${\mathcal{P}_{\lambda}(y_0+w)}$} yields good approximations of $u_0$ in the low noise regime.

To our knowledge, finding sufficient conditions for the identifiablity of $u_0$ is mostly open. Still, let us point out that, in \cite{brediesPerfectReconstructionProperty2019}, an identifiability result is obtained for the recovery of $u_0$ from its image under a linear partial differential operator with unknown boundary conditions. An exact recovery result is also obtained in \cite{hollerExactReconstructionReconstruction2022} in a different setting, where the regularizer is the so-called \emph{anisotropic} total variation.

\subsection{Motivation}
\label{sec_motiv}
In order to motivate our analysis, we present in this subection a simple experiment showcasing the phenomenon we wish to analyze. We consider an unknown image of the form~${u_0=\sum_{i=1}^N a_i\mathbf{1}_{E_i}}$, and define $\Phi$ as the convolution with a Gaussian filter followed by a subsampling on a regular grid of size $50\times 50$. The noise is drawn from a multivariate Gaussian with a zero mean and an isotropic covariance matrix. Given these noisy observations $y_0+w$, we numerically approximate a solution $u_{\lambda,w}$ of \customref{primal_reg}{$\mathcal{P}_{\lambda}(y_0+w)$} using the method introduced in \cite{condat_discrete_2017}, which is a discrete image defined a grid $5$ times finer than the observation grid. We notice that, for two different choices of $u_0$, the approximation of $u_{\lambda,w}$ has a structure which is close to that of $u_0$. Up to discretization artifcats, it is the superposition of the same number of shapes, each being close to one of the unknown shapes.

In the present article, we wish to theoretically analyze this phenomenon. Our aim is to investigate whether, if $u_0$ is the superposition of a few simple shapes, solutions of \customref{primal_reg}{$\mathcal{P}_{\lambda}(y_0+w)$} have the same structure.

\begin{figure}
    \newcommand*{\x}{3.1cm}
    \centering
    \includegraphics[height=\x]{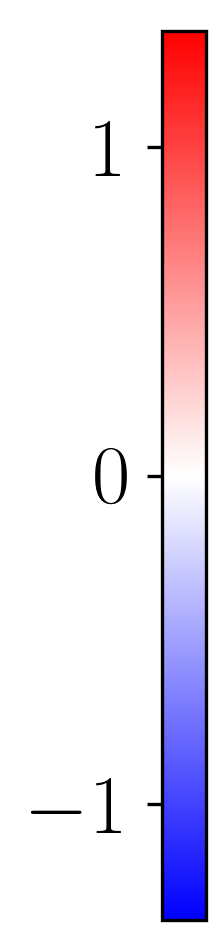}
    \subfloat[\normalsize $u_0$]{\includegraphics[height=\x]{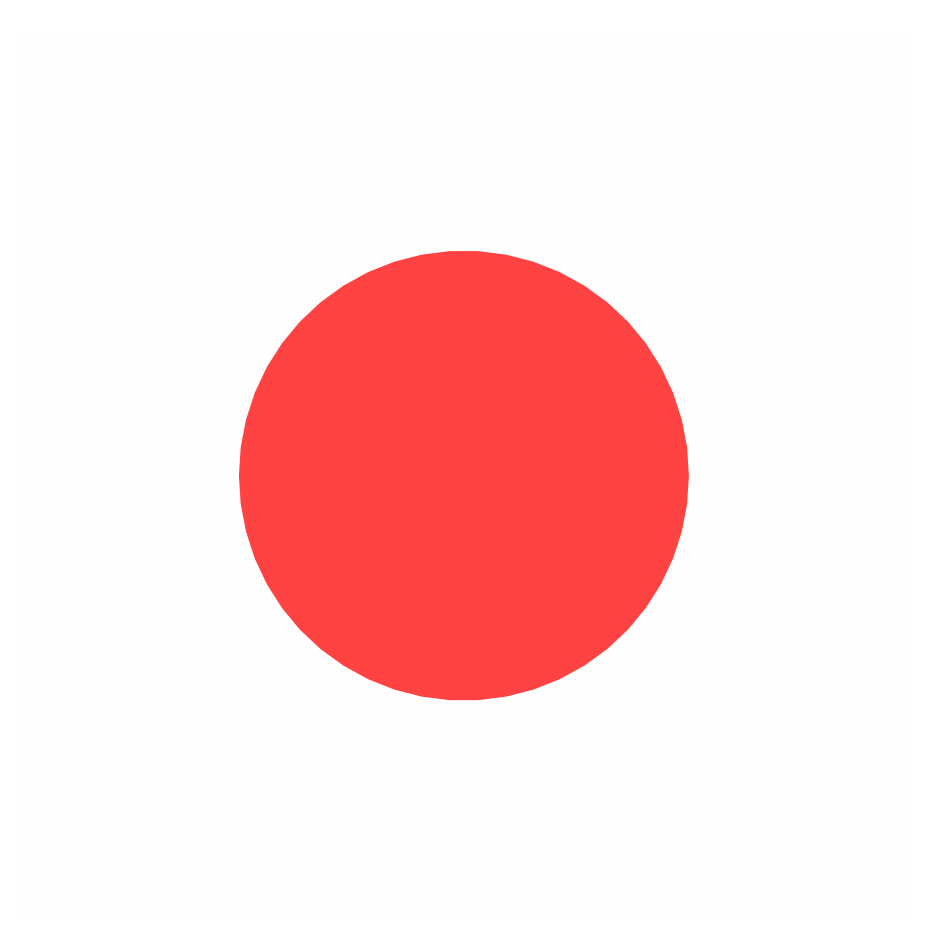}}
    \quad
    \subfloat[\normalsize $y_0+w$]{\includegraphics[height=\x]{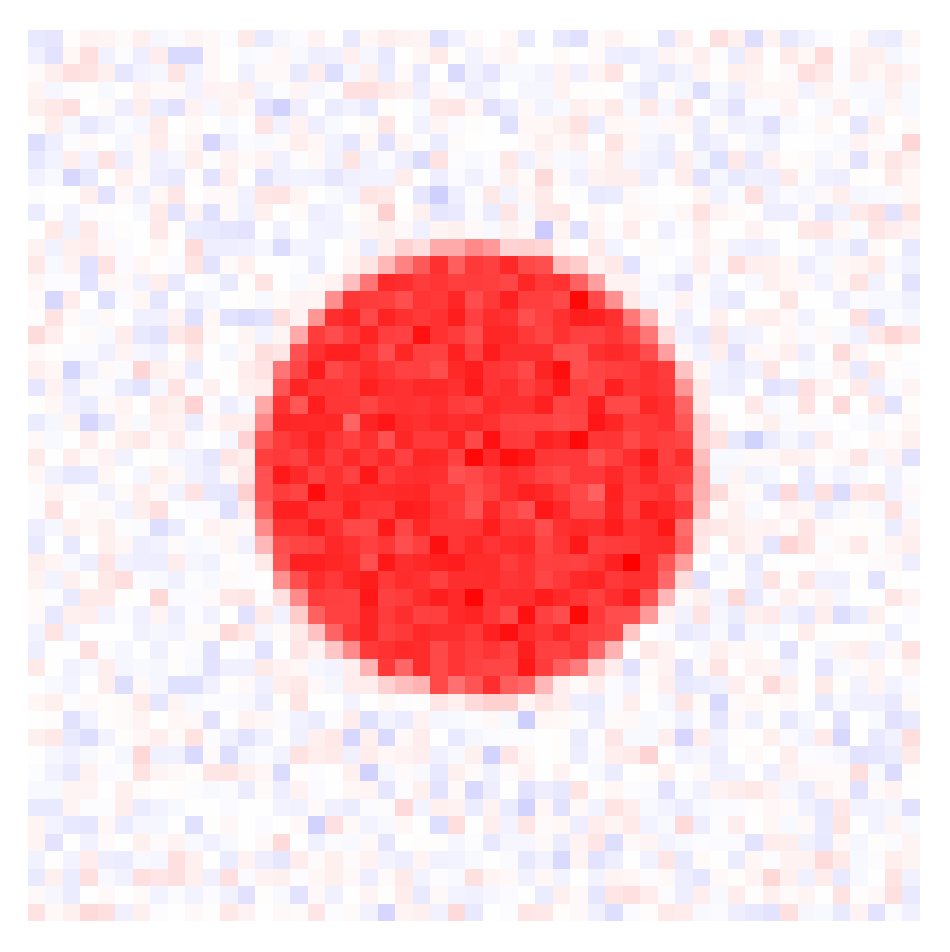}}
    \quad
    \subfloat[\normalsize approx. of $u_{\lambda,w}$]{\includegraphics[height=\x]{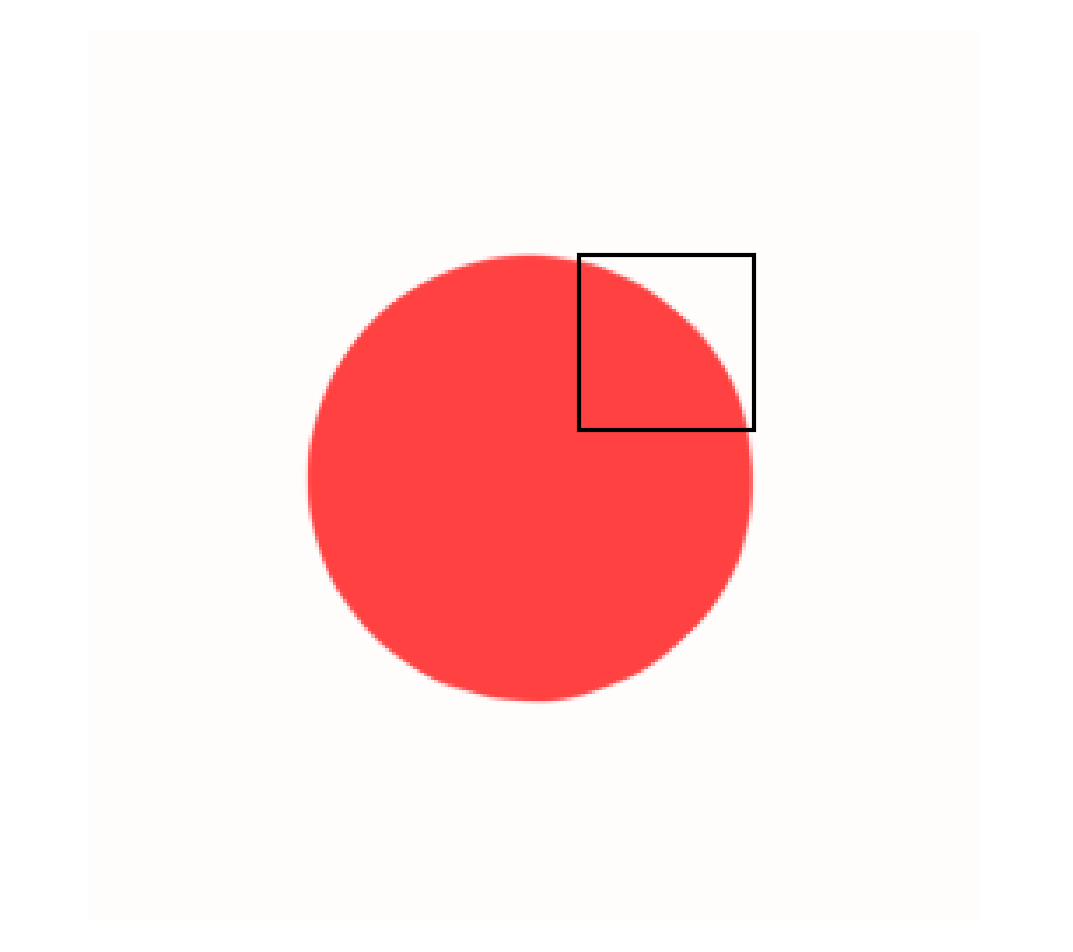}}
    \quad
    \subfloat{\includegraphics[height=\x]{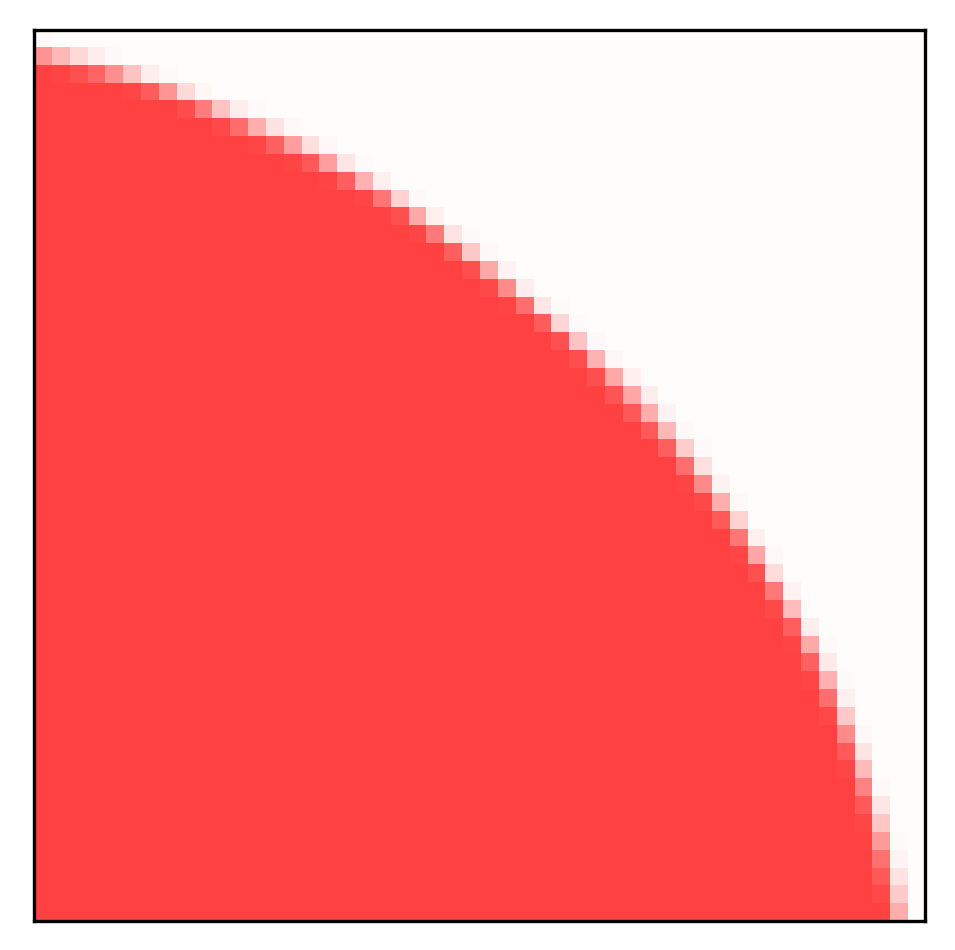}}
    \caption{Numerical resolution of \customref{primal_reg}{$\mathcal{P}_{\lambda}(y_0+w)$} for $u_0=\mathbf{1}_{B(0,R)}$.}
    \label{condat_disk}
 \end{figure}

\begin{figure}
    \newcommand*{\x}{3.1cm}
    \centering
    \includegraphics[height=\x]{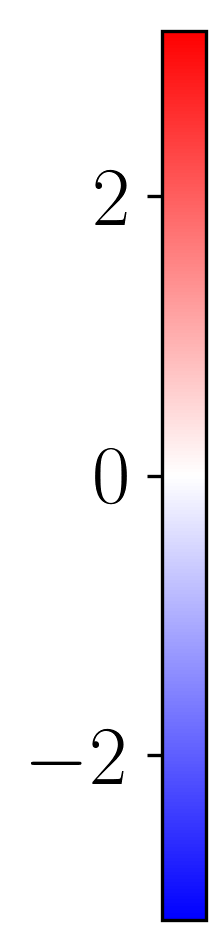}
    \,
    \subfloat[\normalsize $u_0$]{\includegraphics[height=\x]{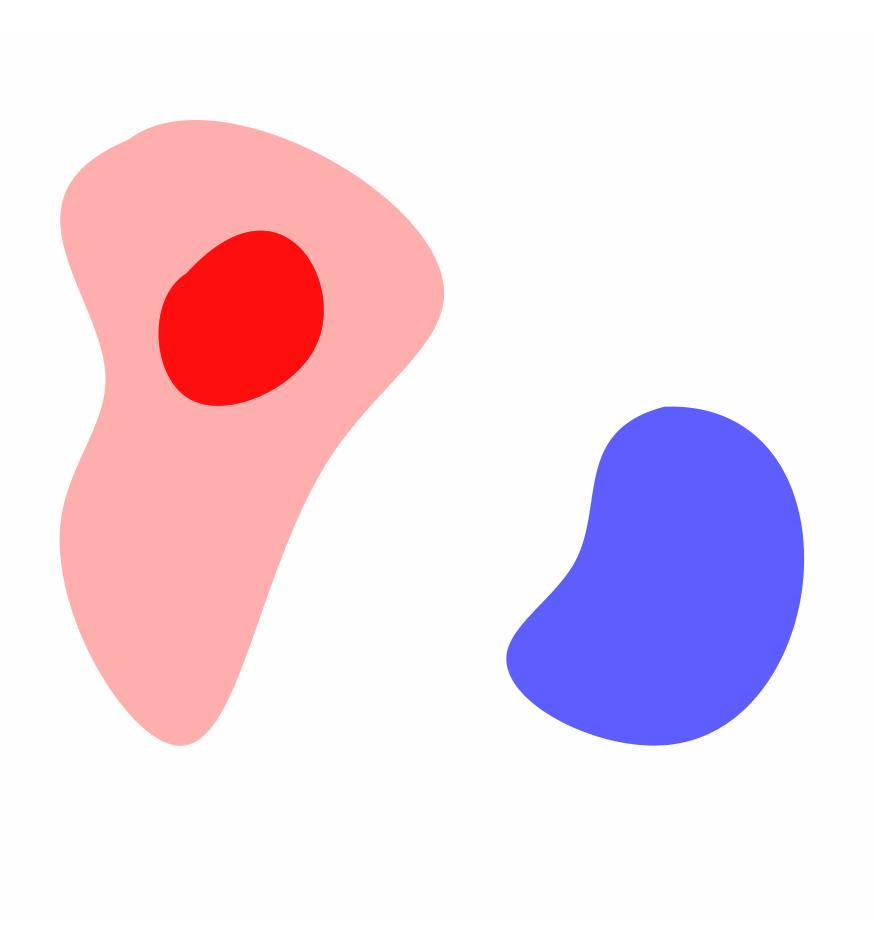}}
    \quad\,
    \subfloat[\normalsize $y_0+w$]{\includegraphics[height=\x]{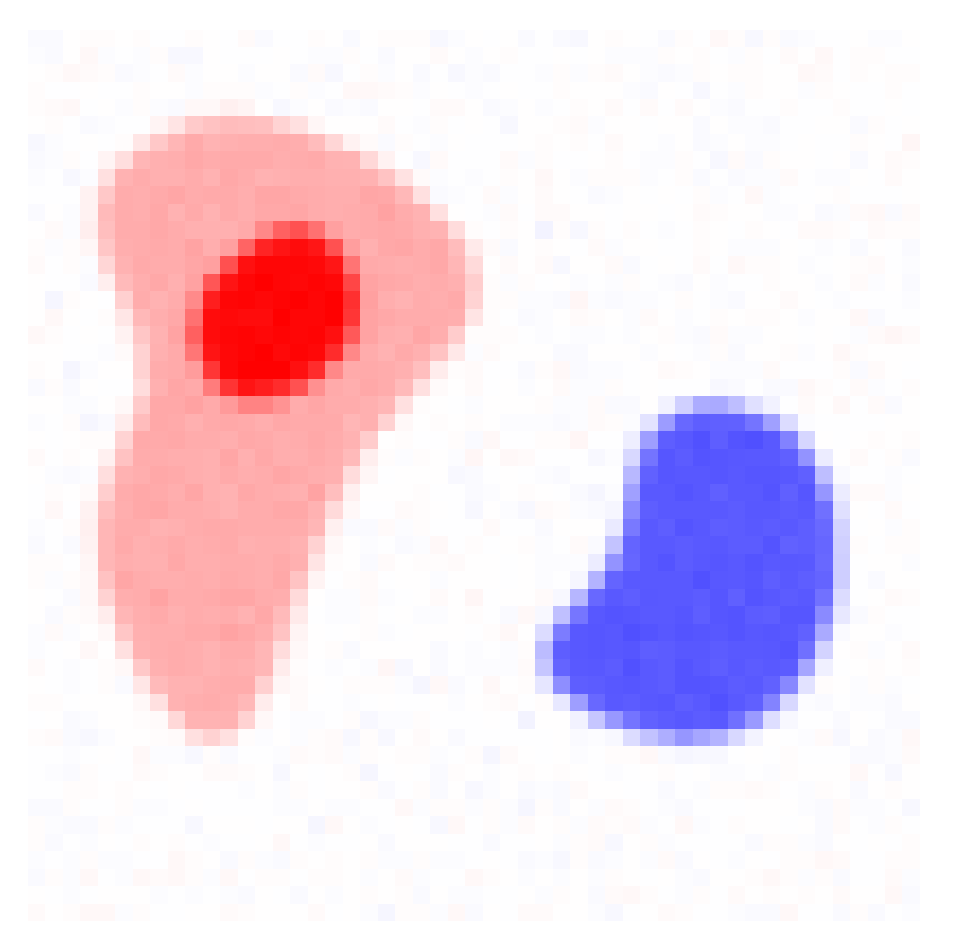}}
    \quad\,
    \subfloat[\normalsize approx. of $u_{\lambda,w}$]{\includegraphics[height=\x]{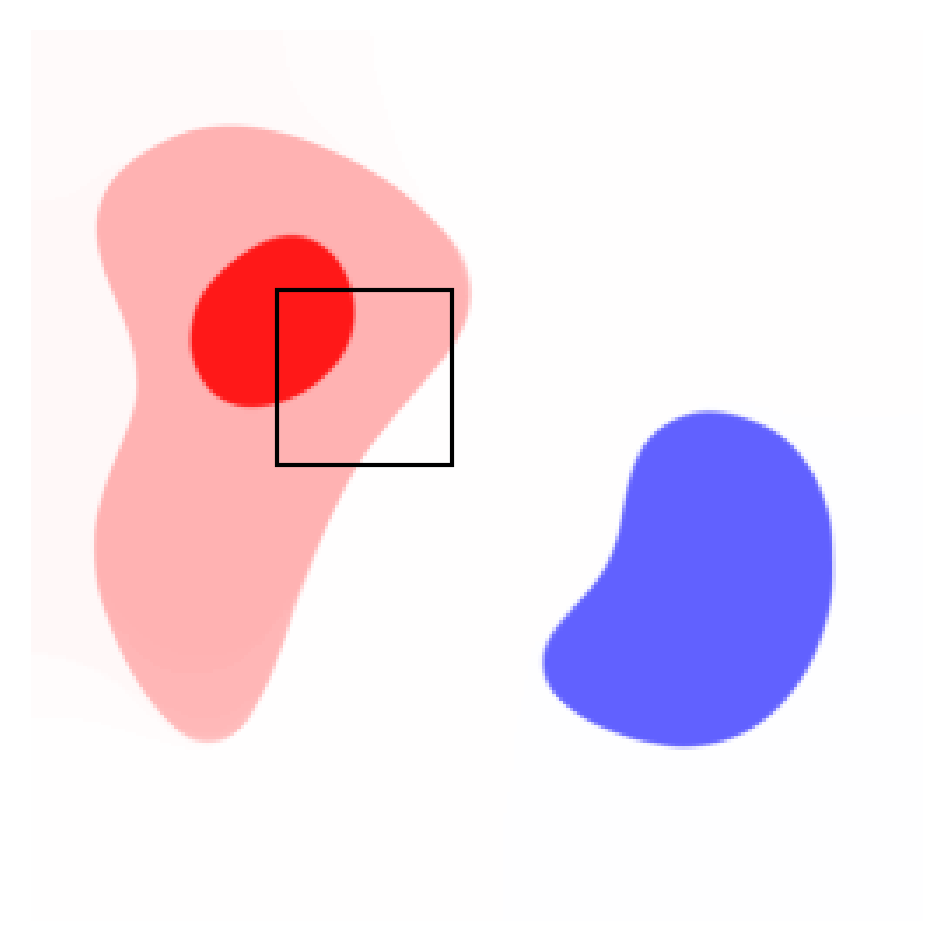}}
    \quad\,
    \subfloat{\includegraphics[height=\x]{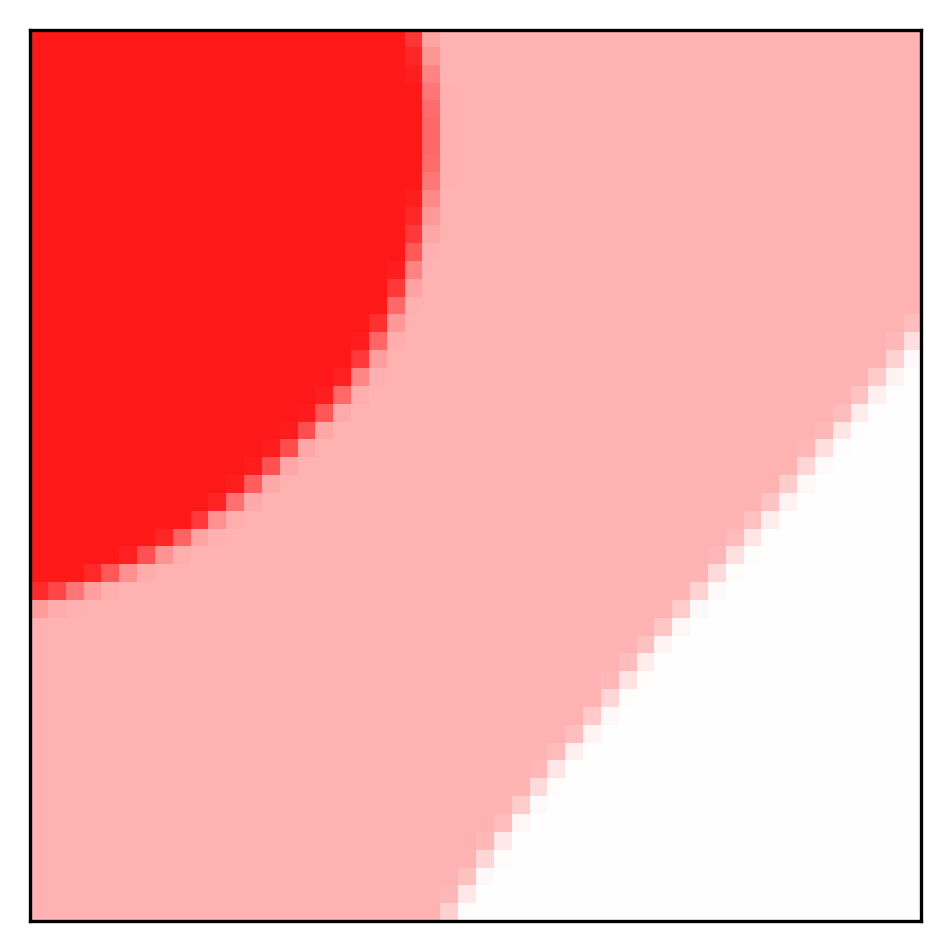}}
    \caption{Numerical resolution of \customref{primal_reg}{$\mathcal{P}_{\lambda}(y_0+w)$} for $u_0=\mathbf{1}_{E_1}+2\mathbf{1}_{E_2}-2\mathbf{1}_{E_3}$.}
    \label{condat_three_shapes}
 \end{figure}

 \subsection{Previous works}
\paragraph{Total variation minimization in imaging.}  The theoretical study of total
variation regularization in imaging was initiated in~\cite{chambolleImageRecoveryTotal1997,ringStructuralPropertiesSolutions2000}. Then, a lot of attention was focused on the denoising case, which can be regarded as one step of the total variation gradient flow, see~\cite{bellettiniTotalVariationFlow2002,alterEvolutionCharacteristicFunctions2005}. Its connection to the Cheeger problem was observed in~\cite{alterCharacterizationConvexCalibrable2005,alteruniq09} and was the key to  understanding  the properties of the Cheeger sets of convex bodies. Let us also mention the 
landmark result~\cite{casellesDiscontinuitySetSolutions2007}, which shows that the jump set of the reconstructions 
is included in the jump set of the noisy input. Independently, Allard gave a precise description of the properties of the minimizers in the series of articles~\cite{allard_total_2008,allard_total_2008-1,allard_total_2009}. We refer to
\cite{chambolleIntroductionTotalVariation2010} for an introduction to total variation references and a much more comprehensive list of bibliographical references.

\paragraph{Piecewise constant images.} Let us emphasize that the above-mentioned representation principle \cite{boyerRepresenterTheoremsConvex2019,brediesSparsitySolutionsVariational2019} only applies to inverse problems with a finite number of measurements. It does not cover the case of denoising or related tasks involving infinite-dimensional observations. For such cases, several authors have proposed specific approaches to promote piecewise constant solutions \cite{fornasier2006nonlinear,fornasier2007restoration,fonseca2010exact,cristoferi2019piecewise}. These methods rely on the minimization of non-convex functionals which resemble the $\ell^0$ norm, while the total variation relates to the $\ell^1$ norm.

\paragraph{Noise robustness.} The general convergence results presented in
\cite{burgerConvergenceRatesConvex2004,hofmannConvergenceRatesResult2007} apply
to the case of total variation regularization, and loosely speaking provide
(under mild assumptions) strict convergence in $\mathrm{BV}_{\mathrm{loc}}$ of
solutions of \customref{primal_reg}{$\mathcal{P}_{\lambda}(y_0+w)$} towards
solutions of \cref{primal_noreg}. Moreover, in specific cases, the analysis in
\cite{burgerConvergenceRatesConvex2004} ensures that the variation of solutions
to $\customref{primal_reg}{$\mathcal{P}_{\lambda}(y_0+w)$}$ is mostly
concentrated in a neighborhood of the support of $\Diff u_0$. In
\cite{chambolleGeometricPropertiesSolutions2016,iglesiasNoteConvergenceSolutions2018},
improved convergence guarantees are derived by exploiting the optimality of
level sets of solutions of
\customref{primal_reg}{$\mathcal{P}_{\lambda}(y_0+w)$} for the prescribed
curvature problem. The main finding of these works is that, under a few
assumptions, the boundaries of the level sets converge in the Hausdorff sense.

\paragraph{Support recovery in variational sparse regularization.}
In the litterature on variational sparse regularization, a line of work has specifically focused on strong noise robustness guarantees, called \emph{exact support recovery} or \emph{model identification}. They were first studied in the case of finite-dimensional $\ell^1$ regularization. In this setting, it was proved in~\cite{fuchsSparseRepresentationsArbitrary2004} and \cite{troppJustRelaxConvex2006} that, under some conditions (called ``irrepresentability conditions'' in the statistics community, see \cite{zhaoModelSelectionConsistency2006}), solving  the celebrated \textsc{Lasso} problem (the analog of \eqref{primal_reg}) allows to exactly recover the support of the unknown sparse vector. Subsequently, it was shown in the landmark paper \cite{candesRobustUncertaintyPrinciples2006} that these conditions are satisfied in the case of subsampled Fourier measurements. Later works established the exact recovery of the group-support for group-sparse vectors \cite{bachConsistencyGroupLasso2008}, of the rank for low-rank matrices \cite{bachConsistencyTraceNorm2008}, or of more general structures encoded in the singularities of partly-smooth regularizers \cite{vaiterModelConsistencyPartly2018} or mirror-stratifiable functions \cite{fadiliSensitivityAnalysisMirrorStratifiable2018}.

\paragraph{Off-the-grid sparse spikes recovery and non-degeneracy assumptions.} Beyond finite-dimensional examples, several works have investigated the case of total variation regularization on the space of Radon measures, where the unknown signal is a linear combination of Dirac masses (see \cite{decastroExactReconstructionUsing2012} and~\cite{brediesInverseProblemsSpaces2013}). In that setting (contrary to that of the present paper), sufficient identifiability conditions have been
extensively studied (see the landmark
paper \cite{candesMathematicalTheorySuperresolution2014} and its generalization~\cite{poonGeometryOfftheGridCompressed2023}). In \cite{duvalExactSupportRecovery2015}, a slighty stronger condition was introduced, called the \emph{non-degenerate source condition}, which ensures the exact support recovery property in this infinite-dimensional setting. More precisely, at low noise, the solutions to the Beurling \textsc{Lasso} problem (the analog of \eqref{primal_reg} for measures) have the same number of atoms as the unknown measure, each close to one of the unknown atoms. Ensuring this condition is a delicate issue. It was first investigated numerically in~\cite{duvalExactSupportRecovery2015} for the deconvolution by the ideal low-pass filter with cutoff frequency~$f_c$. Numerical evidence suggests that it is satisfied for measures whose atoms are separated by a distance of at least $O(1/f_c)$. The special case of positive atoms is somewhat simpler, and was investigated in dimension one in~\cite{denoyelleSupportRecoverySparse2017,duvalCharacterizationNonDegenerateSource2019} and in higher dimension in~\cite{poonMultidimensionalSparseSuperResolution2019}. Eventually, a general theoretical framework ensuring that condition provided the locations of the Dirac masses are sufficiently separated was provided in~\cite{poonGeometryOfftheGridCompressed2023}. One striking feature of the \emph{non-degenerate source condition}, compared to the standard $\ell^1$ identifiability condition, is that it involves the second order derivatives of some solution of the dual problem. It is not known whether the stability of the support holds under weaker assumptions, but there are counterexamples of measures which are identifiable (i.e., they are the unique solution to the noiseless problem) but their support is not stable at any noise level, see \cite[Chapter 5]{duvalFacesExtremePoints2022}.
Thus, some additional assumption is needed, and the condition on the second derivatives fills that gap. It is therefore natural that, in the present paper, our analysis of the total (gradient) variation support recovery property relies on the non-degeneracy of some second derivatives. 
Let us mention, incidentally, that the recent work~\cite{carioni2023general}, which studies the stability of sparse representations in convex optimization involves a similar criterion\footnote{However, as underlined by the authors of \cite{carioni2023general}, their analysis does not cover our setting.}.

\subsection{Contributions}
The above-mentioned works on total (gradient) variation regularization provide little information about the structure of
solutions to \customref{primal_reg}{${\mathcal{P}_{\lambda}(y_0+w)}$} in the low
noise regime. However, in light of the numerical evidence presented in \Cref{sec_motiv}, the following question is natural: if~$u_0$ is identifiable and is the sum of a few indicator functions, do the solutions
of~\customref{primal_reg}{$\mathcal{P}_{\lambda}(y_0+w)$} have a similar
property? Moreover, are these decompositions stable, \ie{} are they made of
the same number of atoms, and are their atoms related?

In this work, we answer these questions by using two main tools. The first is a
set of results about the faces of the unit ball defined by the total variation,
which provide useful information on the above-mentioned decompositions. The
second is an analysis of the behaviour of solutions to the prescribed curvature
problem under variations of the curvature functional. 

To state our main result,
we introduce a non-degenerate version of the source condition, which relies on a
regularity assumption on the measurement operator~$\Phi$, namely
that~$\phi\in\mathrm{C}^1(\RR^2,\mathcal{H})$.
This covers a wide variety of forward operators that are of particular interest for applications. To name only a few, let us cite the case of the convolution with a Gaussian blur (or more generally, with any $\mathrm{C}^1$ filter) possibly followed by a subsampling,  or a motion blur, of particular interest in computational photography. Let us also mention the case of subsampled Fourier measurements, of particular interest in X-ray tomography, magnetic resonance imaging and radio interferometry. However, it does not cover the case of denoising,\footnote{Although it is counterintuitive, we conjecture that  the denoising case is less favorable to support recovery. For instance, it is known that when $u_0$ is the indicator of a square (with rounded corners), the solution of \eqref{primal_reg} has infinitely many level lines for any $\lambda>0$ small enough, see \cite{chambolleIntroductionTotalVariation2010}.} which corresponds to~$\mathcal{H}=\LD$ and~$\Phi=\mathrm{Id}$. In all the following, except in \Cref{sec_prelim} in which we review existing results, we assume that~$\phi\in\mathrm{C}^1(\RR^2,\mathcal{H})$.

\begin{figure}
   \centering
   \subfloat[\normalsize$u_0$]{\input{supprec1}}
   \qquad
   \subfloat[\normalsize$u_{\lambda,w}$]{\input{supprec2}}
   \caption{Illustration of the result stated in \Cref{main_th}. Here $u_0$ is equal to $0$ in $A$ and $B$. The values taken by $u_{\lambda,w}$ in $A^{\lambda,w}$ and $B^{\lambda,w}$ are close (but not necessarily equal) to $0$.}
   \label{tikz_mainres}
\end{figure}
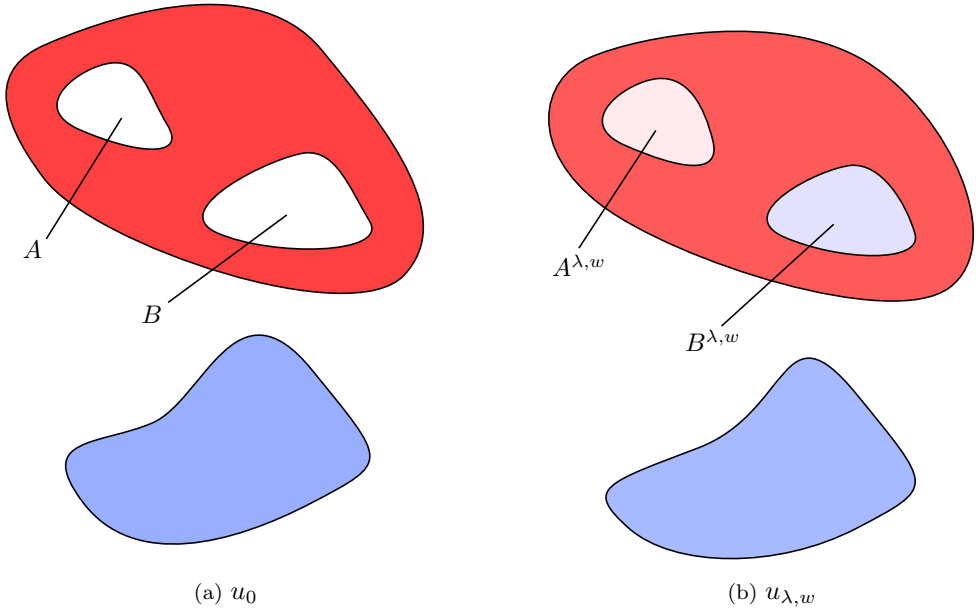

Our main result, which is \Cref{main_th}, informally states that, if the unknown
image modeled by~$u_0$ is the superposition of a few simple shapes and the
non-degenerate source condition holds, then, in the low noise regime, every
solution $u_{\lambda,w}$ of
\customref{primal_reg}{$\mathcal{P}_{\lambda}(y_0+w)$} is made of the same
number of shapes as $u_0$, each shape in $u_{\lambda,w}$ converging smoothly to
the corresponding shape in~$u_0$ as the noise goes to zero (see
\Cref{tikz_mainres} for an illustration).

%% file: supprec1.tex
\scalebox{0.8}{
\begin{tikzpicture}[x=0.75pt,y=0.75pt,yscale=-1,xscale=1]

\draw  [color={rgb, 255:red, 0; green, 0; blue, 0 }  ,draw opacity=1 ][fill={rgb, 255:red, 255; green, 0; blue, 0 }  ,fill opacity=0.75 ] (218.14,33.24) .. controls (249.77,19.25) and (346.79,-18.06) .. (391.09,36.03) .. controls (435.38,90.12) and (475.45,142.35) .. (442.76,177.79) .. controls (410.07,213.23) and (249.77,159.14) .. (218.14,117.17) .. controls (186.5,75.2) and (186.5,47.22) .. (218.14,33.24) -- cycle ;
\draw  [fill={rgb, 255:red, 255; green, 255; blue, 255 }  ,fill opacity=1 ] (260.32,44.43) .. controls (279.3,42.56) and (282.46,62.15) .. (295.12,83.6) .. controls (307.77,105.05) and (274.03,101.32) .. (240.28,86.39) .. controls (206.54,71.47) and (241.34,46.29) .. (260.32,44.43) -- cycle ;
\draw  [fill={rgb, 255:red, 255; green, 255; blue, 255 }  ,fill opacity=1 ] (381.09,100.98) .. controls (400.07,99.11) and (408.5,122.43) .. (421.16,143.88) .. controls (433.81,165.33) and (359.99,166.26) .. (326.25,151.34) .. controls (292.5,136.42) and (362.1,102.84) .. (381.09,100.98) -- cycle ;
\draw  [fill={rgb, 255:red, 115; green, 146; blue, 255 }  ,fill opacity=0.73 ] (285.82,270.67) .. controls (317.45,256.69) and (337.49,178.35) .. (384.95,237.1) .. controls (432.4,295.85) and (434.51,294.92) .. (384.95,321.03) .. controls (335.38,347.15) and (274.22,363) .. (242.58,321.03) .. controls (210.94,279.07) and (254.18,284.66) .. (285.82,270.67) -- cycle ;
\draw    (266.5,79) -- (219.5,155) ;
\draw    (369.5,140) -- (295.5,195) ;

\draw (203,154.4) node[scale=1.2] [anchor=north west][inner sep=0.75pt]    {$A$};
\draw (277,195.4) node[scale=1.2] [anchor=north west][inner sep=0.75pt]    {$B$};
\end{tikzpicture}}

%% file: supprec2.tex
\scalebox{0.8}{
\begin{tikzpicture}[x=0.75pt,y=0.75pt,yscale=-1,xscale=1]

\draw  [color={rgb, 255:red, 0; green, 0; blue, 0 }  ,draw opacity=1 ][fill={rgb, 255:red, 255; green, 0; blue, 0 }  ,fill opacity=0.65 ] (216.14,24.24) .. controls (247.77,10.25) and (336.5,-6) .. (389.09,27.03) .. controls (441.67,60.07) and (477.02,141.58) .. (440.76,168.79) .. controls (404.5,196) and (269.5,152) .. (223.5,117) .. controls (177.5,82) and (184.5,38.22) .. (216.14,24.24) -- cycle ;
\draw  [fill={rgb, 255:red, 226; green, 226; blue, 255 }  ,fill opacity=1 ] (379.09,91.98) .. controls (398.07,90.11) and (413.82,115.76) .. (419.16,134.88) .. controls (424.5,154) and (371.25,151.92) .. (337.5,137) .. controls (303.75,122.08) and (360.1,93.84) .. (379.09,91.98) -- cycle ;
\draw  [fill={rgb, 255:red, 115; green, 146; blue, 255 }  ,fill opacity=0.63 ] (283.82,269.67) .. controls (342.13,247.35) and (335.49,177.35) .. (382.95,236.1) .. controls (430.4,294.85) and (432.51,293.92) .. (382.95,320.03) .. controls (333.38,346.15) and (268.5,345) .. (240.58,320.03) .. controls (212.66,295.07) and (225.5,292) .. (283.82,269.67) -- cycle ;
\draw  [fill={rgb, 255:red, 255; green, 235; blue, 235 }  ,fill opacity=1 ] (258.32,37.43) .. controls (277.3,35.56) and (287.74,54.19) .. (293.12,76.6) .. controls (298.5,99) and (272.03,94.32) .. (238.28,79.39) .. controls (204.54,64.47) and (239.34,39.29) .. (258.32,37.43) -- cycle ;
\draw    (257.5,70) -- (209.5,144) ;
\draw    (368.5,129) -- (298.5,193) ;

\draw (191,146.4) node[scale=1.2] [anchor=north west][inner sep=0.75pt]    {$A^{\lambda ,w}$};
\draw (274,193.4) node[scale=1.2] [anchor=north west][inner sep=0.75pt]    {$B^{\lambda ,w}$};
\end{tikzpicture}
}

%% file: prelim.tex

\section{Preliminaries}
\label{sec_prelim}
\subsection{Smooth sets and normal deformations}
\label{prelim_smooth}
Our analysis mainly concerns the level sets of the solutions to \eqref{primal_noreg} and~\eqref{primal_reg}. As we strongly rely on their regularity, we recall here several definitions and properties related to smooth sets and their normal deformations. We refer to~\cite{delfourShapesGeometriesMetrics2011} for more details.

\paragraph{Smooth set.} Let $E \subset \RR^2$ be an open set such that $\partial E\neq \emptyset$ and $k \in \NN^*$, where $\NN^*\eqdef\{1,2,...\}$ denotes the set of positive integers. We say that the set~$E$ is of class $\mathrm{C}^k$ if,
for every $x \in \partial E$, there exists~$r_x>0$, a rotation matrix $R_x$, and a function~$u_x \in \Cder^k([-r_x,r_x])$ such that
\begin{align*}
	\left\{
		\begin{aligned}
	 R_x^{-1}(E-x)\cap C(0,r_x) &= \enscond{(z,t) \in C(0,r_x)}{t<u_x(z)}\eqdef\mathrm{hypograph}(u_x)\,,\\
	 R_x^{-1}(\partial E-x)\cap C(0,r_x) &= \enscond{(z,t) \in C(0,r_x)}{t=u_x(z)}\eqdef \mathrm{graph}(u_x)\,,
	\end{aligned}
	\right.
\end{align*}
where $C(0,r)\eqdef(-r,r)^2$.
In that case, one can choose $u_x(0)=0$ and $\nabla u_x(0)=0$. Moreover, if~$\partial E$ is compact, $r_x$ can be taken independent of $x$, and the family $\{u_{x}\}_{x\in \partial E}$ uniformly equicontinuous~(see~\cite[Theorem 5.2]{delfourShapesGeometriesMetrics2011}). In local coordinates, the outward unit normal to $E$ at $(z,u_x(z))$ is given by
\begin{equation*}
	\nu_E(z,u_x(z))=\frac{1}{\sqrt{1+u_x'(z)^2}}\begin{pmatrix}-u_x'(z)\\1\end{pmatrix}\,.
\end{equation*}
It is a geometric quantity, which does not depend on the choice of $r$, $x$ and $u_x$. Likewise, the signed curvature of $E$ at $(z,u_x(z))$ is given by
\begin{equation*}
	H_E(z,u_x(z))=\left(\frac{-u_x'}{\sqrt{1+u_x'^2}}\right)'(z)=\frac{-u_x''(z)}{(1+u'(z)^2)^{3/2}}\,.
\end{equation*}

\begin{remark}
	The same definitions and properties hold when replacing $\Cder^k$ with $\Cder^{k,\ell}$ the space of~$k$-times continuously differentiable functions whose $k$-th derivative is $\ell$-H\"{o}lder ($0< \ell\leq 1$).
\end{remark}

\paragraph{Lebesgue equivalence classes and smooth sets.} If $E\subset \RR^2$ is an open set of class $\mathrm{C}^1$ and~$x\in\RR^2$, then its Lebesgue density exists everywhere and it is given by
\begin{align}
	\theta_E(x)\eqdef \underset{r\to 0^+}{\mathrm{lim}}~\frac{|E\cap B(x,r)|}{|B(x,r)|}
	=\begin{cases}
		1 & \mbox{if $x \in E$},\\
		1/2 & \mbox{if $x \in \partial E$},\\
		0 & \mbox{if $x \in \RR^2\setminus \overline{E}$.}
	\end{cases}
	\label{def_density}
\end{align}
When working with measurable sets, it is common to regard them \emph{modulo Lebesgue negligible sets}.
The above equality shows that if a measurable set $\tilde{E} \in \RR^2$ is equivalent to a $\mathrm{C}^1$ open set~$E$, then $E$ is unique and can be recovered as the set of Lebesgue points of $\tilde{E}$, that is~${\{\theta_{\tilde{E}}=1\}}$. In the following, we usually work with Lebesgue equivalence classes. When $E$ has a~$\mathrm{C}^1$ open representative, we say that $E$ is of class $\mathrm{C}^1$, and we denote by $\partial E$ the topological boundary of that representative.

\paragraph{Convergence.} Let us define the square of axis $\nu\in\mathbb{S}^1$ and side $r>0$ centered at $x\in\RR^2$:
\begin{equation}
	C(x,r,\nu)\eqdef x+R_{\nu}\,C(0,r)\,.
	\label{def_cylinder}
\end{equation} 
where $R_\nu\colon \mathbb{R}^2\rightarrow \mathbb{R}^2$ is the rotation that maps $(0,1)$ to $\nu$.

Let $E$ be a set of class $\mathrm{C}^k$ such that~$\partial E$ is compact. We say that a sequence $(E_n)_{n\in\NN}$ converges to~$E$ in $\mathrm{C}^k$ if there exists~${r>0}$ and $n_0\in\NN$ such that
\begin{itemize}
	\item for every $n\geq n_0$ we have $\partial E_n\subset \bigcup\limits_{x\in\partial E}C(x,r,\nu_E(x))$
	\item for every $n\geq n_0$ and $x\in\partial E$ there exists $u_{n,x}\in \mathrm{C}^{k}([-r,r])$ such that:
	$$\left\{\begin{aligned} R_{\nu_E(x)}^{-1}\left(\partial E_n -x\right)\cap C(0,r)&=\mathrm{graph}(u_{n,x})\\
		R_{\nu_E(x)}^{-1}\left(\mathrm{int}\,E_n -x\right)\cap C(0,r)&=\mathrm{hypograph}(u_{n,x})\end{aligned}\right.$$
	\item denoting $(u_x)_{x\in\partial E}$ some functions satisfying
	$$\left\{\begin{aligned} R_{\nu_E(x)}^{-1}\left(\partial E -x\right)\cap C(0,r)&=\mathrm{graph}(u_{x})\\
		R_{\nu_E(x)}^{-1}\left(\mathrm{int}\,E -x\right)\cap C(0,r)&=\mathrm{hypograph}(u_{x})\end{aligned}\right.$$
	we have $\underset{n\to +\infty}{\mathrm{lim}}~\underset{x\in\partial E}{\mathrm{sup}}~\|u_{n,x}-u_x\|_{\mathrm{C}^{k}([-r,r])}=0$
\end{itemize}

\paragraph{Normal deformation.} We state below some useful results regarding normal deformations of a smooth set $E$. First, let us stress that such sets are parametrized by real-valued functions on~$\partial E$, which leads us to use the notion of tangential gradient, tangential Jacobian, and the spaces $\mathrm{C}^k(\partial E)$, $\mathrm{L}^p(\partial E)$ and $\mathrm{H}^1(\partial E)$ (along with their associated norms). We refer to the reader to~\cite[Section 5.4.1, 5.4.3 and 5.9.1]{henrotShapeVariationOptimization2018} for precise definitions.
\begin{lemma}
	If $E$ is a bounded set of class $\mathrm{C}^k$ ($k\geq 2$), then there exists $C>0$ such that, for every $\varphi$ in~$\mathrm{C}^{k-1}(\partial E)$, the mapping $\varphi\,\nu_E$ can be extended to $\xi_\varphi\in\mathrm{C}^{k-1}(\RR^2,\RR^2)$ with
	$$\|\xi_\varphi\|_{\mathrm{C}^{k-1}(\RR^2,\RR^2)}\leq C\,\|\varphi\|_{\mathrm{C}^{k-1}(\partial E)}\,.$$
	\label{lemma_extension}
\end{lemma}

\begin{proposition}
	Let $E$ be a bounded open set of class $\mathrm{C}^k$ (with $k\geq 2$). There exists $c>0$ such that, for every~${\varphi\in\mathrm{C}^{k-1}(\partial E)}$ with $\|\varphi\|_{\mathrm{C}^{k-1}(\partial E)}\leq c$, there is a unique bounded open set of class~$\mathrm{C}^{k-1}$, denoted $E_\varphi$, satisfying
	\begin{equation}
		\partial E_\varphi=(Id+\varphi\,\nu_E)(\partial E)\,.
		\label{normal_deform_boundary_eq}
	\end{equation}
	Moreover, there exists an extension $\xi_\varphi$ of $\varphi\,\nu_E$ such that $E_\varphi=(Id+\xi_\varphi)(E)$ and $$\|\xi_\varphi\|_{\mathrm{C}^{k-1}(\RR^2,\RR^2)}<1\,.$$
	In particular, $E_\varphi$ is $\mathrm{C}^{k-1}$-diffeomorphic to $E$.
	\label{normal_deform_boundary}
\end{proposition}

\begin{proposition}
	If $(E_n)_{n\geq 0}$ converges to a bounded set $E$ in $\mathrm{C}^k$ with $k\geq 2$, then for $n$ large enough there exists $\varphi_n\in\mathrm{C}^{k-1}(\partial E)$ such that $E_n=E_{\varphi_n}$, and $\|\varphi_n\|_{\mathrm{C}^{k-1}(\partial E)}\to 0$.
	\label{prop_conv_normal}
\end{proposition}

\subsection{Functions of bounded variation and sets of finite perimeter}
\label{prelim_bv}
We recall here a few properties of functions of bounded variation and sets of finite perimeter. More details can be found in the monogaphs~\cite{ambrosioFunctionsBoundedVariation2000,maggiSetsFinitePerimeter2012}.
\paragraph{The total variation.} The total variation of a function $u\in\mathrm{L}^1_{\mathrm{loc}}(\RR^2)$ is given by
\begin{equation*}
\TV(u)\eqdef\mathrm{sup}\,\bigg\{-\int_{\RR^2}u\,\mathrm{div}\,z\,\bigg\rvert\,z\in\mathrm{C}^{\infty}_c(\RR^2,\RR^2),~\|z\|_{\infty}\leq 1\bigg\}\,.
\end{equation*}
If $\TV(u)$ is finite, then $u$ is said to have bounded variation, and its distributional gradient $\Diff u$ is a finite Radon measure. In that case, we have $|\Diff u|(\RR^2)=\TV(u)$. In all the following, we consider $\TV$ as a mapping from $\LD$ to $\RR\cup\{+\infty\}$. This mapping is convex, proper and lower semi-continuous.

\paragraph{Sets of finite perimeter.} If a measurable set $E\subset\RR^2$ is such that $P(E)\eqdef \TV(\mathbf{1}_E)<+\infty$, it is said to be of finite perimeter. If $E$ is an open set of class $\mathrm{C}^1$, then $P(E)$ is simply the length of its topological boundary, $P(E)= \Hh^1(\partial E)$, where $\Hh^1$ denotes the one-dimensional Hausdorff measure.

\paragraph{Coarea formula}
Functions with bounded variation and sets of finite perimeter are related through the coarea formula~\cite[Thm. 3.40]{ambrosioFunctionsBoundedVariation2000}. For $u \in \mathrm{L}^1_{\mathrm{loc}}(\RR^2)$ and~${t\in\RR}$, we consider the level sets of $u$,
\begin{equation}\label{eq_levelsets}
	U^{(t)}\eqdef
	\begin{cases}
		\{x\in\RR^2\,\rvert\,u(x)\geq t\} & \text{if } t\geq 0\,,\\
		\{x\in\RR^2\,\rvert\,u(x)\leq t\} & \text{otherwise.}
	\end{cases}
\end{equation}
It is worth noting that, if $u\in\LD$, then $|U^{(t)}|<+\infty$ for all $t\neq 0$.
The coarea formula states that
\begin{align}\label{eq_coarea}
	\forall u \in \LD,\quad	\TV(u)&= \int_{-\infty}^{+\infty} P(U^{(t)}) \diff  t.
\end{align}

\paragraph{The isoperimetric inequality.} For every set of finite perimeter $E$, the isoperimetric inequality states that
\begin{equation}
	\sqrt{\mathrm{min}(|E|,|E^c|)}\leq c_2\,P(E)\,,
	\label{isoperimetric_ineq_geom}
\end{equation}
with equality if and only if $E$ is a ball, and where $c_2\eqdef 1/\sqrt{4\pi}$ is the isoperimetric constant~(see e.g. \cite[Chapter 14]{maggiSetsFinitePerimeter2012}). In particular, if $E$ is a set of finite perimeter, either $E$ or~$E^c$ has finite measure. As a consequence of~\eqref{isoperimetric_ineq_geom} and the coarea formula, the following Poincar\'e-type inequality holds (see~\cite[Theorem 3.47]{ambrosioFunctionsBoundedVariation2000}),
\begin{align}\label{poincare_ineq}
	\forall u \in \LD,\quad \norm{u}_{2} \leq c_2\,\TV(u).
\end{align}

\paragraph{Indecomposable and simple sets.} A set of finite perimeter $E\subset \RR^2$ is said to be decomposable if there exists a partition of $E$ in two sets of positive Lebesgue measure $A$ and $B$ with~${P(E)=P(A)+P(B)}$. We say that~$E$ is indecomposable if it is not decomposable. We say that a measurable set $E$  is simple if~$E=\RR^2$, or $|E|<+\infty$ and both $E$ and $\RR^2\setminus E$ are indecomposable. The importance of simple sets stems from their connection with the extreme points of the total variation unit ball.
\begin{proposition}[{\cite{flemingFunctionsGeneralizedGradient1957,ambrosioConnectedComponentsSets2001}}]
	The extreme points of the convex set
	$$\{\TV\leq 1\}\eqdef \left\{u\in\LD\,\big\rvert\,\TV(u)\leq 1\right\}$$
	are the functions of the form $\pm\mathbf{1}_E/P(E)$, where $E$ is a simple set with $0<|E|<+\infty$.
	\label{extr_points_tv_ball}
\end{proposition}
However, it is worth noting that the \emph{exposed} points of $\{\TV\leq 1\}$ in the $\mathrm{L}^2$ topology (i.e. points that are the only maximizer over $\{\TV\leq 1\}$ of a continuous linear form on $\LD$) have much more structure (see~\Cref{faces_tv,prelim_pc}).

\subsection{Subdifferential of the total variation}
Let us now collect several results on the subdifferential of $\TV$, which are useful to derive and analyze the dual problems of \cref{primal_noreg} and \cref{primal_reg}. Since $\TV\colon\LD\to\RR\cup\{+\infty\}$ is the support function of the convex set
$$C\eqdef\left\{\mathrm{div}\,z\,\big\rvert\,z\in\mathrm{C}_c^{\infty}(\RR^2,\RR^2),~\|z\|_{\infty}\leq 1\right\},$$
its subdifferential at $0$ is the closure of $C$ in~$\LD$, that is
\begin{equation}
	\partial\TV(0)=\overline{C}=\left\{\text{div}\,z\,\big\rvert\,z\in
	\mathrm{L}^{\infty}(\mathbb{R}^2,\mathbb{R}^2),~\text{div} z\in
	\LD,~||z||_{\infty}\leq 1\right\}.
	\label{charac_subdiff_z}
\end{equation}
We also have the following useful identity:
\begin{equation}\label{eq_subdiff_zero}
	\begin{aligned}
		\partial \TV(0)= \left\{\eta\in \LD\,\bigg\rvert\,\forall u\in \LD,\,\left|\int_{\mathbb{R}^2}\eta\,u\right|\leq \TV(u)\right\}.
	\end{aligned}
\end{equation}
Finally, the sudifferential of $\TV$ at some $u\in\LD$ is given by:
\begin{equation}\label{eq_subdiff_qcq}
	\partial \TV(u)=\left\{\eta\in\partial \TV(0)\,\bigg\rvert\,\int_{\RR^2}\eta\,u=\TV(u)\right\}\,.
\end{equation}
Hence, if $\eta\in\partial\TV(u)$, then $\eta$ is an element of $\partial\TV(0)=\overline{C}$ for which the supremum in the definition of the total variation is attained.

\subsection{Dual problems and dual certificates}
\label{sec_dual}
The backbone of our main result is the relation between the solutions of~\cref{primal_noreg} or~\cref{primal_reg} and the solutions of their dual problems. We gather here several properties of these dual problems which can be found in~\cite{chambolleGeometricPropertiesSolutions2016} (for the denoising case) and~\cite[Section 2]{iglesiasNoteConvergenceSolutions2018} (for the general case).

\paragraph{Dual problems.} The Fenchel-Rockafellar dual problems to~\cref{primal_noreg} and \cref{primal_reg} are respectively
\begin{equation}
	\underset{p \in \mathcal{H}}{\text{sup}}~  \langle p,y_0\rangle_{\mathcal{H}} ~ \text{ s.t. }\Phi^*p\in\partial \TV(0)\,,
	\tag{$\mathcal{D}_0(y_0)$}
	\label{dual_noreg}
\end{equation}

\begin{equation}
	\underset{p \in \mathcal{H}}{\text{sup}}~ \langle p,y\rangle_{\mathcal{H}} - \frac{\lambda}{2}||p||_{\mathcal{H}}^2
	~ \text{ s.t. } \Phi^*p\in\partial \TV(0)\,.
	\tag{$\mathcal{D}_{\lambda}(y)$}
	\label{dual_reg}
\end{equation}

The existence of a solution to~\cref{dual_noreg} does not always hold.  On the contrary, \cref{dual_reg} can be reformulated as the problem of projecting $y/\lambda$ onto the closed convex set~${\{p\in\mathcal{H}\,\rvert\,\Phi^*p\in\partial \TV(0)\}}$, which has a unique solution.

\paragraph{Strong duality.}
From \cite[Theorem 1]{iglesiasNoteConvergenceSolutions2018}, the values of \cref{primal_noreg} and \cref{dual_noreg} are equal. Moreover, if there exists a solution $p$ to~\cref{dual_noreg}, then for every solution $u$ of \cref{primal_noreg} we have
\begin{equation}\Phi^*p\in\partial \TV(u)\,.\label{opt_noreg}\end{equation}
Conversely, if $(u,p)\in \LD\times \mathcal{H}$ with $\Phi u=y_0$ and \cref{opt_noreg} holds, then $u$ and $p$ respectively solve~\cref{primal_noreg} and \cref{dual_noreg}.
From the perspective of inverse problems, given some unknown image $u_0 \in\LD$ and observation $y_0=\Phi u_0$, it is therefore sufficient to assume the existence of~$p \in \mathcal{H}$ with~${\Phi^*p \in \partial \TV(u_0)}$ to ensure that $u_0$ is a solution to~\cref{primal_noreg}. This property is known as the \emph{source condition}~\cite{neubauerTikhonovRegularisationNonlinear1989,burgerConvergenceRatesConvex2004}. If, moreover, $\Phi$ is injective on the cone~$\{u\in\LD\,\rvert\,\Phi^*p \in \partial \TV(u)\}$, then $u_0$ is the unique solution to ~\cref{primal_noreg}.

As in the noiseless case, the values of \cref{primal_reg} and \cref{dual_reg} are equal. Moreover, denoting  by~$p$ the unique solution to \cref{dual_reg}, for every solution $u$ of \cref{primal_reg} we have
		\begin{equation}
			\left\{\begin{aligned}
				&\Phi u=y-\lambda\,p\,,\\
				&\Phi^*p\in\partial \TV(u)\,.
			\end{aligned}\right.
			\label{opt_reg}
		\end{equation}
		Conversely, if \cref{opt_reg} holds, then $u$ and $p$ respectively solve \cref{primal_reg} and \cref{dual_reg}. Although there might not be a unique solution to \cref{primal_reg}, \Cref{opt_reg} yields that all of them have the same image by~$\Phi$ and the same total variation.

\paragraph{Dual certificates.}If $\eta=\Phi^*p$ and $\eta\in\partial\TV(u)$, we call $\eta$ a dual certificate for $u$ with respect to~\cref{primal_noreg}, as its existence certifies the optimality of $u$ for \cref{primal_noreg}, provided  $y_0=\Phi u$. Similarly, if~${\eta=-\Phi^*(\Phi u-y)/\lambda}$ and $\eta\in\partial\TV(u)$, we call $\eta$ a dual certificate for $u$ with respect to
\cref{primal_reg}. There could be multiple dual certificates associated to \cref{primal_noreg}. One of them, the minimal norm certificate, plays a crucial role in the analysis of the low noise regime. A quick look at the objective of~\cref{dual_reg} indeed suggests that, as $\lambda$ goes to $0$, its solution converges to the solution to the limit problem~\cref{dual_noreg} with minimal norm. This is \Cref{conv_dual_vec} below.
	\begin{definition}
		If there exists a solution to \cref{dual_noreg}, the minimal norm dual certificate associated to \cref{primal_noreg}, denoted~$\eta_0$, is defined as
		\begin{equation*}
			\eta_0=\Phi^*p_0 ~\text{ with }~p_0=\mathrm{argmin}~\|p\|_{\mathcal{H}}~\text{ s.t. }~ p \text{ solves } \cref{dual_noreg}\,.
		\end{equation*}
		\label{min_norm_certif}
	\end{definition}
If $\lambda>0$, we denote $p_{\lambda,w}$ the unique solution to \customref{dual_reg}{$\mathcal{D}_\lambda(y_0+w)$}, and~${\eta_{\lambda,w}=\Phi^*p_{\lambda,w}}$ the associated dual certificate. Noise robustness results extensively rely on the behaviour of $\eta_{\lambda,w}$ as $\lambda$ and $w$ go to zero. This behaviour is described by the following results.
\begin{proposition}[{\protect\cite[Prop.6]{chambolleGeometricPropertiesSolutions2016},\cite[Prop. 3]{iglesiasNoteConvergenceSolutions2018}}]
	If there exists a solution to \cref{dual_noreg}, then $p_{\lambda,0}$ converges strongly to $p_0$ as~${\lambda\to 0}$.
	\label{conv_dual_vec}
\end{proposition}

Since $p_{\lambda,w}$ is the projection of $(y_0+w)/\lambda$ onto the closed convex set $\{p\in\mathcal{H}\,\rvert\,\Phi^*p\in\partial \TV(0)\}$, the non-expansiveness of the projection mapping yields
\begin{equation}
	\forall (\lambda,w)\in \RR_+^*\times \mathcal{H},~\|p_{\lambda,w}-p_{\lambda,0}\|_{\mathcal{H}}\leq \frac{\|w\|_{\mathcal{H}}}{\lambda}\,,
	\label{proj_dual_vec}
\end{equation}
and hence
\begin{equation*}
	\forall (\lambda,w)\in \RR_+^*\times \mathcal{H},~\|\eta_{\lambda,w}-\eta_{\lambda,0}\|_{\LD}\leq \frac{\|\Phi^*\|\,\|w\|_{\mathcal{H}}}{\lambda}\,.
\end{equation*}
As a result, if $\lambda\to 0$ and $\|w\|_{\mathcal{H}}/\lambda\to 0$, the dual certificate $\eta_{\lambda,w}$ converges strongly in $\LD$ to the minimal norm certificate $\eta_0$.

\subsection{Noise robustness results}
Let us now review existing noise robustness results, which we use in various parts of this work. From \Cref{subdiff_lvlsets} and the results of \Cref{sec_dual}, we know that the levels sets of solutions to~\customref{primal_reg}{$\mathcal{P}_{\lambda}(y_0+w)$} are solution to the prescribed curvature problem associated to $\eta_{\lambda,w}$. In~\cite{chambolleGeometricPropertiesSolutions2016,iglesiasNoteConvergenceSolutions2018}, this fact is exploited to obtain uniform properties of the level sets in the low noise regime. We collect the byproducts of this analysis in the following lemma.
	\begin{lemma}[{\cite[Section 5]{chambolleGeometricPropertiesSolutions2016}}]
	Let $(\eta_n)_{n\geq 0}\subset \partial\TV(0)$ be a sequence converging strongly in $\LD$ to $\eta_{\infty}$, and let $\mathcal{E}$ be defined by
	\begin{equation*}
		\mathcal{E}\eqdef\left\{E\subset\RR^2,~0<|E|<+\infty\,\bigg\rvert\,\exists n\in \NN\cup\{\infty\},~P(E)=\left|\int_{E}\eta_n\right|\right\}.
	\end{equation*}
	Then the following holds:
	\begin{enumerate}
		\item $\underset{E\in\mathcal{E}}{\mathrm{inf}}~P(E)>0$ and $\underset{E\in\mathcal{E}}{\mathrm{sup}}~P(E)<+\infty$,
		\item $\underset{E\in\mathcal{E}}{\mathrm{inf}}~|E|>0$ and $\underset{E\in\mathcal{E}}{\mathrm{sup}}~|E|<+\infty$,
		\item there exists $R>0$ such that, for every $E\in\mathcal{E}$, it holds $E\subset B(0,R)$,
		\item there exists $r_0>0$ and $C\in(0,1/2)$ such that for every $r\in (0,r_0]$ and $E\in\mathcal{E}$:
		\begin{equation*}
			\forall x\in\partial E,~C\leq\frac{|E\cap B(x,r)|}{|B(x,r)|}\leq 1-C\,.
		\end{equation*}
	\end{enumerate}
	\label{lemma_bounds}
\end{lemma}
In the above-mentioned works, \Cref{lemma_bounds} is used to obtain the convergence result of \Cref{conv_strict_bv}. It indeed allows to show that, in the low noise regime, the solutions to \customref{primal_reg}{$\mathcal{P}_{\lambda}(y_0+w)$} have bounded support, and therefore belong to $\mathrm{L}^1(\RR^2)$. This can in turn be used to show their strict convergence in $\mathrm{BV}(\RR^2)$ towards a solution $u_*$ of \Cref{primal_noreg}, which in particular imply the weak-* convergence of their gradient (see e.g. \cite[Proposition 3.13 and Definition~3.14]{ambrosioFunctionsBoundedVariation2000}). Finally, one obtains the convergence of their level set towards those of~$u_*$ in the Hausdorff sense, which corresponds to the uniform convergence of the associated distance functions (see Theorem 6.1 and its proof in \cite{ambrosioFunctionsBoundedVariation2000}).
\begin{proposition}
	Assume \cref{dual_noreg} has a solution, $\lambda_n\to 0$ and
	$$\frac{\|w_n\|_{\mathcal{H}}}{\lambda_n}\leq \frac{1}{4c_2\,\|\Phi^*\|}\,.$$
	Then, if $u_n$ is a solution of \customref{primal_reg}{$\mathcal{P}_{\lambda_n}(y_0+w_n)$} for all $n\in\NN$, we have that $(\mathrm{Supp}(u_n))_{n\geq 0}$ is bounded and that, up to the extraction of a subsequence (not relabeled), $(u_n)_{n\geq 0}$ converges strictly in~$\mathrm{BV}(\RR^2)$ to a solution $u_*$ of \cref{primal_noreg}. Moreover, for almost every $t\in\RR$, we have:
	$$\big|U_n^{(t)}\triangle U_*^{(t)}\big|\longrightarrow 0~~~\mathrm{and}~~~\partial U_n^{(t)}\longrightarrow\partial U_{*}^{(t)}\,,$$
	where the last limit holds in the Hausdorff sense\footnote{See e.g. \cite[Chapter 4]{rockafellarVariationalAnalysis1998} for a definition.}.
	\label{conv_strict_bv}
\end{proposition}

%% file: faces_tv.tex

\section{The exposed faces of the total variation unit ball}
\label{faces_tv}
We recall that, in the remaining of this work, we assume that
$\phi\in\mathrm{C}^1(\RR^2,\mathcal{H})$, so that $\Phi^*$ is continuous from
$\mathcal{H}$ to $\mathrm{C}^1(\RR^2)$.

In order to take advantage of the extremality relations~\eqref{opt_noreg}
and~\eqref{opt_reg}, it is important to understand the properties of $u$ implied
by the relation $\eta \in \partial \TV(u)$, for a given $\eta \in \partial
\TV(0)$. In other words, our goal is to study the set
\begin{align}\label{argminsousdiff}
	\partial \TV^*(\eta)= \left\{u \in \LD\,\big\rvert\,\eta \in \partial \TV(u)\right\}=\uArgmax{u \in \LD} \left(\int_{\RR^2} \eta u - \TV(u)\right),
\end{align}
where $\TV^*$ denotes the Fenchel conjugate of $\TV$.

\subsection{Subgradients and exposed faces}
\label{subgrad_exposed}
It is possible to relate~\eqref{argminsousdiff} to the faces of the total
variation unit ball, in connection with Fleming's result (i.e.
\Cref{extr_points_tv_ball}). We say that a set $\mathcal{F}\subseteq \{\TV\leq
1\}$ is an exposed face of~${\{\TV\leq 1\}}$ if there exists $\eta \in \LD$ such
that
\begin{align}
	\mathcal{F} =\uArgmax{u \in \{\TV\leq 1\}} \int_{\RR^2}\eta u\,.
\end{align}
Exposed faces are closed convex subsets of $\{\TV\leq 1\}$, and they are faces
in the classical sense: if~$u \in \mathcal{F}$ and $I\subset \{\TV\leq 1\}$ is
an open line segment containing $u$, then $I \subset \mathcal{F}$. We refer the
reader to~\cite[Chapter 18]{rockafellarConvexAnalysis1970} for more detail on
faces and exposed faces. To emphasize the dependency on $\eta$ we sometimes
write $\mathcal{F}_\eta$ for $\mathcal{F}$, and we note that
\begin{align*}
	\Ff_0=\{\TV\leq 1\} \quad\mbox{and}\quad \Ff_{t\eta}=\Ff_{\eta} \ \mbox{for all $t>0$}.
\end{align*}

It is also worth considering the corresponding value, which is sometimes called
the $G$-norm\footnote{The $G$-norm is the polar of the total variation
(see~\cite[Ch. 15]{rockafellarConvexAnalysis1970}). It is possible to prove that
the~$G$-norm is indeed a norm on $\LD$, in particular $0<\norm{\eta}_G<+\infty$
for all $\eta \in \LD\setminus \{0\}$, see~\cite{haddad_texture_2007}.} in the
literature~\cite{meyer_oscillating_2001,aujol_image_2005,kindermannDenoisingBVduality2006,haddad_texture_2007},
\begin{align*}
	\forall \eta \in \LD,\quad \norm{\eta}_G \eqdef \sup_{u \in \{\TV\leq 1\}} \int_{\RR^2} \eta u\,.
\end{align*}
In view of~\eqref{eq_subdiff_zero}, we see that $\eta \in \partial \TV(0)$ if
and only if $\norm{\eta}_G\leq 1$. Assuming that $\norm{\eta}_G\leq 1$, the
condition $\int_{\RR^2}\eta u= \TV(u)$ in~\eqref{eq_subdiff_qcq} is equivalent
to $u=0$ or
\begin{align}
	\int_{\RR^2}\eta\left(\frac{u}{\TV(u)}\right)=1,
\end{align}
the latter equality implying that $\norm{\eta}_G=1$ and $u/\TV(u) \in
\Ff_{\eta}$. As a result, we obtain the following description,
\begin{align}
	\TV^*(\eta) &=
	\begin{cases}
		\emptyset & \mbox{if $\norm{\eta}_G>1$},\\
		\{0\} \cup \left(\bigcup_{t>0} (t\mathcal{F}_\eta) \right),&\mbox{if $\norm{\eta}_G=1$,} \\
		\{0\} & \mbox{if $\norm{\eta}_G<1$}.
	\end{cases}
	  \label{subdiff_face_exposee}
\end{align}
To summarize the connection between subgradients and exposed faces, if $\Ff$ is
a face of $\{\TV\leq 1\}$ exposed by some $\eta \in \LD\setminus \{0\}$, its
conic hull $\RR_+\Ff$ is equal to $\partial \TV^*(\eta/\norm{\eta}_G)$.
Conversely, if~$\norm{\eta}_G=1$, then $\partial \TV^*(\eta)\cap \{\TV=1\}$ is
$\Ff_\eta$, the face of $\{\TV\leq 1\}$ exposed by $\eta$.

In the rest of this section, we fix some  $\eta \in \LD$ such that
$\norm{\eta}_G=1$ (the only interesting case), and we study $\TV^*(\eta)$.
Equivalently, we describe all the faces of $\{\TV\leq 1\}$ exposed by nonzero
vectors.

\subsection{Exposed versus non-exposed faces}
The extreme points of $\{\TV\leq 1\}$ are described
by~\Cref{extr_points_tv_ball}: those are the (signed, renormalized) indicators
of simple sets. The $k$-dimensional faces ($k \in \NN$) are more complex, and
they involve functions which are piecewise constant on some partition of $\RR^2$
(see for instance~\cite{duvalFacesExtremePoints2022}, or  the
monographs~\cite{fujishigeSubmodularFunctionsOptimization2005,bachLearningSubmodularFunctions2013}
in a finite-dimensional setting). Even though it is known that the
$k$-dimensional faces have a finite number of extreme points and are thus
polytopes (see \cite[Theorem 2.1]{duvalFacesExtremePoints2022}), the
corresponding partition can be singular and
counter-intuitive~\cite{boyerinpreparation}.

However, the faces involved in~\eqref{subdiff_face_exposee} are the
\emph{exposed} faces, and we emphasize in this section that they have a simpler
structure than arbitray faces, especially if~$\eta \in \Cder^1$, as is the case
in the extremality conditions~\eqref{opt_noreg} and~\eqref{opt_reg}. We prove
in~\Cref{number_extr_points} that, under this assumption, the~$k$-dimensional
exposed faces of $\{\TV\leq 1\}$ are $k$-simplices.

At the core of our discussion is the reformulation of the subdifferential
property into a geometric variational problem using level sets and the coarea
formula (see~\eqref{eq_levelsets} and~\eqref{eq_coarea}).
	\begin{proposition}[\cite{kindermannDenoisingBVduality2006,chambolleGeometricPropertiesSolutions2016}]
		\label{subdiff_lvlsets}
		Let $u\in\LD$ be such that~${\TV(u)<+\infty}$, and let~${\eta\in\LD}$.
		Then the following conditions are equivalent.
		\begin{enumerate}[label=(\roman*)]
			\item $\eta\in\partial \TV(u)\,.$
			\item $\eta\in\partial \TV(0)$ and the level sets of $u$ satisfy
			\begin{equation*}
				\forall t>0,~P(U^{(t)})=\int_{U^{(t)}}\eta~~~\text{and }~~
				\forall t<0,~P(U^{(t)})=-\int_{U^{(t)}}\eta\,.
			\end{equation*}
			\item The level sets of $u$ satisfy
			\begin{equation*}
				\begin{aligned}
					&\forall t>0,~U^{(t)}\in\underset{E\subset\RR^2,\,|E|<+\infty}{\mathrm{Argmin}}\left(P(E)-\int_{E}\eta\right),\\
					&\forall t<0,~U^{(t)}\in\underset{E\subset\RR^2,\,|E|<+\infty}{\mathrm{Argmin}}\left(P(E)+\int_{E}\eta\right).
				\end{aligned}
			\end{equation*}
		\end{enumerate}
	\end{proposition}

\subsection{The prescribed curvature problem}
\label{prelim_pc}
The geometric variational problem appearing in~\Cref{subdiff_lvlsets},
\begin{equation}
\underset{E\subset\RR^2,\,|E|<+\infty}{\mathrm{inf}}~J(E)\eqdef P(E)-\int_{E}\eta\,,
\tag{$\mathcal{PC}(\eta)$}
\label{prescribed_curv}
\end{equation}
is called the \emph{prescribed curvature problem} associated to $\eta$. This
terminology stems from the fact that, if $\eta$ is sufficiently regular, every
solution to \cref{prescribed_curv} has a (scalar) distributional curvature~(see
\cite[Section 17.3]{maggiSetsFinitePerimeter2012} for a definition) equal to
$\eta$. This problem plays a crucial role in the analysis of total variation
regularization, as explained below. For now, let us gather some properties of
that problem.

\paragraph{Existence of minimizers.} Solutions to~\cref{prescribed_curv} exist
provided $\eta \in \partial \TV(0)$. Indeed, the objective $J$ is  nonnegative,
and equal to zero for $E=\emptyset$. From \Cref{subdiff_lvlsets}, we also know
there is a non-empty solution as soon as ${\eta\in\partial\TV(u)}$ for some
$u\in\LD\setminus \{0\}$.

\paragraph{Boundedness.} By \cite[Lemma
4]{chambolleGeometricPropertiesSolutions2016}, all solutions of
\cref{prescribed_curv} are included in some common ball, i.e. there exists $R>0$
such that, for every solution $E$ of~\cref{prescribed_curv}, we have~${E\subset
B(0,R)}$.

\paragraph{Regularity of the solutions.} The regularity of the solutions
to~\cref{prescribed_curv} is well understood. If~$\eta$ is only assumed to be
square integrable, the solutions can be singular but they  have some weak form
of regularity, as shown in \cite{gonzalesBoundariesPrescribedMean1993}. In
particular, it is known that the square $C=[0,1]^2$ is not a solution
to~\cref{prescribed_curv} for any $\eta \in \LD$ (see, e.g.\
\cite{meyer_oscillating_2001}). As a result, \emph{the function $\bun_{C}/P(C)$
is an extreme point of $\{\TV\leq 1\}$ which is not exposed.} More regularity
can be obtained by strengthening the integrability and smoothness of $\eta$. If,
in addition to being square integrable,
${\eta\in\mathrm{L}^{\infty}_{\mathrm{loc}}(\RR^2)}$ (which ensures
that~${\eta\in\mathrm{L}^{\infty}(B(0,R))}$), then any solution to
\eqref{prescribed_curv} is a strong quasi-minimizer of the perimeter, and,
consequently, is equivalent to an open set of class~$\mathrm{C}^{1,1}$~(see e.g.
\cite[Definition 4.7.3 and Theorem~4.7.4]{ambrosioCorsoIntroduttivoAlla2010}).
Furthermore, if $\eta$ is continuous, then the boundary of any solution is
locally the graph of a function $u$ which solves (in the sense of distributions)
the Euler-Lagrange equation associated to~\eqref{prescribed_curv}, that is (up
to a translation and a rotation):
\begin{equation}
\left(\frac{u'}{\sqrt{1+{u'}^2}}\right)'\-(z)=\eta(z, u(z))\,.
\label{eulerlag_prescribed_curv}
\end{equation}
This in turn implies that $u$ is $\mathrm{C}^2$ ($\mathrm{C}^{k+2,\alpha}$ if
$\eta\in\mathrm{C}^{k,\alpha}(\RR^2)$) and solves
\cref{eulerlag_prescribed_curv} in the classical sense.

\subsection{Indicator functions corresponding to a given face}
We fix $\eta \in \LD$ such that $\norm{\eta}_G=1$ (hence $\eta \in \partial
\TV(0)$), and we study the face  $\mathcal{F}$ of $\{\TV\leq 1\}$ exposed by
$\eta$. We assume in addition that $\eta \in \Cder^1(\RR^2)$.

As we know that the extreme points of $\mathcal{F}$ must be (signed,
renormalized) indicators of simple sets, it is natural to focus on such
functions. The main result we prove in this section is the following.
\begin{restatable}{proposition}{extrpointsdisjsupp}
	For any extreme point $u$ of $\mathcal{F}$, there exists a unique pair
	$(s,E)$, where $E$ is a simply connected open set of class
	$\Cder^3$ and ${s\in \{-1,1\}}$, such that $u=s
	\bun_{E}/P(E)$.	

	If $u_1$ and $u_2$ are two distinct extreme points of $\mathcal{F}$, and
	$\{(s_i,E_i)\}_{i=1,2}$ are their corresponding decompositions,
	then $\partial E_1\cap\partial E_2=\emptyset$.
	\label{extr_points_disj_supp}
	\label{EXTR_POINTS_DISJ_SUPP}
\end{restatable}

To obtain this result, we study the elements of $\mathcal{F}$ which are
(proportional to) indicator functions, and we introduce the collection
\begin{equation}
	\begin{aligned}
		\mathcal{E}&\eqdef\mathcal{E}^+\cup\mathcal{E}^-\cup\{\emptyset,\RR^2\}\,,~\mbox{where}\\
		\mathcal{E}^+&\eqdef\left\{E\subset\RR^2\,\bigg\rvert\,|E|<+\infty,~0<P(E)<+\infty,~\frac{\mathbf{1}_E}{P(E)}\in\mathcal{F}\right\},\\
		\mathcal{E}^-&\eqdef\left\{E\subset\RR^2\,\bigg\rvert\,|E^c|<+\infty,~0<P(E^c)<+\infty,~\frac{-\mathbf{1}_{E^c}}{P(E^c)}\in\mathcal{F}\right\}.
	\end{aligned}
	\label{def_E}
\end{equation}
If $\abs{E}<+\infty$ (resp. $\abs{E}=+\infty$), \Cref{subdiff_lvlsets} above
shows that $E \in \mathcal{E}^+$ (resp. $E \in \mathcal{E}^-$) if and only if
$E$ is a solution to \eqref{prescribed_curv} (resp. $E^c$ is a solution to
$\mathcal{PC}(-\eta)$).

\subsubsection{Structure of \texorpdfstring{$\mathcal{E}$}{E}}
\label{structure_set_E}
The collection $\mathcal{E}$ has the remarkable property of being closed under
union and intersection.
\begin{proposition}\label{prop_stability_union} Let $E\in\mathcal{E}$ and
	$F\in\mathcal{E}$. Then $E\cap F\in\mathcal{E}$ and $E\cup F\in\mathcal{E}$.
	\label{prop_stab_E}
\end{proposition}
In fact, $\mathcal{E}$ is even closed under \emph{countable} union and
intersection, but we do not need this property here.
\begin{proof}
	If $E\in\mathcal{E}^+$ and $F\in\mathcal{E}^+$ the submodularity of the
	perimeter (see e.g. \cite[Proposition
	1]{ambrosioConnectedComponentsSets2001}) yields:
	$$P(E\cap F)+P(E\cup F)\leq P(E)+P(F)=\int_{E}\eta+\int_{F}\eta=\int_{E\cap
	F}\eta+\int_{E\cup F}\eta\,.$$ We hence obtain:
	$$\left(P(E\cap F)-\int_{E\cap F}\eta\right)+\left(P(E\cup F)-\int_{E\cup
	F}\eta\right)\leq 0\,.$$ By~\eqref{eq_subdiff_zero}, the above two terms are
	nonnegative, which yields~${E\cap F\in\mathcal{E}^+}$~(unless~${E\cap
	F=\emptyset}$) and~$E\cup F\in\mathcal{E}^+$. The same argument applies to
	the complements, when both $E$ and $F$ are in $\mathcal{E}^-$. Now, if $E
	\in \mathcal{E}^+$ and $F \in \mathcal{E^-}$,
   \begin{equation*}
	\begin{aligned}
		P(E\cap F)+P((E\cup F)^c)&=P(E\cap F)+P(E\cup F)\\
		&\leq P(E)+P(F)\\
		&=P(E)+P(F^c)\\
		&=\int_{E}\eta-\int_{F^c}\eta\\
		&=\int_{E\cap F}\eta-\int_{(E\cup F)^c}\eta
	\end{aligned}
\end{equation*}
Reasoning as above, we obtain that $E\cap F\in\mathcal{E}^+$ (unless $E\cap
F=\emptyset$) and $E\cup F\in\mathcal{E}^-$ (unless~${E\cup F=\RR^2}$).
\end{proof}

\subsubsection{Relative position of elements of \texorpdfstring{$\mathcal{E}$}{E}}
\label{relative_pos}

In view of the regularity results of~\Cref{prelim_pc} and the assumption that
$\eta \in \Cder^1$,  the solutions of the prescribed curvature problem
associated to $\pm \eta$ (and hence the elements of $\mathcal{E}$) are
equivalent to open sets of class $\mathrm{C}^3$. This property, together with
\Cref{prop_stability_union}, imposes strong constraints on the intersection of
the boundaries of elements of $\mathcal{E}$, as the next proposition shows.

\begin{proposition}
	Let $E,F\in\mathcal{E}$. Then
	\begin{equation}
		\partial E\cap\partial F=\{\nu_E=\nu_F\}\cup\{\nu_E=-\nu_F\}\,.
		\label{decomp_inter_boundary}
	\end{equation}
	Moreover, the sets $\{\nu_E=-\nu_F\}$ and $\{\nu_E=\nu_F\}$ are both open
	and closed in $\partial E$ and $\partial F$.
	\label{boundary_elts_E}
\end{proposition}
Let us recall that the topological boundaries and the normals mentioned above
are those of the unique open representative of $E$ (resp. $F$), see
\Cref{prelim_smooth}. The proof of \Cref{boundary_elts_E} follows from the next
two Lemmas.
\begin{lemma}
	Let $E$ and $F$ be two sets of class $\mathrm{C}^1$. If $E\cap F$ and $E\cup
	F$ (modulo a Lebesgue-negligible set) are trivial or $\mathrm{C}^1$, then
	$$\partial E\cap\partial F=\{\nu_E=\nu_F\}\cup\{\nu_E=-\nu_F\}\,$$
	and~${\{\nu_E=-\nu_F\}}$ is both open and closed in $\partial E$ and
	$\partial F$
	\label{first_lemma}
\end{lemma}
\begin{proof}
	The regularity of $E,~F,~E\cap F$ and $E\cup F$ implies that the densities
	$\theta_E,~\theta_F,~\theta_{E\cap F}$ and~$\theta_{E\cup F}$~(defined in
	\Cref{prelim_smooth}) are well-defined on $\RR^2$ and take values in
	$\{0,1/2,1\}$. Moreover, since
	$$|E\cap B(x,r)|+|F\cap B(x,r)|=|(E\cap F)\cap B(x,r)|+|(E\cup F)\cap
	B(x,r)|\,,$$ we have:
	$$\theta_E+\theta_F=\theta_{E\cap F}+\theta_{E\cup F}\,.$$ Since $E$ and $F$
	are $\mathrm{C}^1$, for every $x\in\partial E\cap\partial F$ we have
	$\theta_E(x)=\theta_F(x)=1/2$, which yields $$\theta_{E\cap
	F}(x)+\theta_{E\cup F}(x)=1\,.$$ Since $\theta_{E\cap F}\leq \theta_{E\cup
	F}$, we obtain $(\theta_{E\cap F}(x),\theta_{E\cup F}(x))=(0,1)$ or
	$(\theta_{E\cap F}(x),\theta_{E\cup F}(x))=(1/2,1/2)$.

Now, by a blow-up argument as in the proof of \cite[Theorem
16.3]{maggiSetsFinitePerimeter2012}, we note that:
	$$\frac{(E\cap F)-x}{r}\cap B(0,1)=\left[\frac{E-x}{r}\cap B(0,1)\right]\cap
	\left[\frac{F-x}{r}\cap B(0,1)\right]\xrightarrow{r\to
	0^+}B_{\nu_E(x)}^-\cap B_{\nu_F(x)}^-\,,$$ with $B_\nu^-\eqdef\enscond{x\in
	B(0,1)}{\langle x,\nu\rangle\leq 0}$, and where the convergence is in measure.
	Since, for any measurable set $A$,  $$ \left|\frac{A-x}{r}\cap
	B(0,1)\right|=\frac{|A\cap B(x,r)|}{r^d}\,,$$ we deduce that, if
	$\theta_{E\cap F}(x)=0$, then $|B^-_{\nu_E(x)}\cap B^-_{\nu_F(x)}|=0$, hence
	$\nu_E(x)=-\nu_F(x)$. Similarly, if~$\theta_{E\cap F}(x)=1/2$, then
	$|B^-_{\nu_E(x)}\cap B^-_{\nu_F(x)}|=\pi/2$ and $\nu_E(x)=\nu_F(x)$.

Let us now prove that $\{\nu_E=-\nu_F\}$ is both open and closed in $\partial E$
(similar arguments hold for $\partial F$). Since $\nu_E$ and~$\nu_F$ are
continuous, the set $\{\nu_E=-\nu_F\}$ is closed. Now, we show
that~$\{\nu_E=-\nu_F\}$ is open in $\partial E$. Let  $x \in \partial E$. Since
both $E$ and $F$ are of class $\Cder^1$, the above blow-up argument shows that
\begin{align*}
	x \in \{\nu_E=-\nu_F\} \Longleftrightarrow \theta_{E\cap F}(x)=0\ \mbox{and}\ \theta_{E\cup F}(x)=1.
\end{align*}
Since $E\cap F$ is (equivalent to) the empty set or an open set of class
$\Cder^1$, the set $\{\theta_{E\cap F}=0\}$ is open. Similarly, since $E\cup F$
is (equivalent to) $\RR^2$  or an open set of class $\Cder^1$, the set
$\{\theta_{E\cup F}=1\}$ is open. As a result $\{\nu_E=-\nu_F\}=(\partial E)\cap
\{\theta_{E\cap F}=0\}\cap \{\theta_{E\cup F}=1\}$ is open in $\partial E$.
\end{proof}

\begin{figure}
	\centering
	\input{regularity_inter_union}
	\caption{An example of two smooth sets $E$ (blue region) and $F$ (hatched
	region) such that~$E\cap F$ and $E\cup F$ are smooth but $\{\nu_E=\nu_F\}$
	is neither open in $\partial E$ nor in $\partial F$.}
	\label{fig_reg_union_inter}
\end{figure}
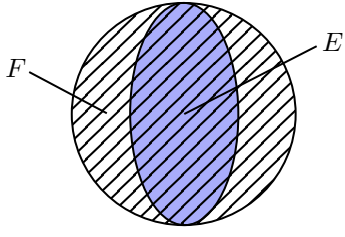

In the next lemma, we prove that $\{\nu_E=\nu_F\}$ is both open and closed in
$\partial E$ and $\partial F$. Contrary to the results of \Cref{first_lemma},
this does not hold in general if we only assume that $E$ and~$F$ are sets of
class $\mathrm{C}^1$ such that $E\cap F$ and $E\cup F$ are also of class
$\mathrm{C}^1$, as the example given in~\Cref{fig_reg_union_inter} shows.

\begin{lemma}
	If $E,F\in\mathcal{E}$ then $\{\nu_E=\nu_F\}$ is both open and closed in
	$\partial E$ and $\partial F$.
\end{lemma}
\begin{proof}
	The set $\{\nu_E=\nu_F\}$ is closed, by continuity of the normals. Let us
	prove that it is open in $\partial E$ and $\partial F$. We begin with the
	case $E, F \in \mathcal{E}^+$. Let~${x\in\{\nu_E=\nu_F\}}$, and let us
	denote~$\nu\eqdef\nu_E(x)=\nu_F(x)$. There exists~$r>0$ such that, in
	$C(x,r,\nu)$, $\partial E$ and $\partial F$ coincide with the graphs of two
	$\mathrm{C}^3$ functions which solve the prescribed curvature equation
	\begin{align}
	\frac{u''(z)}{(1+u'(z)^2)^{3/2}}=H(z,u(z))~~\mathrm{with}~~H(z,t)=\eta(x+R_\nu(z,t))
		 \label{eq_edo_prescribed}
	\end{align}
	on $(-r,r)$, with Cauchy data $u(0)=u'(0)=0$. The prescribed curvature
	equation can be reduced to a first order ODE on $\RR\times \RR^2$ defined by
	the mapping
	$$(t,x)\in \RR\times \RR^2\mapsto \begin{pmatrix}x_2\\
	H(t,x_1)(1+x_2^2)^{3/2}\end{pmatrix}$$ Since this mapping is locally
	Lipschitz continuous with respect to its second variable, the
	Cauchy-Lipschitz theorem ensures that the two functions mentioned above
	coincide on $(-r,r)$. In particular, as
\begin{align*}
	\nu_E\left(x+R_\nu(z,u(z))\right)=(u'(z),-1)/\sqrt{1+(u'(z))^2}=\nu_F\left(x+R_\nu(z,u(z))\right),
\end{align*}  the outer unit normals coincide in $C(x,r,\nu)$. As a result,
 $$\{\nu_E=\nu_F\}\cap C(x,r,\nu)=\partial E\cap C(x,r,\nu)=\partial F\cap
 C(x,r,\nu)\,.$$ which shows that $\{\nu_E=\nu_F\}$ is open in $\partial E$ and
 $\partial F$.

	If $E,F \in \mathcal{E}^-$, then $E^c$ and $F^c$ are solutions to
	$\mathcal{PC}(-\eta)$,  hence we may apply the above argument to $E^c$ and
	$F^c$, with obvious adaptations, to deduce that $\{\nu_{E^c}=\nu_{F^c}\}=
	\{\nu_{E}=\nu_{F}\}$ is open in $\partial E$ and $\partial F$.

	If $E \in \mathcal{E}^+$ and $F \in \mathcal{E}^-$, let
	$\nu\eqdef\nu_E(x)=\nu_F(x)=-\nu_{F^c}(x)$. As above, for $r>0$ small
	enough, $\partial E$ coincides in $C(x,r,\nu)$ with the graph of some
	function $u$ which satisfies~\eqref{eq_edo_prescribed},
	with~${u(0)=0}$,~$u'(0)=0$. On the other hand, $F^c$ is a solution to
	$\mathcal{PC}(-\eta)$, so that for $r>0$ small enough, it coincides in
	$C(x,r,-\nu)$ with the solution to
	\begin{align}
		\frac{v''(z)}{(1+v'(z)^2)^{3/2}}=G(z,v(z))~~\mathrm{with}~~G(z,t)=-\eta(x+R_{-\nu}(z,t)),
	\end{align}
	with $v(0)=0$, $v'(0)=0$. Since $C(x,r,-\nu)=C(x,r,\nu)$ and
	$R_{-\nu}(z,t)=R_{\nu}(z,-t)$ for all~${t \in (-r,r)}$, we observe that
	$\partial F$ coincides in $C(x,r,\nu)$ with the graph of some
	function~${\tilde{u} =-v}$ which
	satisfies~\eqref{eq_edo_prescribed} with $\tilde{u}(0)=0$,
	$\tilde{u}'(0)=0$. We conclude as before that $u$ and $\tilde{u}$ coincide
	in~$(-r,r)$, so that the set $\{\nu_E=\nu_F\}$ is open in $\partial E$ and
	$\partial F$.
\end{proof}

Now, we conclude this section by proving \Cref{extr_points_disj_supp}, whose
statement is recalled below.
\extrpointsdisjsupp*
\begin{proof}
	If $u$ is an extreme point of $\mathcal{F}$, it must be an extreme point of
	$\{\TV\leq 1\}$, hence Fleming's result (\Cref{extr_points_tv_ball}) implies
	that~${u=s \bun_{E}/P(E)}$ for some simple set $E\subset \RR^2$
	with~${0<\abs{E}<+\infty}$. Now, by~\Cref{subdiff_lvlsets}, $E$ is a
	solution to~$\mathcal{PC}\left(s \eta\right)$, so that $E$ is
	(equivalent to) an open set of class $\Cder^3$. Since $E$ is simple, that
	open set is  the interior of a rectifiable Jordan curve, as a consequence
	of~\cite[Theorem~7]{ambrosioConnectedComponentsSets2001}. Then, the
	Jordan-Schoenflies theorem implies that $E$ is homeomorphic to a
	disk, hence simply connected.

	Now, let $u_1$ and $u_2$ be two distinct extreme points. First, we note that
	$E_1\neq E_2$. Otherwise, we would have $s_1=-s_2$,
	hence $0=\frac{1}{2}s_1 \bun_{E_1}/P(E_1) +
	\frac{1}{2}s_2\bun_{E_2}P(E_2)\in \mathcal{F}$, so that
\begin{align*}
	0=\max_{u \in \{\TV\leq 1\}} \int_{\RR^2}\eta u =\norm{\eta}_G= 1,
\end{align*}
a contradiction.

Hence, $E_1\neq E_2$, and we recall that $E_i\in\mathcal{E}^+$ if
$s_i=1$ and $E_i^c\in\mathcal{E}^-$ if $s_i=-1$.
From~\Cref{boundary_elts_E}, we  know that, in any case, $\partial
E_1\cap\partial E_2$ is open and closed in $\partial E_1$ and $\partial E_2$.
Since~$\partial E_1$ and~$\partial E_2$ are Jordan curves (in particular, they
are connected) this implies $\partial E_1\cap\partial E_2=\emptyset$
or~$\partial E_1=\partial E_2$. Now, if we had $\partial E_1=\partial E_2$, the
Jordan curve theorem would yield $E_1=E_2$, which is impossible. As a result, we
obtain $\partial E_1\cap\partial E_2=\emptyset$.
\end{proof}

\subsection{Structure of finite-dimensional exposed faces}
As a consequence of \Cref{extr_points_disj_supp}, we obtain the following
result.
\begin{corollary}
	\label{linear_indep_extr}
	Every family of pairwise distinct extreme points of $\mathcal{F}$ is
	linearly independent.
\end{corollary}
\begin{proof}
	Let $(u_i=s_i\mathbf{1}_{E_i}/P(E_i))_{i\in I}$ be a  family of
	pairwise distinct extreme points of $\mathcal{F}$. If there exists
	$\lambda\in\RR^I$ (if $I$ is infinite, we assume that $\lambda$ vanishes
	except on a finite set) such that~${\sum_{i\in I}\lambda_i u_i=0}$,
	then~${\sum_{i\in I}\lambda_i \Diff u_i=0}$. Since, for every $i\neq j$, we
	have $${\mathrm{Supp}(\Diff u_i)\cap\mathrm{Supp}(\Diff u_j)=\partial
	E_i\cap \partial E_j=\emptyset}\,,$$ we obtain that the measures $(\Diff
	u_i)_{i\in I}$ have disjoint support, which yields $\lambda_i=0$ for
	every~${i\in I}$.
\end{proof}
We eventually deduce the main result of this section.
\begin{theorem}
	If $\mathrm{dim}(\mathcal{F})=d<+\infty$ then $\mathcal{F}$ has exactly
	$d+1$ extreme points. It is a $d$-simplex.
	\label{number_extr_points}
\end{theorem}
\begin{proof}
	Let $u_1,...,u_m$ be distinct extreme points of $\mathcal{F}$. From
	\Cref{linear_indep_extr}, we know that $u_1,...,u_m$ are linearly
	independent, which is hence also the case of $u_2-u_1,...,u_m-u_1$. But
	since this last family is contained in the direction space of
	$\mathrm{Aff}(\mathcal{F})$, we obtain $m-1\leq
	\mathrm{dim}(\mathcal{F})=d$.

	Conversely, $\mathcal{F}$ has at least $d+1$ extreme points, otherwise, by
	Carath\'eodory's theorem, it would be contained in a $d-1$-dimensional
	affine space, a contradiction.
\end{proof}

\begin{figure}[htpb]
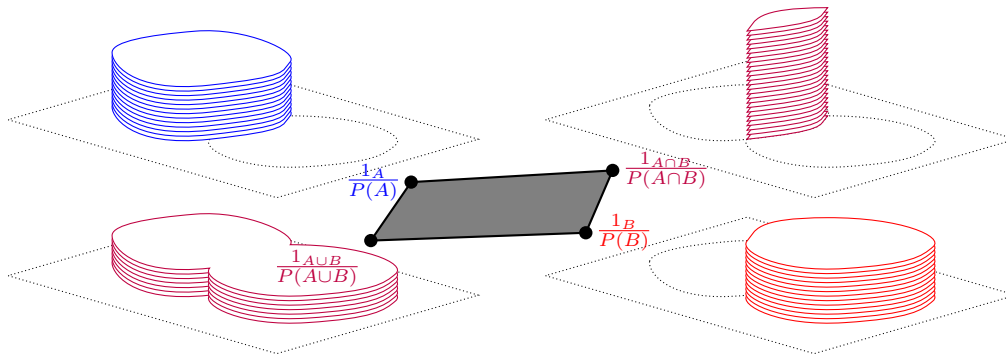

	\centering
\tdplotsetmaincoords{73}{45}
	\includestandalone[scale=0.5]{face4_3d}
	\caption{An example of 2-face of the total variation unit ball $\{\TV\leq 1\}$. Such a face has more than 3 extreme points and therefore it is not exposed by any $\Cder^1$ function, by virtue of~\Cref{number_extr_points}~(notice that $A\cap B$ and $A\cup B$ are not smooth).}
	\label{fig_face4}
\end{figure}

\begin{figure}[htpb]
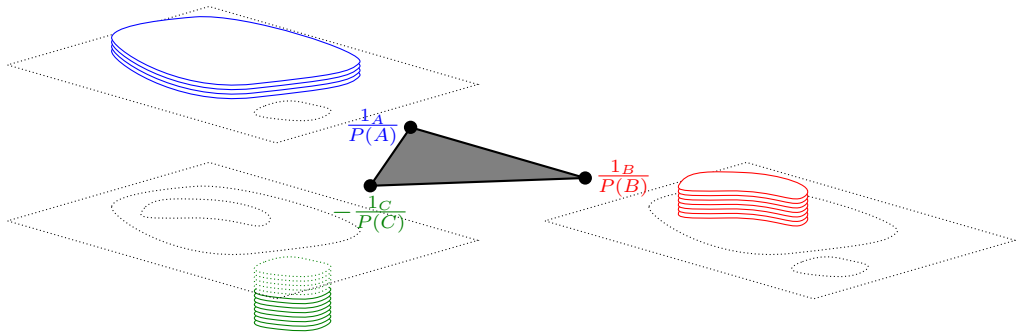

	\centering
\tdplotsetmaincoords{73}{45}
	\includestandalone[scale=0.5]{face3_3d}
	\caption{An example of 2-face of the total variation unit ball $\{\TV\leq 1\}$ that is exposed by some~$\Cder^1$ function. Such a face has  exactly 3 extreme points.}
	\label{fig_face3}
\end{figure}

\Cref{number_extr_points} is illustrated in~\Cref{fig_face4}
and~\Cref{fig_face3}. The 2-face depicted in~\Cref{fig_face4} has more than $3$
extreme points, therefore it is not exposed by any $\Cder^1$ function. On the
contrary, \Cref{fig_face3} illustrates a typical 2-face of $\{\TV\leq 1\}$
exposed by some $\Cder^1$ function: it is a triangle~(2-simplex).

If $\mathrm{dim}(\mathcal{F})=d<+\infty$, Carath\'eodory's theorem implies that every function $u\in\mathcal{F}$ is of the form
\begin{align}
u=\sum\limits_{i \in  I} a_i\mathbf{1}_{E_i}\,, \quad \mbox{where $1\leq \mathrm{card}\,I \leq d+1$},\ a\in(\RR\setminus \{0\})^{I},
	 \label{face_carat_u}
\end{align}
and $\{E_{i}\}_{i\in I}$ is a collection of simple sets with positive finite measure that satisfy  $$\TV(u)=\sum\limits_{i \in I}|a_i|P(E_i)\,.$$ As a consequence of~\Cref{linear_indep_extr} and~\Cref{number_extr_points}, the decomposition~\eqref{face_carat_u} is unique (the~$E_i's$ must form a subcollection of the extreme points of $\mathcal{F}$, and the corresponding $a_i$'s are then uniquely determined).

Coming back to our inverse problem, we deduce that, if some dual certificate $\eta=\Phi^*p$, with~$p$ a solution to \cref{dual_noreg}, exposes some face $\mathcal{F}$ of $\{\TV\leq \min \eqref{primal_noreg}\}$ with dimension $d$, every solution to~\eqref{primal_noreg} has the form \eqref{face_carat_u}. If, moreover, the operator
$$\! \begin{aligned}[t]
	\Phi_{\mathcal{F}} : \RR^{d+1} &  \to \mathcal{H} \\
	a & \mapsto \Phi\Bigg(\sum\limits_{i=1}^{d+1} a_i\mathbf{1}_{E_i}\Bigg)\,,
\end{aligned}$$
is injective, the solution is unique. We see that, in that case,  total (gradient) variation minimization behaves similarly to  $\ell^1$ (synthesis) minimization~\cite{chenAtomicDecomposition1999} or  total variation (of Radon measures) minimization~\cite{candesMathematicalTheorySuperresolution2014}, in the sense that the only faces $\mathcal{F}$ that are involved are simplices. In the next sections, we show that, under some stability assumption given below, that similarity also holds at low noise: not only~$\mathcal{F}$~(the face of the unknown), but all the faces involved in the solutions of~\eqref{primal_reg} for small~$\lambda$ and small $\norm{w}_{\mathcal{H}}$ are simplices, with the same dimension. Hence, with low noise and regularization, the problem~\eqref{primal_reg} behaves like the \textsc{Lasso}~\cite{tibshiraniRegressionShrinkageSelection1996} or the Beurling \textsc{Lasso}~\cite{brediesInverseProblemsSpaces2013,azaisSpikeDetectionInaccurate2015,duvalExactSupportRecovery2015}. 

In our context, the equivalent notion to $k$-sparse vectors (or measures) is the following class of piecewise constant functions.
	\begin{definition}[$k$-simple functions]
		If $k\in \NN^*$, we say that a function $u:\RR^2\to\RR$ is~$k\text{-simple}$ if there exists a collection $\{E_i\}_{1\leq i\leq k}$ of simple sets of class $\Cder^1$ with positive finite measure such that~${\partial E_i\cap \partial E_j=\emptyset}$ for every~$i\neq j$, and $a\in\RR^k$ such that
		\begin{equation*}
			u=\sum\limits_{i=1}^k a_i\mathbf{1}_{E_i}\,.
		\end{equation*}
	\end{definition}
In particular $1$-simple functions are (proportional to) indicators of simple
sets.

The next step is thus to study the stability of $k$-simple functions with
respect to noise and regularization: if $u_0$ is $k$-simple and identifiable,
with $w$ and $\lambda$ small enough, are the solutions of
\customref{primal_reg}{$\mathcal{P}_\lambda(y_0+w)$} $k$-simple? What is the
number of atoms appearing in their decomposition, and how are they related to
those appearing in the decomposition of $u_0$?

%% file: regularity_inter_union.tex
\tikzset{
pattern size/.store in=\mcSize, 
pattern size = 5pt,
pattern thickness/.store in=\mcThickness, 
pattern thickness = 0.3pt,
pattern radius/.store in=\mcRadius, 
pattern radius = 1pt}
\makeatletter
\pgfutil@ifundefined{pgf@pattern@name@_b80mficos}{
\pgfdeclarepatternformonly[\mcThickness,\mcSize]{_b80mficos}
{\pgfqpoint{0pt}{0pt}}
{\pgfpoint{\mcSize+\mcThickness}{\mcSize+\mcThickness}}
{\pgfpoint{\mcSize}{\mcSize}}
{
\pgfsetcolor{\tikz@pattern@color}
\pgfsetlinewidth{\mcThickness}
\pgfpathmoveto{\pgfqpoint{0pt}{0pt}}
\pgfpathlineto{\pgfpoint{\mcSize+\mcThickness}{\mcSize+\mcThickness}}
\pgfusepath{stroke}
}}
\makeatother
\tikzset{every picture/.style={line width=0.75pt}} 

\begin{tikzpicture}[x=0.75pt,y=0.75pt,yscale=-0.5,xscale=0.5]

\draw  [fill={rgb, 255:red, 169; green, 173; blue, 251 }  ,fill opacity=1 ] (340.75,43) .. controls (370.58,42.85) and (395,92.77) .. (395.3,154.48) .. controls (395.6,216.2) and (371.67,266.35) .. (341.84,266.49) .. controls (312.02,266.64) and (287.59,216.73) .. (287.29,155.01) .. controls (286.99,93.29) and (310.92,43.15) .. (340.75,43) -- cycle ;
\draw    (341.3,154.75) -- (473.5,86.5) ;
\draw    (186.5,112.5) -- (264,154) ;
\draw  [pattern=_b80mficos,pattern size=6pt,pattern thickness=0.75pt,pattern radius=0pt, pattern color={rgb, 255:red, 0; green, 0; blue, 0}] (229,154.75) .. controls (229,93.03) and (279.03,43) .. (340.75,43) .. controls (402.47,43) and (452.5,93.03) .. (452.5,154.75) .. controls (452.5,216.47) and (402.47,266.5) .. (340.75,266.5) .. controls (279.03,266.5) and (229,216.47) .. (229,154.75) -- cycle ;

\draw (476,70.4) node [anchor=north west][inner sep=0.75pt]    {$E$};
\draw (160,96.4) node [anchor=north west][inner sep=0.75pt]    {$F$};

\end{tikzpicture}

%% file: face4_3d.tex
\usetikzlibrary{shapes,calc,positioning}
\tdplotsetmaincoords{73}{45}

\begin{tikzpicture}[scale=0.5, tdplot_main_coords,axis/.style={->,dashed},thick]
	\def \L{20}
	\def \l{15}
	\begin{scope}[shift={(-20,0)}]
	\coordinate (A) at (36,7,0);
	\myplan{\L}{\l}
	\foreach \z in {0} 
	\draw[dotted]   (23,5,\z) .. controls (23.95,3.4,\z) and (28.25,2.1,\z) .. (30.15,5.1,\z) .. controls (32.05,8.1,\z) and (31.5,9.5,\z) .. (30.65,11.2,\z) .. controls (29.8,12.9,\z) and (25.75,13.8,\z) .. (23.0,11.0,\z) .. controls (20.25,8.2,\z) and (22.05,6.6,\z) .. (23.0,5.0,\z) -- cycle ;
		\foreach \z in {0,0.25,...,3}
		\filldraw[draw=blue,fill=white]   (15.05,5.7,\z) .. controls (16.2,4.4,\z) and (21.15,0.9,\z) .. (23.0,5.0,\z) .. controls (24.85,9.1,\z) and (23.95,9.4,\z) .. (23.0,11.0,\z) .. controls (22.05,12.6,\z) and (16.9,13.8,\z) .. (15.55,11.4,\z) .. controls (14.2,9.0,\z) and (13.9,7.0,\z) .. (15.05,5.7,\z) -- cycle ;
\end{scope}

\begin{scope}[shift={(20,0)}]
	\coordinate (B) at (9,7,0);
	\myplan{\L}{\l}
		\foreach \z in {0}
		\draw[dotted]   (15.05,5.7,\z) .. controls (16.2,4.4,\z) and (21.15,0.9,\z) .. (23.0,5.0,\z) .. controls (24.85,9.1,\z) and (23.95,9.4,\z) .. (23.0,11.0,\z) .. controls (22.05,12.6,\z) and (16.9,13.8,\z) .. (15.55,11.4,\z) .. controls (14.2,9.0,\z) and (13.9,7.0,\z) .. (15.05,5.7,\z) -- cycle ;
	\foreach \z in {0,0.25,...,3} 
	\filldraw[draw=red,fill=white]   (23,5,\z) .. controls (23.95,3.4,\z) and (28.25,2.1,\z) .. (30.15,5.1,\z) .. controls (32.05,8.1,\z) and (31.5,9.5,\z) .. (30.65,11.2,\z) .. controls (29.8,12.9,\z) and (25.75,13.8,\z) .. (23.0,11.0,\z) .. controls (20.25,8.2,\z) and (22.05,6.6,\z) .. (23.0,5.0,\z) -- cycle ;
\end{scope}
\begin{scope}[shift={(00,-20)}]
	\coordinate (U) at (22,18,0);
	\myplan{\L}{\l}
	\foreach \z in {0,0.25,...,1.5} 
	\filldraw[fill=white,draw=purple]   (15.55,11.4,\z) .. controls (14.2,9.0,\z) and (13.9,7.0,\z) .. (15.05,5.7,\z) .. controls (16.2,4.4,\z) and (21.15,0.9,\z) .. (23.0,5.0,\z) .. controls (23.95,3.4,\z) and (28.25,2.1,\z) .. (30.15,5.1,\z) .. controls (32.05,8.1,\z) and (31.5,9.5,\z) .. (30.65,11.2,\z) .. controls (29.8,12.9,\z) and (25.75,13.8,\z) .. (23.0,11.0,\z) .. controls (22.05,12.6,\z) and (16.9,13.8,\z) .. (15.55,11.4,\z) -- cycle ;
\end{scope}

\begin{scope}[shift={(0,20)}]
	\coordinate (I) at (22,-4,0);
	\myplan{\L}{\l}
		\foreach \z in {0}
		\draw[dotted]   (15.05,5.7,\z) .. controls (16.2,4.4,\z) and (21.15,0.9,\z) .. (23.0,5.0,\z) .. controls (24.85,9.1,\z) and (23.95,9.4,\z) .. (23.0,11.0,\z) .. controls (22.05,12.6,\z) and (16.9,13.8,\z) .. (15.55,11.4,\z) .. controls (14.2,9.0,\z) and (13.9,7.0,\z) .. (15.05,5.7,\z) -- cycle ;
	\foreach \z in {0} 
	\draw[dotted]   (23,5,\z) .. controls (23.95,3.4,\z) and (28.25,2.1,\z) .. (30.15,5.1,\z) .. controls (32.05,8.1,\z) and (31.5,9.5,\z) .. (30.65,11.2,\z) .. controls (29.8,12.9,\z) and (25.75,13.8,\z) .. (23.0,11.0,\z) .. controls (20.25,8.2,\z) and (22.05,6.6,\z) .. (23.0,5.0,\z) -- cycle ;
	\foreach \z in {0,0.25,...,6} 
	\filldraw[fill=white,draw=purple]   (23,5,\z) .. controls (24.85,9.1,\z) and (23.95,9.4,\z) .. (23.0,11.0,\z) .. controls (20.25,8.2,\z) and (22.05,6.6,\z) .. (23.0,5.0,\z) -- cycle ;
\end{scope}

\draw[fill=gray,draw=black,ultra thick] (A) -- (I) -- (B) -- (U) -- cycle;
\fill (A) circle (10pt);
\fill (B) circle (10pt);
\fill (U) circle (10pt);
\fill (I) circle (10pt);
\node[left,scale=2,blue] at (A) {$\frac{1_A}{P(A)}$};
\node[right,scale=2,red] at (B) {$\frac{1_B}{P(B)}$};
\node[above,right,scale=2,purple] at (I) {$\frac{1_{A\cap B}}{P(A\cap B)}$};
\node[below left,scale=2,purple] at (U) {$\frac{1_{A\cup B}}{P(A\cup B)}$};
\end{tikzpicture}

%% file: face3_3d.tex
\usetikzlibrary{shapes,calc,positioning}
\tdplotsetmaincoords{73}{45}

\tikzset{every picture/.style={line width=0.75pt}} 

\tikzset{every picture/.style={line width=0.75pt}} 

%
%
%
%
%

\begin{tikzpicture}[scale=0.5, tdplot_main_coords,axis/.style={->,dashed},thick]
	\def \L{20}
	\def \l{15}
	\begin{scope}[shift={(-20,0)}]
	\coordinate (A) at (36,7,0);
	\myplan{\L}{\l}
		\foreach \z in {0,0.25,...,0.75} 
		\filldraw[draw=blue,fill=white]  (15.05,5.7,\z) .. controls (16.2,4.4,\z) and (24.477,0.877,\z) .. (26.327,4.977,\z) .. controls (28.177,9.077,\z) and (28.477,10.077,\z) .. (27.527,11.677,\z) .. controls (26.577,13.277,\z) and (16.9,13.8,\z) .. (15.55,11.4,\z) .. controls (14.2,9.0,\z) and (13.9,7.0,\z) .. (15.05,5.7,\z) -- cycle ;
	\foreach \z in {0} 
	\draw[dotted,fill=white]   (27.827,4.277,\z) .. controls (28.127,2.177,\z) and (29.388,2.807,\z) .. (30.388,3.307,\z) .. controls (31.388,3.807,\z) and (31.688,4.407,\z) .. (31.288,5.507,\z) .. controls (30.888,6.607,\z) and (30.288,6.407,\z) .. (29.527,6.277,\z) .. controls (28.765,6.147,\z) and (27.527,6.377,\z) .. (27.827,4.277,\z) -- cycle ;
\end{scope}

\begin{scope}[shift={(20,0)}]
	\coordinate (B) at (9,7,0);
	\myplan{\L}{\l}
	\foreach \z in {0} 
		\draw[dotted]  (15.05,5.7,\z) .. controls (16.2,4.4,\z) and (24.477,0.877,\z) .. (26.327,4.977,\z) .. controls (28.177,9.077,\z) and (28.477,10.077,\z) .. (27.527,11.677,\z) .. controls (26.577,13.277,\z) and (16.9,13.8,\z) .. (15.55,11.4,\z) .. controls (14.2,9.0,\z) and (13.9,7.0,\z) .. (15.05,5.7,\z) -- cycle ;
	\foreach \z in {0} 
	\draw[dotted]   (27.827,4.277,\z) .. controls (28.127,2.177,\z) and (29.388,2.807,\z) .. (30.388,3.307,\z) .. controls (31.388,3.807,\z) and (31.688,4.407,\z) .. (31.288,5.507,\z) .. controls (30.888,6.607,\z) and (30.288,6.407,\z) .. (29.527,6.277,\z) .. controls (28.765,6.147,\z) and (27.527,6.377,\z) .. (27.827,4.277,\z) -- cycle ;
	\foreach \z in {0,0.25,...,1.5} 
	\filldraw[draw=red,fill=white]   (16.427,6.777,\z) .. controls (17.377,5.177,\z) and (17.627,5.177,\z) .. (19.427,7.177,\z) .. controls (21.227,9.177,\z) and (23.777,7.777,\z) .. (22.927,9.477,\z) .. controls (22.077,11.177,\z) and (19.227,11.277,\z) .. (17.727,10.177,\z) .. controls (16.227,9.077,\z) and (15.477,8.377,\z) .. (16.427,6.777,\z) -- cycle ;
\end{scope}
\begin{scope}[shift={(00,-20)}]
	\coordinate (U) at (22,18,0);
	\foreach \z in {-3,-2.75,...,-1.25}  
	\filldraw[fill=white,draw=darkgreen]   (27.827,4.277,\z) .. controls (28.127,2.177,\z) and (29.388,2.807,\z) .. (30.388,3.307,\z) .. controls (31.388,3.807,\z) and (31.688,4.407,\z) .. (31.288,5.507,\z) .. controls (30.888,6.607,\z) and (30.288,6.407,\z) .. (29.527,6.277,\z) .. controls (28.765,6.147,\z) and (27.527,6.377,\z) .. (27.827,4.277,\z) -- cycle ;
	\foreach \z in {-1.0,-0.75,...,0}  
	\filldraw[fill=white,draw=darkgreen,dotted]   (27.827,4.277,\z) .. controls (28.127,2.177,\z) and (29.388,2.807,\z) .. (30.388,3.307,\z) .. controls (31.388,3.807,\z) and (31.688,4.407,\z) .. (31.288,5.507,\z) .. controls (30.888,6.607,\z) and (30.288,6.407,\z) .. (29.527,6.277,\z) .. controls (28.765,6.147,\z) and (27.527,6.377,\z) .. (27.827,4.277,\z) -- cycle ;
	\myplan{\L}{\l}
	\foreach \z in {0} 
		\draw[dotted]  (15.05,5.7,\z) .. controls (16.2,4.4,\z) and (24.477,0.877,\z) .. (26.327,4.977,\z) .. controls (28.177,9.077,\z) and (28.477,10.077,\z) .. (27.527,11.677,\z) .. controls (26.577,13.277,\z) and (16.9,13.8,\z) .. (15.55,11.4,\z) .. controls (14.2,9.0,\z) and (13.9,7.0,\z) .. (15.05,5.7,\z) -- cycle ;
	\foreach \z in {0} 
	\draw[dotted]   (16.427,6.777,\z) .. controls (17.377,5.177,\z) and (17.627,5.177,\z) .. (19.427,7.177,\z) .. controls (21.227,9.177,\z) and (23.777,7.777,\z) .. (22.927,9.477,\z) .. controls (22.077,11.177,\z) and (19.227,11.277,\z) .. (17.727,10.177,\z) .. controls (16.227,9.077,\z) and (15.477,8.377,\z) .. (16.427,6.777,\z) -- cycle ;
\end{scope}

\draw[fill=gray,draw=black,ultra thick] (A) --  (B) -- (U) -- cycle;
\fill (A) circle (10pt);
\fill (B) circle (10pt);
\fill (U) circle (10pt);
\node[left,scale=2,blue] at (A) {$\frac{1_A}{P(A)}$};
\node[right,scale=2,red] at (B) {$\frac{1_B}{P(B)}$};
\node[below,scale=2,darkgreen] at (U) {$-\frac{1_{C}}{P(C)}$};
\end{tikzpicture}

%% file: curv_pb.tex

\section{Stability analysis of the prescribed curvature problem}
\label{sec_curv_pb}
The simple sets appearing in the decomposition of any solution to \customref{primal_reg}{$\mathcal{P}_\lambda(y_0+w)$} are all solutions of the prescribed curvature problem associated to $\eta_{\lambda,w}$. In \Cref{sec_dual}, we have also seen that, under a few assumptions,~$\eta_{\lambda,w}$ converges to the minimal norm certificate $\eta_0$ when $w$ and $\lambda$ go to zero. It is therefore natural to investigate how solutions of the prescribed curvature problem behave under variations of the curvature functional.

In this section, we consider the prescribed curvature problem \cref{prescribed_curv} associated to some function~${\eta\in\partial\TV(0)\cap\mathrm{C}^1(\RR^2)}$. We investigate how the solution set of \cref{prescribed_curv} behaves when $\eta$ varies. To be more specific, given two sufficiently close curvature functionals $\eta$ and $\eta'$, we address the following two questions.
\begin{enumerate}[label=(\roman*)]
	\item Are the solutions to \customref{prescribed_curv}{$\mathcal{PC}(\eta')$} close to some solutions to \customref{prescribed_curv}{$\mathcal{PC}(\eta)$}?
	\item How many solutions to \customref{prescribed_curv}{$\mathcal{PC}(\eta')$} are there in a neighborhood of a given solution to~\customref{prescribed_curv}{$\mathcal{PC}(\eta)$}?
\end{enumerate}
We answer the first question using the notion of quasi-minimizers of the perimeter, as well as first order optimality conditions for \cref{prescribed_curv}. Then, under a strict stability assumption on solutions to \cref{prescribed_curv}, we answer the second question using second order shape derivatives.

\paragraph{Convergence result. } First, we tackle Question (i) with the following proposition which states that any neighborhood~(in terms of $\mathrm{C}^2$-normal deformations) of the solution set of~\customref{prescribed_curv}{${\mathcal{PC}(\eta_0)}$} contains the solution set of \cref{prescribed_curv} provided~$\eta$ is sufficiently
close to $\eta_0$ in $\mathrm{C}^1(\RR^2)$ and $\LD$. The proof, which relies on standard compactness results for quasi-minimizers of the perimeter, is postponed to \Cref{appendix_proof}.

\begin{proposition}
		Let $\eta_0\in\partial\TV(0)\cap\mathrm{C}^{1}(\RR^2)$. For every $\epsilon>0$ there exists
		$r>0$ such that for every~${\eta\in\partial\TV(0)}\cap\mathrm{C}^1(\RR^2)$ with~${\|\eta-\eta_0\|_{\LD}+\|\eta-\eta_0\|_{\mathrm{C}^1(\RR^2)}\leq r}$, the
		following holds: every non-empty solution $F$ of \cref{prescribed_curv} is a $\mathrm{C}^2$-normal deformation of size at most $\epsilon$ of a non-empty solution $E$ of \customref{prescribed_curv}{$\mathcal{PC}(\eta_0)$}, i.e., using the notation of \Cref{normal_deform_boundary}, $F=E_\varphi$ with~${\|\varphi\|_{\mathrm{C}^2(\partial E)}\leq \epsilon}$.
		\label{conv_result}
\end{proposition}

\subsection{Stability result}
\label{sec_stability}

Question (ii) is closely linked to the stability of minimizers to \cref{prescribed_curv}, that is to the behaviour of the objective $J$ in a neighborhood of a solution. To analyze this behaviour, we use the general framework presented in~\cite{dambrineStabilityShapeOptimization2019}, which relies on the notion of second order shape derivative. In this section, unless otherwise stated, $E$ denotes a non-empty bounded open set of class $\mathrm{C}^2$.

\paragraph{Approach.} The natural path to obtain our main stability result, which is \Cref{prop_stability}, is to prove that~$J$ is in some sense of class $\mathrm{C}^2$, i.e. that its second order shape derivative is continuous at zero (see \Cref{lemma_cont_shape_hessian_3} for a precise statement). Although it is likely to be known, we could not find this result in the literature. We postpone its proof to \Cref{proof_stab_prob_pc}. To obtain~\Cref{prop_stability}, we had to use a stronger condition than the ``improved continuity condition'' $(\mathbf{IC}_{\mathrm{H}^1,\mathrm{C}^{2}})$ of~\cite{dambrineStabilityShapeOptimization2019}, which is satisfied by our functional. The latter only requires some uniform control of second order directional derivatives at zero, which is weaker than the result of \Cref{lemma_cont_shape_hessian_3}.

\paragraph{Structure of shape derivatives.} We introduce the following mapping, where $E_\varphi$ denotes the normal deformation of $E$ associated to $\varphi$, defined in \Cref{normal_deform_boundary}:
$$\! \begin{aligned}[t]
\mathrm{j}_E :  \mathrm{C}^{1}(\partial E)&  \to \RR \\
\varphi   & \mapsto J(E_\varphi)\,.
\end{aligned}$$
With this notation, the following result holds.
	\begin{proposition}[{{See e.g. \cite[Chapter 5]{henrotShapeVariationOptimization2018}}}]
		If $\eta\in\mathrm{C}^1(\RR^2)$, then $\mathrm{j}_E$ is twice Fr\'echet differentiable at $0$ and, for every~$\psi\in\mathrm{C}^{1}(\partial E)$, we have:
		\begin{equation*}
			\begin{aligned}
				\mathrm{j}'_E(0).(\psi)&=\int_{\partial E}\left[H-\eta\right]\psi\,d\mathcal{H}^1\\
				\mathrm{j}''_E(0).(\psi,\psi)&=\int_{\partial E}\left[|\nabla_{\tau}\psi|^2-\left(H\,\eta+\frac{\partial \eta}{\partial \nu}\right)\psi^2\right]d\mathcal{H}^1
			\end{aligned}
		\end{equation*}
		where $H$ denotes the curvature of~$E$ and $\nabla_{\tau}\psi\eqdef\nabla \psi-(\nabla\psi\cdot \nu)\,\nu$ is the tangential gradient of $\psi$ with respect to $E$.
	\end{proposition}
From the expression of $\mathrm{j}'_E(0)$ and $\mathrm{j}''_E(0)$ given above, we immediately notice that $\mathrm{j}'_E(0)$ can be extended to a continuous linear form on~$\mathrm{L}^1(\partial E)$, and~$\mathrm{j}''_E(0)$ to a continuous bilinear form on~${\mathrm{H}^1(\partial E)}$.

\paragraph{Strict stability.} Following \cite{dambrineStabilityShapeOptimization2019}, we say that a non-empty open solution~$E$ of \customref{prescribed_curv}{$\mathcal{PC}(\eta)$} is strictly stable if~$\mathrm{j}''_E(0)$ is coercive in~${\mathrm{H}^1(\partial E)}$, i.e. if the following property holds:
\begin{equation*}
	\exists \alpha>0,~\forall\psi\in\mathrm{H}^1(\partial E),\quad\mathrm{j}''_E(0).(\psi,\psi)\geq \alpha\,\|\psi\|_{\mathrm{H}^1(\partial E)}^2\,.
\end{equation*}
As noticed by Dambrine and Lamboley, this strict stability condition is a key ingredient (together with several assumptions) to ensure that $E$ is a strict local minimizer of $J$ (see Theorem 1.1 in the above-mentioned reference), and is hence the only minimizer among the sets $E_\varphi$ with $\varphi$ in a neighborhood of $0$. It plays a crucial role in our answer to Question~(ii).

\paragraph{Continuity results.} Now, we state two important results concerning the convergence of $\mathrm{j}''_E$
towards~$\mathrm{j}''_{0,E}$ and the continuity of~${\varphi\mapsto \mathrm{j}''_E(\varphi)}$, where $\mathrm{j}_E$ and $\mathrm{j}_{0,E}$ are
the functionals respectively associated to $\eta$ and $\eta_0$. Their proof is postponed to \Cref{proof_stab_prob_pc}.  In all the following, if $X$ is a~(real) vector space, we denote by $\mathcal{Q}(X)$ the set of quadratic forms over $X$, and define $\|\cdot\|_{\mathcal{Q}(X)}$ as follows:
$$\|q\|_{\mathcal{Q}(X)}\eqdef \underset{x\in X\setminus\{0\}}{\mathrm{sup}}~\frac{\abs{q(x,x)}}{\|x\|_X^2}$$
	\begin{proposition}
		\label{lemma_cont_shape_hessian_3}
		If $\eta\in\mathrm{C}^1(\RR^2)$, the mapping
		$$\! \begin{aligned}[t]
		\mathrm{j}''_E : \mathrm{C}^{2}(\partial E) &  \to \mathcal{Q}(\mathrm{H}^1(\partial E))  \\
		\varphi & \mapsto \mathrm{j}''_E(\varphi)
		\end{aligned}$$
		is continuous at $0$.
	\end{proposition}
\begin{proposition}
	\label{lemma_conv_shape_hessian}
	Let $\eta_0\in\mathrm{C}^1(\RR^2)$. There exists $\epsilon>0$ such that
	$$\underset{\|\eta-\eta_0\|_{\mathrm{C}^1(\RR^2)}\to 0}{\mathrm{lim}}~\underset{\|\varphi\|_{\mathrm{C}^2(\partial E)}\leq \epsilon}{\mathrm{sup}}~\left\|\mathrm{j}''_E(\varphi)-\mathrm{j}''_{0,E}(\varphi)\right\|_{\mathcal{Q}(\mathrm{H}^1(\partial E))}=0\,.$$
\end{proposition}

\paragraph{Stability result.} We are now able to state the final result of this section, which states that if~$E$ is a strictly stable solution to~\customref{prescribed_curv}{$\mathcal{PC}(\eta_0)$}, there is at most one $\varphi$ in a neighborhood of $0$ such that~$E_{\varphi}$ is a solution to~\customref{prescribed_curv}{$\mathcal{PC}(\eta)$}, provided $\|\eta-\eta_0\|_{\mathrm{C}^1(\RR^2)}$ is small engouh.
\begin{proposition}
	Let $\eta_0\in\partial\TV(0)\cap \mathrm{C}^{1}(\RR^2)$ and $E$ be a strictly stable solution to \customref{prescribed_curv}{$\mathcal{PC}(\eta_0)$}. Then there
	exists~${\epsilon>0}$ and $r>0$ such that for every $\eta\in\partial\TV(0)$ with $\|\eta-\eta_0\|_{\mathrm{C}^1(\RR^2)}\leq r$ there is at most one~${\varphi\in\mathrm{C}^2(\partial E)}$ such that~${\|\varphi\|_{\mathrm{C}^2(\partial E)}\leq \epsilon}$ and $E_{\varphi}$ solves \customref{prescribed_curv}{$\mathcal{PC}(\eta)$}.
	\label{prop_stability}
\end{proposition}
\begin{proof}
	The fact that $E$ is a strictly stable solution to \customref{prescribed_curv}{$\mathcal{PC}(\eta_0)$} and the results above give the existence of $\epsilon>0$, $r>0$ and
	$\alpha>0$ such that, for every $(\varphi,\eta)\in\mathrm{C}^2(\partial E)\times \mathrm{C}^1(\RR^2)$ with~$\|\varphi\|_{\mathrm{C}^2(\partial E)}\leq
	\epsilon$ and~${\|\eta-\eta_0\|_{\mathrm{C}^1(\RR^2)}\leq r}$, we have:
	$$\underset{\psi\in\mathrm{H}^1(\partial E)\setminus\{0\}}{\mathrm{sup}}~\frac{\mathrm{j}''_{E}(\varphi).(\psi,\psi)}{\|\psi\|_{\mathrm{H}^1(\partial E)}^2}\geq \alpha$$
	As a result, $\mathrm{j}''_E(\varphi)$ is coercive (and hence positive definite) for every $\varphi$ such that $\|\varphi\|_{\mathrm{C}^2(\partial E)}\leq
	\epsilon$. We therefore obtain that $\mathrm{j}_E$ is strictly convex on this set and the result follows.
\end{proof}

\paragraph{Summary.} Combining the results of \Cref{conv_result,prop_stability}, we have proved that, provided~$\eta$ is sufficiently
close to $\eta_0$ in $\mathrm{C}^1(\RR^2)$ and $\LD$, every solution to \cref{prescribed_curv} belongs to a neighborhood~(in terms of $\mathrm{C}^2$-normal deformations) of a solution to \customref{prescribed_curv}{${\mathcal{PC}(\eta_0)}$}, and that, under a strict stability assumption, each of these neighborhoods contains at most one solution to \cref{prescribed_curv}. In~\Cref{sec_suff_stab} below, we discuss this strict stability assumption in greater details. Then, in~\Cref{main_th}, we prove (under suitable assumptions) that, if $\eta=\eta_{\lambda,w}$ is the dual certificate associated to~\customref{primal_reg}{$\mathcal{P}_{\lambda}(y_0+w)$} and $\eta_0$ the minimal norm dual certificate associated to \cref{primal_noreg}, then each neighborhood of a solution to \customref{prescribed_curv}{${\mathcal{PC}(\eta_0)}$} contains exactly one solution to~\customref{prescribed_curv}{${\mathcal{PC}(\eta_{\lambda,w})}$}.

\subsection{A sufficient condition for strict stability}
\label{sec_suff_stab}
As mentioned above, we here discuss how to ensure that a non-empty open solution to \cref{prescribed_curv} is strictly stable. We derive a sufficient condition for this property to hold, and then discuss to what extent it is necessary.

\paragraph{Setting.} We fix $\eta\in\partial\TV(0)\cap\mathrm{C}^1(\RR^2)$ and $E$ a non-empty open solution to~\cref{prescribed_curv}. We recall that, necessarily, $H_E=\eta$ on $\partial E$, and  the quadratic form $\mathrm{j}''_E(0)$ is positive semi-definite. By definition, the set $E$ is a strictly stable solution to \cref{prescribed_curv} if and only if $\mathrm{j}''_E(0)$ is coercive in~$\mathrm{H}^1(\partial E)$, with
\begin{equation*}
	\forall\psi\in\mathrm{H}^1(\partial E),~\mathrm{j}''_E(0).(\psi,\psi)=\int_{\partial E}\left[|\nabla_{\tau_E}\psi|^2-\left(H_E^2+\frac{\partial \eta}{\partial \nu_E}\right)\psi^2\right]d\mathcal{H}^1\,.
\end{equation*}

\paragraph{Equivalence of coercivity and positive definiteness.} As explained (in a more general context) in~\cite{dambrineStabilityShapeOptimization2019}, the bilinear form $\mathrm{j}''_E(0)$ is in fact coercive if and only if it is positive definite. Our functional $J$ fits the assumptions of~Lemma 3.1 in the above reference. Indeed, writing $\mathrm{j}''_E(0)=\ell_m+\ell_r$ with 
\begin{equation*}
	\begin{aligned}
		\ell_m(\psi,\psi)&=\int_{\partial E}|\nabla_{\tau_E}\psi|^2\,,\\
		\ell_r(\psi,\psi)&=-\int_{\partial E}\left[\left(H_E^2+\frac{\partial \eta}{\partial \nu_E}\right)\psi^2\right]d\mathcal{H}^1\,,
	\end{aligned}
\end{equation*}
we see that $\mathrm{j}''_E(0)$ satisfies $(\mathbf{C}_{\mathrm{H}^{s_2}})$ with $s_1=s_2=1$.
We consequently obtain the following result, which we do not use in the following but which is interesting in itself.
\begin{lemma}[{\cite[Lemma 3.1]{dambrineStabilityShapeOptimization2019}}] The following two propositions are equivalent:
	\begin{enumerate}[label=(\roman*)]
		\item $\mathrm{j}''_E(0)$ is positive definite, i.e. $$\forall\psi\in\mathrm{H}^1(\partial E)\setminus\{0\},\quad\mathrm{j}''_E(0).(\psi,\psi)>0\,,$$
		\item $\mathrm{j}''_E(0)$ is coercive, i.e. $$\exists \alpha>0,~\forall\psi\in\mathrm{H}^1(\partial E),\quad \mathrm{j}''_E(0).(\psi,\psi)\geq \alpha\|\psi\|_{\mathrm{H}^1(\partial E)}^2\,.$$
	\end{enumerate}
\label{eq_coerciv_pos}
\end{lemma}

\paragraph{A sufficient condition for coercivity.} Using the expression of $\mathrm{j}''_E(0)$, the following result can be directly obtained.
\begin{proposition} If
	\begin{equation}
		\underset{x\in\partial E}{\mathrm{sup}}\left[H_E(x)^2+\frac{\partial\eta}{\partial\nu_E}(x)\right]<0\,,
		\label{suff_coerciv_eq}
	\end{equation}
	then $\mathrm{j}''_E(0)$ is coercive.
	\label{suff_coerciv}
\end{proposition}

\paragraph{Necessity of the condition?} A natural question is whether the condition in \Cref{suff_coerciv} is necessary. We conjecture that it is not the case. Indeed, assuming that $E$ is simple and $\gamma$ is an arc-length parametrization of $\partial E$, we have:
\begin{equation}
	\forall \psi\in\mathrm{H}^1(\partial E),~\mathrm{j}''_E(0).(\psi,\psi)=\int_{I}(\psi\circ \gamma)'^2- \bigg(\bigg(H_E^2+\frac{\partial\eta}{\partial \nu_E}\bigg)\circ\gamma\bigg)\,(\psi\circ\gamma)^2\,,
	\label{change_var}
\end{equation}
where $I\eqdef (0,\mathcal{H}^1(\partial E))$. The existence of $\psi\neq 0$ such that $\mathrm{j}''_E(0).(\psi,\psi)=0$\footnote{We recall that $\mathrm{j}''_E(0)$ is positive semi-definite, and hence that it is positive definite (or by \Cref{eq_coerciv_pos}, equivalently, coercive) if and only if $\mathrm{j}''_E(0).(\psi,\psi)=0$ implies $\psi=0$.} is hence equivalent to the existence of a non-zero minimizer of $\int_{I} \varphi'^2+ V\varphi^2$ under periodicity constraint, where $${V\eqdef-\left(H_E^2+\frac{\partial\eta}{\partial\nu_E}\right)\circ\gamma}\,.$$ The first order optimality condition associated to this problem writes~${-\varphi''+V\varphi=0}$. The coercivity of $\mathrm{j}''_E(0)$ can hence be related to the spectrum of the Schr\"odinger operator with periodic boundary conditions associated to $V$. It is known that there exist potentials $V$ which are not positive and yet correspond to positive definite Schr\"odinger operators. Therefore, it might be possible to construct examples where \eqref{suff_coerciv_eq} does not hold and yet $\mathrm{j}''_E(0)$ is positive definite.

 However, as we explain in \Cref{almost_ness_coerciv}, if $H_E^2+\frac{\partial\eta}{\partial \nu_E}\geq \alpha$ on a connected portion $\Gamma$ of~$\partial E$ and $\alpha\geq (\pi/\mathcal{H}^1(\Gamma))^2$, we are able to prove that~$\mathrm{j}''_E(0)$ is not coercive. Let us consider a $\mathrm{C}^1$ simple open curve $\Gamma$ with finite length. We define the first Dirichlet eigenvalue of the Laplacian associated to $\Gamma$\footnote{We refer the reader to e.g.~\cite{kuttlerEigenvaluesLaplacianTwo1984} for the more classical case of open bounded sets.} by:
\begin{equation}
	\lambda_1(\Gamma)\eqdef\underset{\psi\in\mathrm{H}_0^{1}(\Gamma)\setminus\{0\}}{\mathrm{inf}}~
	\frac{\|\nabla_{\tau_\Gamma}\psi\|^2_{\mathrm{L}^2(\Gamma)}}{\|\psi\|^2_{\mathrm{L}^2(\Gamma)}}\,.
	\label{eigenvalue_curve}
\end{equation}
Using a change of variable as in \Cref{change_var}, one can see that the infimum in \cref{eigenvalue_curve} is attained and is actually equal to the Dirichlet eigenvalue of the interval~${I=(0,\mathcal{H}^1(\Gamma))\subset\RR}$, which is $(\pi/ \mathcal{H}^1(\Gamma))^2$. Using this fact, we can now prove the following result.
	\begin{proposition}
		If there exists $\alpha>0$ such that $H_E^2+\frac{\partial\eta}{\partial\nu_E}\geq \alpha$ on a  connected subset $\Gamma$ of $\partial E$ with $\alpha\geq (\pi/\mathcal{H}^1(\Gamma))^2$, then $\mathrm{j}''_E(0)$ is not coercive.
		\label{almost_ness_coerciv}
	\end{proposition}
\begin{proof}
	Since the infimum in the definition of $\lambda_1(\Gamma)$ is attained, we have the existence of a nonzero function~${\varphi\in\mathrm{H}^1_0(\Gamma)}$ such that
	\begin{equation*}
		\frac{\|\nabla_{\tau_\Gamma}\varphi\|^2_{\mathrm{L}^2(\Gamma)}}{\|\varphi\|^2_{\mathrm{L}^2(\Gamma)}}=\lambda_1(\Gamma)=\left(\frac{\pi}{\mathcal{H}^1(\Gamma)}\right)^2\leq \alpha\,.
	\end{equation*}
	We hence obtain
	\begin{equation*}
		\int_{\Gamma}\left[|\nabla_{\tau_\Gamma}\varphi|^2-\left(H_E^2+\frac{\partial\eta}{\partial\nu_E}\right)\varphi^2\right]d\mathcal{H}^1\leq \int_{\Gamma}\left[|\nabla_{\tau_\Gamma}\varphi|^2-\alpha\,\varphi^2\right]d\mathcal{H}^1\leq 0\,.
	\end{equation*}
	We can then extend $\varphi$ to $\psi\in\mathrm{H}^1(\partial E)$ whose support is compactly included in $\Gamma$, which yields
	\begin{equation*}
		\int_{\partial E}\left[|\nabla_{\tau_E}\psi|^2-\left(H_E^2+\frac{\partial\eta}{\partial\nu_E}\right)\psi^2\right]d\mathcal{H}^1\leq 0\,.
	\end{equation*}
	We can therefore conclude that $\mathrm{j}''_E(0)$ is not coercive.
\end{proof}

%% file: main_res.tex

\section{Exact support recovery}
\label{sec_supp_rec}
To obtain our support recovery result, which is \Cref{main_th}, we first prove a
stability result for the exposed faces of the total variation unit ball, which
is \Cref{th_faces}.

\subsection{Stability of the exposed faces of the total variation unit ball}

\paragraph{Notations and definitions.} In the following, if $\eta$ (resp.
$\eta_n$) belongs to $\partial\TV(0)\cap\mathrm{C}^1(\RR^2)$, we denote by
$\mathcal{F}$ (resp. by $\mathcal{F}_n$) the face of $\{\TV\leq 1\}$ exposed by
$\eta$. We say that~$s\,\mathbf{1}_E/P(E)\in\mathrm{extr}(\mathcal{F})$
is strictly stable if $s=1$ and $E$ is a strictly stable solution to
\customref{prescribed_curv}{$\mathcal{PC}(\eta)$} or $s=-1$ and $E$ is a
strictly stable solution to~\customref{prescribed_curv}{$\mathcal{PC}(-\eta)$}.

\begin{theorem}
	Let $\eta_0\in\partial \TV(0)\cap\mathrm{C}^1(\RR^2)$ be such that
	$\mathcal{F}_0$ has finite dimension, with all its extreme points strictly
	stable. Then for every $\epsilon>0$, there exists $r>0$ such that, for
	every~${\eta\in\partial\TV(0)\cap\mathrm{C}^1(\RR^2)}$ with
	$$\|\eta-\eta_0\|_{\LD}+\|\eta-\eta_0\|_{\mathrm{C}^1(\RR^2)}\leq r\,,$$
	there exists an injective mapping $\theta:\mathrm{extr}(\mathcal{F})\to
	\mathrm{extr}({\mathcal{F}_0})$ such that, for
	every~${u=s\,\mathbf{1}_{F}/P(F)}$ in~$\mathrm{extr}(\mathcal{F})$,
	we have $\theta(u)=s\,\mathbf{1}_E/P(E)$ with
	$$F=E_{\varphi}~\mathrm{and}~\|\varphi\|_{\mathrm{C}^2(\partial E)}\leq
	\epsilon\,.$$ In particular $\mathrm{dim}(\mathcal{F})\leq
	\mathrm{dim}(\mathcal{F}_0)$.
	\label{th_faces}
\end{theorem}

To prove \Cref{th_faces}, we rely on \Cref{lemma_th_faces} below and its
corollary.
\begin{lemma}
Let $(\eta_n)_{n\in\NN^*}$ be a sequence of functions in $\partial\TV(0)\cap
\mathrm{C}^1(\RR^2)$ converging in $\LD$ and $\mathrm{C}^1(\RR^2)$ to $\eta_0$.
Assume that $\mathcal{F}_0$ has finite dimension, with all its extreme points
strictly stable, and that there are infinitely many $n\in\NN^*$ such that
$\mathcal{F}_n$ has at least $m$ pairwise distinct extreme points,
say~$(s_{n,i}\,\mathbf{1}_{E_{n,i}}/P(E_{n,i}))_{1\leq i\leq m}$. Then
there exists $(s_i)_{1\leq i\leq m}$ and pairwise distinct sets
$(E_i)_{1\leq i\leq m}$ such
that~${s_i\,\mathbf{1}_{E_i}/P(E_i)\in\mathrm{extr}({F}_0})$ for all
$i\in\{1,...,m\}$ and, up to the extraction of a (not relabeled) subsequence,
	\begin{equation}
		\forall n\in\NN^*,~\forall i\in\{1,...,m\},\quad
		\left\{\begin{aligned}
			&s_{n,i}=s_i,\\
			&E_{n,i}=(E_{i})_{\varphi_{n,i}}~\mathrm{with}~\underset{n\to+\infty}{\mathrm{lim}}~\|\varphi_{n,i}\|_{\mathrm{C}^2(\partial E_{i})}=0.
		\end{aligned}\right.
		\label{chain_deform}
	\end{equation}
	In particular $m\leq \mathrm{card}(\mathrm{extr}(\mathcal{F}_0))$.
	\label{lemma_th_faces}
\end{lemma}
\begin{proof}
	For every $i\in\{1,...,m\}$, there are infinitely many $n\in\NN^*$ such that
	$s_{n,i}=1$, or infinitely many $n\in\NN^*$ such that
	$s_{n,i}=-1$. Hence, there exists $(s_i)_{1\leq
	i\leq m}$ such that, up to the extraction of a (not relabeled) subsequence ,
	$s_{n,i}=s_i$ for all $n\in\NN^*$ and $i\in\{1,...,m\}$. Now,
	from~\Cref{conv_result}, up to the extraction of a
	subsequence, for every $i\in\{1,...,m\}$, the sequence~$(E_{n,i})_{n\in\NN^*}$ converges in $\mathrm{C}^3$ towards a solution $E_i$
	of \customref{prescribed_curv}{$\mathcal{PC}(s_i\eta_0)$}, which
	yields \Cref{chain_deform}. Moreover, since $E_{n,i}$ is simple and
	diffeomorphic to $E_i$ for $n$ large enough, we obtain that $E_i$ is simple
	and hence~$s_i\,\mathbf{1}_{E_i}/P(E_i)\in\mathrm{extr}(\mathcal{F}_0)$.

	Now, let us prove that the $(E_i)_{1\leq i\leq m}$ are pairwise distinct.
	By contradiction, if~$E_i=E_j$ for some $i\neq j$, then\footnote{The fact that $P(E_i)>0$ follows from $E_i$ being nonempty and not $\mathbb{R}^2$, see \Cref{conv_result}.}
	\begin{align*}
		s_i = \frac{\int_{E_i}\eta_0}{P(E_i)}= \frac{\int_{E_j}\eta_0}{P(E_j)}=s_j\,.
	\end{align*}
Thus, there would exist
	two distinct solutions of
	$\customref{prescribed_curv}{$\mathcal{PC}(\epsilon_i\eta_n)$}$ (namely
	$E_{n,i}$ and $E_{n,j}$) in arbitrarily small neighborhoods of $E_i$, which
	would contradict its strict stability (\Cref{prop_stability}).
\end{proof}

\begin{proof}[Proof of \Cref{th_faces}]
	By contradiction, we assume the existence of some $\epsilon>0$ and of some sequence 
	$(\eta_n)_{n\in\NN^*}$  in
	$\partial\TV(0)\cap\mathrm{C}^1(\RR^2)$ converging in $\LD$ and
	$\mathrm{C}^1(\RR^2)$ to $\eta_0$, and such that, for all~${n\in\NN^*}$, the
	claimed property does not hold. 

	Let  $m=\limsup_{n\to +\infty}\mathrm{card}(\mathrm{extr}(\mathcal{F}_n))$.
	\Cref{lemma_th_faces} ensures that $m\leq \mathrm{card}(\mathrm{extr}(\mathcal{F}_0))$ and that, up to the extraction of a subsequence, there exists an injection $\theta_n:\mathrm{extr}(\mathcal{F}_n)\to
	\mathrm{extr}({\mathcal{F}_0})$ such that for
	every~${u=s\,\mathbf{1}_{F}/P(F)}$ in~$\mathrm{extr}(\mathcal{F}_n)$,
	we have $\theta_n(u)=s\,\mathbf{1}_E/P(E)$ with $F=E_{\varphi_u}$, and
	$$\lim_{n\to +\infty} \left(\max_{u \in \mathrm{extr}(\mathcal{F}_n)}\|\varphi_u\|_{\mathrm{C}^2(\partial E)}\right)= 0\,.$$
	In particular, for all $n$ large enough $\|\varphi_u\|_{\mathrm{C}^2(\partial E)}\leq \epsilon$ for all deformations $\varphi_u$, so that the conclusion of \Cref{th_faces} holds. We hence obtain a contradiction.
	\end{proof}

\subsection{Main result}
We are now able to introduce a non-degenerate version of the source condition,
which ultimately allows us to state our support recovery result.
\begin{definition}[Non-degenerate source condition]
	Let $u_0=\sum_{i=1}^N a_i\,\mathbf{1}_{E_i}$ be a simple function. We say
	that $u_0$ satisfies the non-degenerate source condition if
	\begin{enumerate}
		\item the source condition $\mathrm{Im}\,\Phi^* \cap \partial \TV(u_0)\neq \emptyset$ holds,
		\item for every $i\in\{1,...,N\}$, the set $E_i$ is a strictly stable
		solution to
		\customref{prescribed_curv}{$\mathcal{PC}(\mathrm{sign}(a_i)\eta_0)$},
		\item for every simple set $E\subset \RR^2$ s.t. $|E\triangle E_i|>0$
		for all $i\in\{1,...,N\}$, we have $\left|\int_{E}\eta_0\right|<P(E)$\,.
	\end{enumerate}
	In that case, we say that $\eta_0$ is non-degenerate.
	\label{nondegenerate_source_cond}
\end{definition}

\begin{theorem}
	Assume that $u_0=\sum_{i=1}^N a_i\mathbf{1}_{E_i}$ is a simple function
	satisfying the non-degenerate source condition, and that
	$\Phi_{\mathcal{F}_0}$ is injective. Then there exist constants
	$\alpha,\lambda_0\in\RR_+^*$ such that, for every $(\lambda,w)\in
	\RR_+^*\times\mathcal{H}$ with~$\lambda\leq \lambda_0$
	and~${\|w\|_{\mathcal{H}}/\lambda\leq \alpha}$, every
	solution~$u_{\lambda,w}$ of~\cref{primal_reg} is such that
	\begin{equation}
		u_{\lambda,w}=\sum\limits_{i=1}^N a_i^{\lambda,w}\,\mathbf{1}_{E_i^{\lambda,w}}
		\label{decomp_sol}
	\end{equation}
	with
	\begin{equation}
		\forall i\in\{1,...,N\},\,\left\{
		\begin{aligned}
			&\mathrm{sign}(a_i^{\lambda,w})=\mathrm{sign}(a_i)\\
			&E^{\lambda,w}_i=(E_i)_{\varphi_i^{\lambda,w}}~\mathrm{with}~\varphi^{\lambda,w}_i\in\mathrm{C}^2(\partial E_i)\,.\label{deform_support}
		\end{aligned}\right.
	\end{equation}
	Moreover,
\begin{align}
	\lim_{\substack{(\lambda,w)\to (0,0),\\
	0<\lambda\leq \lambda_0,\\ \|w\|_{\mathcal{H}}\leq \alpha \lambda }} a_{i}^{\lambda,w}=a_i, \quad \mbox{and} \quad	\lim_{\substack{(\lambda,w)\to (0,0),\\
	0<\lambda\leq \lambda_0,\\ \|w\|_{\mathcal{H}}\leq \alpha \lambda }} \big\|\varphi_i^{\lambda,w}\big\|_{\mathrm{C}^2(\partial E_i)} =0.
\end{align}
	\label{main_th}
\end{theorem}
\begin{proof}
	We fix $\delta>0$ small enough to have $(\partial E_i)^{\delta}\cap(\partial E_j)^{\delta}=\emptyset$ for every $i\neq j$, where~$${A^{\delta}\eqdef\cup_{x\in A}B(x,\delta)}\,.$$
	We also fix $\epsilon>0$ small enough to have $\epsilon<|a_i|P(E_i)$ for all $i\in\{1,...,N\}$
	and~$${(Id+\varphi\,\nu_{E_i})(\partial E_i)\subset (\partial E_i)^{\delta}}$$ as
	soon as~$\|\varphi\|_{\mathrm{C}^2(\partial E_i)}\leq \epsilon$. Finally, we
	take~$r>0$ such that the assumptions of \Cref{th_faces} hold.

	Our assumptions imply that $u_0$ is the unique solution to
	\cref{primal_noreg}. Hence, by~\Cref{conv_strict_bv}, we get that $\Diff
	u_{\lambda,w}$ and $|\Diff
	u_{\lambda,w}|$ respectively converge towards $\Diff u_0$ and $|\Diff u_0|$ in the weak-* topology  when $\lambda\to 0$ and
	$\|w\|_{\mathcal{H}}/\lambda \to 0$. 
	Since $|\Diff u_0|$ does not charge the boundary of the open set $(\partial E_i)^{\delta}$ for $1\leq i\leq N$,
	 there exist $\alpha>0$ and
	$\lambda_0>0$ such that for every~${(\lambda,w)\in\RR_+^*\times
	\mathcal{H}}$ with~${\lambda\leq \lambda_0}$
	and~${\|w\|_{\mathcal{H}}/\lambda\leq \alpha}$,
	\begin{equation}
		\forall i\in\{1,...,N\},\quad \left||\Diff u_{\lambda,w}|((\partial E_i)^{\delta})-|\Diff u_0|((\partial E_i)^{\delta})\right|\leq \epsilon.
		\label{bound_grad}
	\end{equation}
	Moreover, possibly reducing $\alpha$ and $\lambda_0$, we may also require that 
	$$\|\eta_{\lambda,w}-\eta_0\|_{\LD}+\|\eta_{\lambda,w}-\eta_0\|_{\mathrm{C}^1(\RR^2)}\leq
	r\,,$$ 
	which is possible by the continuity of $\Phi^*\colon\mathcal{H}\rightarrow \mathrm{C}^1(\RR^2)$, \Cref{conv_dual_vec,proj_dual_vec}.
	Let us fix such a pair $(\lambda,w)$, and write
	$$\mathrm{extr}(\mathcal{F}_{\lambda,w})=\left\{s^{\lambda,w}_i\frac{\mathbf{1}_{E_{i}^{\lambda,w}}}{P(E_{i}^{\lambda,w})}\right\}_{1\leq
	i\leq N_{\lambda,w}}$$
	where $\mathcal{F}_{\lambda,w}$ is the face of $\{\TV\leq 1\}$ exposed by $\eta_{\lambda,w}$.
	By \Cref{th_faces}, there exists an
	injective mapping
	$\theta_{\lambda,w}:\{1,...,N_{\lambda,w}\}\to\{1,...,N\}$ such that
	\begin{equation*}\forall i\in\{1,...,N_{\lambda,w}\},\quad\left\{\begin{aligned}
			&s_{i}^{\lambda,w}=s_{\theta_{\lambda,w}(i)}\,,\\
			&E_{i}^{\lambda,w}=(E_{\theta_{\lambda,w}(i)})_{\varphi_{i}^{\lambda,w}}~\mathrm{with}~\|\varphi_{i}^{\lambda,w}\|_{\mathrm{C}^2(\partial E_{\theta_{\lambda,w}(i)})}\leq \epsilon\,.
		\end{aligned}\right.\end{equation*}

	Let us show that $N_{\lambda,w}=N$. For all $j\in\{1,...,N\}$, since 
	$$|\Diff u_{\lambda,w}|((\partial E_j)^{\delta})\geq|\Diff u_0|((\partial
	E_j)^{\delta})-\epsilon=|a_j|P(E_j)-\epsilon>0\,,$$
	we note that $\mathrm{Supp}(\Diff
	u_{\lambda,w})\cap (\partial E_j)^{\delta}\neq \emptyset$. On the other hand, the sets $\{(\partial
	E_{j})^{\delta}\}_{j=1}^N$ are pairwise disjoint and
	$$\mathrm{Supp}(\Diff u_{\lambda,w})\subseteq
	\bigcup\limits_{i=1}^{N_{\lambda,w}}\partial
	E_i^{\lambda,w} \subseteq \bigcup\limits_{i=1}^{N_{\lambda,w}} (\partial
	E_{\theta_{\lambda,w}(i)})^{\delta}.$$
	Therefore, $\theta_{\lambda,w}$ must be surjective, and   
	 $N_{\lambda,w}=N$. Moreover, up to a permutation,
	$$\forall
	i\in\{1,...,N\},~s_i^{\lambda,w}=s_i~\mathrm{and}~E_i^{\lambda,w}=(E_i)_{\varphi_i^{\lambda,w}}\,.$$
	Now, the fact that $\big\|\varphi_i^{\lambda,w}\big\|_{\mathrm{C}^2(\partial F_i)} \to 0$ as $\left(\lambda, \norm{w}_\mathcal{H}/\lambda \right)\to (0,0)$ follows from \Cref{conv_result}.
	Moreover, since  $s_i^{\lambda,w}=s_i$ implies that 
	$\mathrm{sign}(a_i^{\lambda,w})=\mathrm{sign}(a_i)$, we have
	\begin{equation}
		\begin{aligned}
			||\Diff u_{\lambda,w}|((\partial E_i)^{\delta})-|\Diff u_0|((\partial E_i)^{\delta})|&=|a_i^{\lambda,w}P(E_i^{\lambda,w})-a_iP(E_i)|\,.
		\end{aligned}\label{massegraddelta}
	\end{equation}
	Finally, the weak-* convergence mentioned above implies that the left hand side of \eqref{massegraddelta} vanishes as~${\left(\lambda, \norm{w}_\mathcal{H}/\lambda \right)\to (0,0)}$. Since $P(E_i^{\lambda,w})\to P(E_i)>0$, we deduce that $a_i^{\lambda,w}\to a_i$.
	\end{proof}

\paragraph{Discussion on the non-degenerate source condition.} The non-degenerate source condition we introduce in \Cref{nondegenerate_source_cond} is the natural analog, for total (gradient) variation regularization, of the one introduced in \cite{duvalExactSupportRecovery2015}. As in the latter case, we expect its verification to be a delicate issue, which should be further investigated in future works. In particular, one could hope to prove that the condition is satisfied as soon as the unknown shapes are smooth enough and sufficiently separated from one another. In \Cref{sec_precertif}, we numerically investigate the validity of the non-degenerate source condition in a simplified setting, where $\Phi$ is a convolution operator and the unknown $u_0$ is a radial simple function.

\subsection{Verification of the non-degenerate source condition}
\label{sec_precertif}
Given $u$ an admissible function for~\cref{primal_noreg}, one may prove its optimality by finding some $p\in\mathcal{H}$ such that $\Phi^*p\in\partial\TV(u)$. We adopt here a strategy which is common in the literature on sparse recovery (see, e.g., \cite[Section 4]{duvalExactSupportRecovery2015} and references therein), which is to define a \emph{dual
pre-certificate}, that is, some ``good candidate'' $p\in\mathcal{H}$ for solving $\Phi^*p\in\partial\TV(u)$, usually defined by linearizing the dual problem. In
this subsection, we introduce the natural analog of the \emph{vanishing derivatives
pre-certificate} of \cite{duvalExactSupportRecovery2015}. We show that, if it is a valid dual certificate, then it is the one of minimal norm. In this case, it can hence be used to check the validity of the non-degenerate source condition. We numerically investigate the behaviour of this pre-certificate when the
unknown function is simple and radial. All the plots and experiments contained in this section can be reproduced using the code available online at~\url{https://github.com/rpetit/2023-support-recovery-tv}.

\subsubsection{The vanishing derivatives pre-certificate}

Let $N\in\NN^*$ and $u=\sum_{i=1}^N a_i\mathbf{1}_{E_i}$ be a $N$-simple
function with $a\in(\RR^*)^N$. If $p\in\mathcal{H}$, then $\Phi^*p$ is a dual
certificate associated to $u$ if and only if~${\Phi^*p\in\partial\TV(u)}$, that is, $\Phi^*p\in\partial\TV(0)$ and
\begin{equation*}
	\forall i\in\{1,...,N\},\quad E_i\in\underset{E\subset\RR^2,~|E|<+\infty}{\mathrm{Argmin}}~\left(P(E)-\mathrm{sign}(a_i)\int_{E}\Phi^*p\right)\,.
\end{equation*}
 The optimality conditions at order 0 and 1 respectively yield 
\begin{equation}
	\forall i\in\{1,...,N\},~\int_{E_i}\Phi^*p=\mathrm{sign}(a_i)\,P(E_i)~~\mathrm{and}~~
	\restriction{\Phi^*p}{\partial E_i}=\mathrm{sign}(a_i)\,H_{E_i}\,.
	\label{first_order_opt_dual_certif}
\end{equation}
We can then define a candidate dual certificate as the solution to
\cref{first_order_opt_dual_certif} with minimal norm.
\begin{definition}
	We call vanishing derivatives pre-certificate associated to some $N$-simple
	function~${u=\sum_{i=1}^N a_i\mathbf{1}_{E_i}}$ the
	function~$\eta_v=\Phi^*p_v$ with~$p_v$ the unique solution to
	\begin{equation}
		\underset{p \in \mathcal{H}}{\mathrm{min}}~\|p\|_{\mathcal{H}}^2~~\mathrm{s.t.}~~\forall i\in\{1,...,N\},~\int_{E_i}\Phi^*p=\mathrm{sign}(a_i)\,P(E_i)~~\mathrm{and}~~
		\restriction{\Phi^*p}{\partial E_i}=\mathrm{sign}(a_i)\,H_{E_i}\,.
		\label{vanishing_derivatives}
	\end{equation}
\end{definition}
The admissible set of \Cref{vanishing_derivatives} is weakly closed. Hence, if
the source condition holds (i.e. there exists $\eta\in\mathrm{Im}\,\Phi^*$ such
that $\eta\in\partial\TV(u)$), \Cref{vanishing_derivatives} is feasible
and~$p_v$ is therefore well-defined.

Since any dual certificate satisfies \cref{first_order_opt_dual_certif}, we have
the following result.
\begin{proposition}
	If \Cref{vanishing_derivatives} is feasible and $\eta_v\in\partial\TV(0)$,
	then $\eta_v$ is the minimal norm dual certificate, i.e. $\eta_v=\eta_0$.
\end{proposition}

\subsubsection{Deconvolution of radial simple functions}
\label{deconv_radial}
We now focus on the case where $\mathcal{H}=\LD$ and $\Phi=h\star\cdot$ is the
convolution with the Gaussian kernel~$h$ with variance $\sigma$, and
$E_i=\mathbf{1}_{B(0,R_i)}$ for all $i\in\{1,...,N\}$ with $0<R_1<...<R_N$. Let us
introduce the following mappings
\begin{equation*}
	\begin{aligned}
		\Phi_E : \RR^{N} &  \to \mathcal{H} & \Phi'_E : \RR^{N} &  \to \mathcal{H} & \Gamma_E : \RR^{2N} &  \to \mathcal{H}\\
		a & \mapsto \sum\limits_{i=1}^N a_i\,h\star \mathbf{1}_{E_i}\,, & b & \mapsto \sum\limits_{i=1}^N b_i\,h\star (\mathcal{H}^1\mres\partial E_i)\,, & \begin{pmatrix}a\\b\end{pmatrix} & \mapsto \Phi_E a+\Phi'_E b\,.
	\end{aligned}
\end{equation*}
With these notations, we can show the following.
\begin{lemma}
	If \Cref{vanishing_derivatives} is feasible, $p_v$ is radial\footnote{We say
	that a function $f\in\LD$ is radial if there exists $g:\RR_+\to \RR$ such
	that $f(x)=g(\|x\|)$ for almost every $x\in\RR^2$.} and is the unique
	solution of
		\begin{equation}
		\underset{p \in \mathcal{H}}{\mathrm{min}}~\|p\|_{\mathcal{H}}^2~~\mathrm{s.t.}~~\Gamma_E^* p=\begin{pmatrix}(\mathrm{sign}(a_i)P(E_i))_{1\leq i \leq N}\\(\mathrm{sign}(a_i)2\pi)_{1\leq i\leq N}\end{pmatrix}.
		\label{vanishing_derivatives_bis}
	\end{equation}
	\label{lemma_radial_eta_v}
\end{lemma}
To prove this, we introduce the radialization $\tilde{p}$ of any function
$p\in\LD$, defined by:
\begin{equation}
	\mbox{for a.e. $x\in\RR^2$},~\tilde{p}(x)=\frac{1}{2\pi}\int_{\mathbb{S}^1}p(\|x\|e)\,d\mathcal{H}^1(e)\,.\label{radialization}
\end{equation}
The radialization operator is self-adjoint and $\|\tilde{u}\|_{\LD}\leq
\|u\|_{\LD}$ for every $u\in\LD$.
\begin{proof}[Proof of \Cref{lemma_radial_eta_v}]
	Let us show that, if $p$ is admissible for~\Cref{vanishing_derivatives},
	then so is $\tilde{p}$. Using the fact that~$\|\tilde{p}\|_{\LD}\leq
	\|p\|_{\LD}$ and the uniqueness of the solution to~\eqref{vanishing_derivatives}, this will conclude that (if it exists) $p_v$ is radial. 
	
	Note that the radialization~\eqref{radialization}  can be seen a Bochner integral in $\LD$:
\begin{align*}
\forall p \in \LD,\quad	\tilde{p}=\frac{1}{2\pi}\int_{\mathbb{S}^1} (p\circ {R}_e)\mathrm{d}\mathcal{H}^1(e),
\end{align*}
where $R_e$ is the rotation which maps $(0,1)$ to $e$, and the integral is well-defined since $e\mapsto p\circ R_e$ is continuous from $\mathbb{S}^1$ to $\LD$. Moreover,  since $h$ is radial, for all $e \in \mathbb{S}^1$ and $x \in \mathbb{R}^2$,
\begin{align*}
	\left\langle p\circ R_e, \varphi(x)  \right\rangle = \int_{\mathbb{R}^2}p\circ R_e(x-t)h(t)\mathrm{d} t = \int_{\mathbb{R}^2}p(R_e(x)-t')h(t')\mathrm{d} t' = \left\langle p,\varphi(R_e(x))\right\rangle.
\end{align*}
As a result, if $p$ is admissible for~\Cref{vanishing_derivatives}, since the sets $(E_i)_{1\leq i\leq N}$ are radial, we get 
\begin{align*}
	\int_{E_i} \langle \tilde{p},\varphi(x)\rangle \mathrm{d}x &= \frac{1}{2\pi}\int_{\mathbb{S}^1} \left(\int_{E_i}  \langle {p},\varphi(R_e(x))\rangle \mathrm{d}x\right)\mathrm{d}\mathcal{H}^1(e) = \int_{E_i}\langle {p},\varphi(x)\rangle \mathrm{d}x=\mathrm{sign}(a_i)\,P(E_i), \\
	\intertext{and for all $x \in \partial E_i,$}
		\left\langle \tilde{p}, \varphi(x)  \right\rangle &= \frac{1}{2\pi}\int_{\mathbb{S}^1}\left\langle(p, \varphi(R_e(x))\right\rangle\mathrm{d}\mathcal{H}^1(e) = \frac{1}{R_i}.  
\end{align*}
Hence $\tilde{p}$ is admissible too. 
The reformulation~\eqref{vanishing_derivatives_bis} follows from the fact that the convolution with $h$ is self-adjoint.
\end{proof}

\begin{proposition}\label{proppseudoinv}
	The operator $\Gamma_E$ is injective. Moreover, if \Cref{vanishing_derivatives_bis} is feasible, then
	\begin{equation*}
		\eta_v=\Phi^*\Gamma^{+,*}_E\begin{pmatrix}(\mathrm{sign}(a_i)P(E_i))_{1\leq i \leq N}\\(\mathrm{sign}(a_i)2\pi)_{1\leq i\leq N}\end{pmatrix}.
	\end{equation*}
	where $\Gamma_E^{+,*}=\Gamma_E(\Gamma_E^*\Gamma_E)^{-1}$.
\end{proposition}

\begin{proof}
	First, we prove that $\Gamma_E$ is injective. 
	Let $(a,b)\in \RR^{2N}$ be such that $\Phi_E a+\Phi'_E b=0$. We get
	that~$h\star(\sum_{i=1}^N a_i \mathbf{1}_{E_i}+b_i \mathcal{H}^1\mres \partial
	E_i)=0$, which, using the injectivity of $h\star\cdot$, yields $$\sum_{i=1}^N
	a_i\mathbf{1}_{E_i}+b_i\mathcal{H}^1\mres \partial E_i=0\,.$$ Integrating both
	sides of this equality against a test function compactly supported in the open set~$B(0,R_N)\setminus \overline{B(0,R_{N-1})}$ shows that~$a_N=0$. Apply this argument
	repeatedly also allows to obtain $a_1=...=a_N=0$. Then, since the measures
	$(\mathcal{H}^1\mres\partial E_i)_{1\leq i\leq N}$ have disjoint support, we
	obtain $b_1=...=b_N=0$.
	
	Now, \Cref{vanishing_derivatives_bis} reformulates $p_v$ as the least-norm solution of a linear system, therefore 
\begin{align*}
	p_v=(\Gamma_E^{*})^+\begin{pmatrix}(\mathrm{sign}(a_i)P(E_i))_{1\leq i \leq N}\\(\mathrm{sign}(a_i)2\pi)_{1\leq i\leq N}\end{pmatrix}
\end{align*}
where $(\Gamma_E^{*})^+$ is the Moore-Penrose pseudoinverse (see \cite{engl1996regularization}) of the (closed-range) operator $\Gamma_E^{*}\colon \mathcal{H}\rightarrow \mathbb{R}^{2N}$.
Since $\Gamma_E$ is injective (hence $\Gamma_E^*$ is surjective), it is standard that the normal equations imply that $(\Gamma_E^{*})^+=\Gamma_E(\Gamma_E^*\Gamma_E)^{-1} =(\Gamma_E^{+})^*$. 
\end{proof}

\begin{figure}
	\centering
	\includegraphics[width=.99\textwidth]{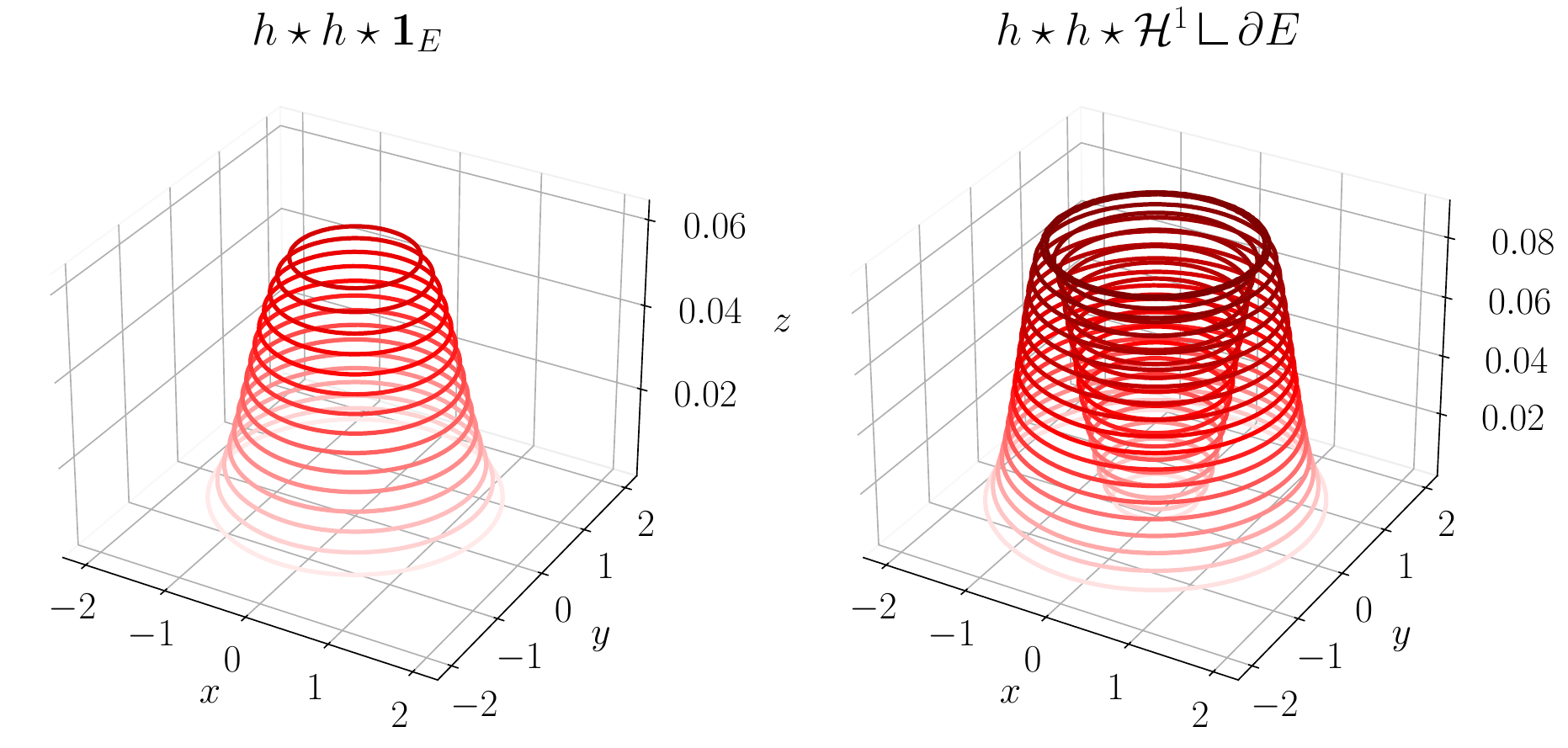}
	\caption{Plots of $h\star h\star \mathbf{1}_E$ and $h\star
	h\star(\mathcal{H}^1\mres\partial E)$ for $E=\mathbf{1}_{B(0,1)}$ and $h$
	the Gaussian kernel with variance $\sigma=0.2$.}
	\label{basis_func_vanishing}
\end{figure}

\Cref{proppseudoinv} asserts that there exist Lagrange multipliers $(a,b)\in\RR^{2N}$ such that
$$p_v=\sum\limits_{i=1}^N \left(a_i\,h\star\mathbf{1}_{E_i}+b_i\,
h\star(\mathcal{H}^1\mres\partial E_i)\right)\,.$$ 
We provide in
\Cref{basis_func_vanishing} a plot of~${\Phi^*(h\star\mathbf{1}_E)}=h\star
h\star \mathbf{1}_E$ and~${\Phi^*(\mathcal{H}^1\mres\partial E)}=h\star
h\star(\mathcal{H}^1\mres\partial E)$ for $E=B(0,1)$, which are the two~``basis
functions'' from which $\eta_v$ is built.

\paragraph{Ensuring $\bm{\eta_v}$ is a valid dual certificate.} From
\Cref{charac_subdiff_z}, we know that to show $\eta_v\in\partial\TV(0)$, it is
sufficient to find~${z\in\mathrm{L}^{\infty}(\RR^2,\RR^2)}$ such that
$\|z\|_{\infty}\leq 1$ and $\mathrm{div}\,z=\eta_v$. Since~$p_v$ is radial, so
is~$\eta_v$. It is hence natural to look for a radial vector field $z$ (i.e.
such that there exists $z_r:\RR_+\to \RR$ with~${z(x)=z_r(\|x\|)\,x/\|x\|}$ for
almost every $x\in\RR^2$). In this case we have $\mathrm{div}\,z=\eta_v$ if and
only if, for every $r>0$:
\begin{equation*}
	\begin{aligned}
		\tilde{\eta}_v(r)=\frac{1}{r}\frac{\partial}{\partial r}(r\,z_r)(r)&\iff r\,\tilde{\eta}_v(r)=\frac{\partial}{\partial r}(r\,z_r)(r)\\
		&\iff z_r(r)=\frac{1}{r}\int_{0}^r \tilde{\eta}_v(s)\,s\,ds\,,
	\end{aligned}
\end{equation*}
where, abusing notation, we have denoted by $\tilde{\eta}_v(r)$ the value of
$\eta_v(x)$ for any $x$ such that $\|x\|=r$. Thus, one only needs
to ensure that the mapping $f_v$ defined by
\begin{equation}
	\begin{aligned}[t]
		f_v : \RR_+ &  \to \RR \\
		r  & \mapsto \frac{1}{r}\int_{0}^r \tilde{\eta}_v(s)\,s\,ds
	\end{aligned}
	\label{f_v}
\end{equation}
satisfies $\|f_v\|_{\infty}\leq 1$ to show $\eta_v\in\partial\TV(0)$. 
\begin{remark}
	Looking for a radial vector field is not
	restrictive. In fact, if a vector field $z$ is suitable, then so is the
	radial vector field $\tilde{z}$ defined by
	$$\tilde{z}(x)\eqdef\tilde{z}_r(\|x\|)\frac{x}{\|x\|}~~~\mathrm{with}~~~\tilde{z}_r(r)\eqdef\frac{1}{2\pi}\int_{0}^{2\pi}z_r(r,\theta)\,d\theta\,,$$
	where $z_r$ denotes the radial component of $z$. Indeed, we have
	$|\tilde{z}(r)|\leq 1$ for all $r$ with equality if and only
	if~${z_r(r,\theta)=1}$ for almost every $\theta$ or $z_r(r,\theta)=-1$ for
	almost every $\theta$. Moreover
	\begin{equation*}
		\begin{aligned}
			\eta_v(r)=\frac{1}{2\pi}\int_{0}^{2\pi}\eta_v(r)\,dr&=\frac{1}{r}\frac{\partial}{\partial r}\left(r\,\frac{1}{2\pi}\int_{0}^{2\pi}z_r(r,\theta)\,d\theta\right)+\frac{1}{r}\,\frac{1}{2\pi}\int_{0}^{2\pi}\frac{\partial z_{\theta}}{\partial\theta}(r,\theta)\,d\theta\\
			&=\frac{1}{r}\frac{\partial}{\partial r}\left(r\,\tilde{z}_r\right)=\mathrm{div}\,\tilde{z}.
		\end{aligned}
	\end{equation*}
\end{remark}
\paragraph{Verification of the non-degenerate source condition.} Finally, we can
investigate the validity of the non-degenerate source condition. In this
setting, it holds if and only if the following three conditions are
simultaneously satisfied:
\begin{equation}
	\begin{aligned}
	&\forall R\in\RR_+\setminus \{R_1,...,R_N\},~&&|f_v(R)|<1\,,\\
	&\forall i\in\{1,...,N\},~&&f_v(R_i)=\mathrm{sign}(a_i)\,,\\
	&\forall i\in\{1,...,N\},~&&E_i\text{ is a strictly stable solution to }\customref{prescribed_curv}{$\mathcal{PC}(\mathrm{sign}(a_i)\eta_v)$}\,.
	\end{aligned}
\end{equation}
As explained in \Cref{sec_suff_stab}, the last property holds provided that
\begin{equation*}
	\forall i\in\{1,...,N\},~-\mathrm{sign}(a_i)\underset{x\in\partial E_i}{\mathrm{sup}}\left[H_{E_i}^2(x)+\frac{\partial\eta_v}{\partial\nu_{E_i}}(x)\right]>0\,.
\end{equation*}
In our case $H_{E_i}$ is constant equal to $1/R_i$, and, since $\eta_v$ is
radial, $\frac{\partial\eta_v}{\partial\nu_{E_i}}$ is constant on $\partial
E_i$. Proving that
\begin{equation}
	\forall i\in\{1,...,N\},~-\mathrm{sign}(a_i)\left[\frac{1}{R_i^2}+\frac{\partial\eta_v}{\partial r}(R_i)\right]>0
	\label{quantity_second_deriv}
\end{equation}
is hence sufficient. Moreover, a direct computation also shows that, if
\Cref{vanishing_derivatives_bis} is feasible, then
\begin{equation*}
	\forall i\in\{1,...,N\},~f_v''(R_i)=\frac{1}{R_i^2}+\frac{\partial \eta_v}{\partial r}(R_i)\,,
\end{equation*}
so that \cref{quantity_second_deriv} can be directly checked by looking at the
graph of $f_v$.

\begin{figure}
	\centering
	\includegraphics[width=.99\textwidth]{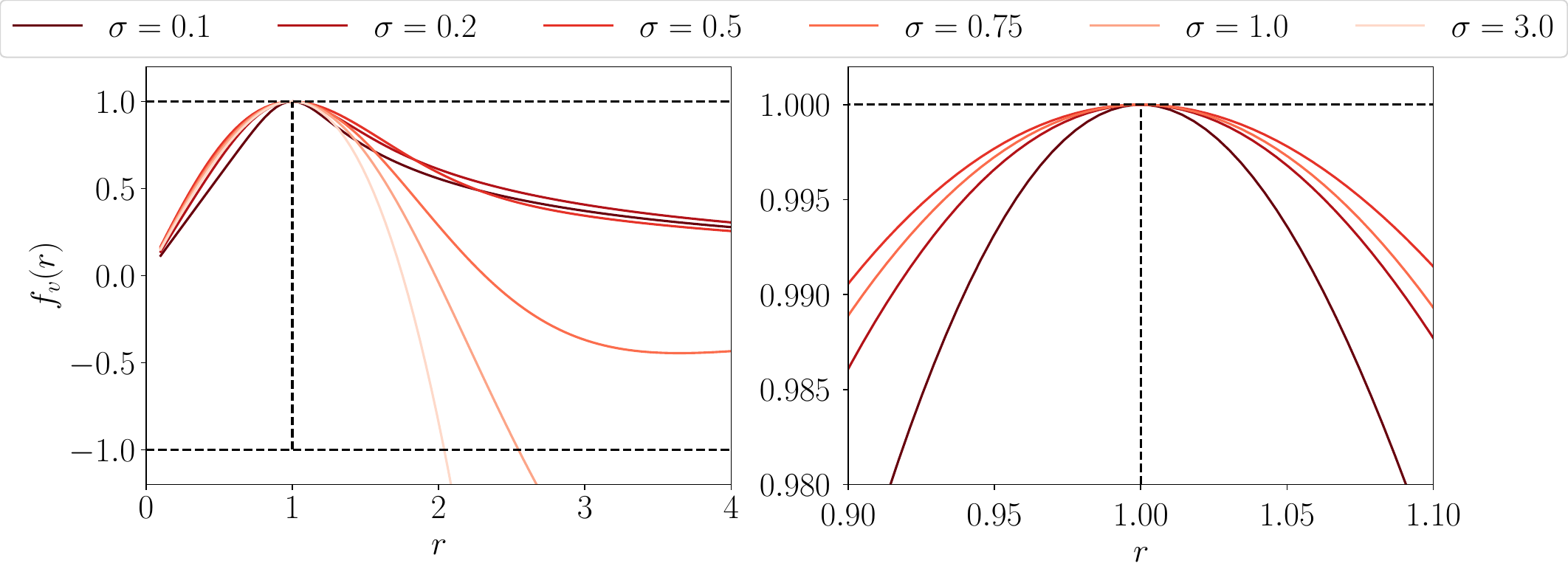}
	\caption{Graph of $f_v$ defined in \cref{f_v} when $N=1$, $R_1=1$ and
	$\mathrm{sign}(a_1)=1$ (left: global graph, right: zoom around $r=1$).}
	\label{pre_certif}
\end{figure}

\begin{figure}
	\centering
	\includegraphics[width=.7\textwidth]{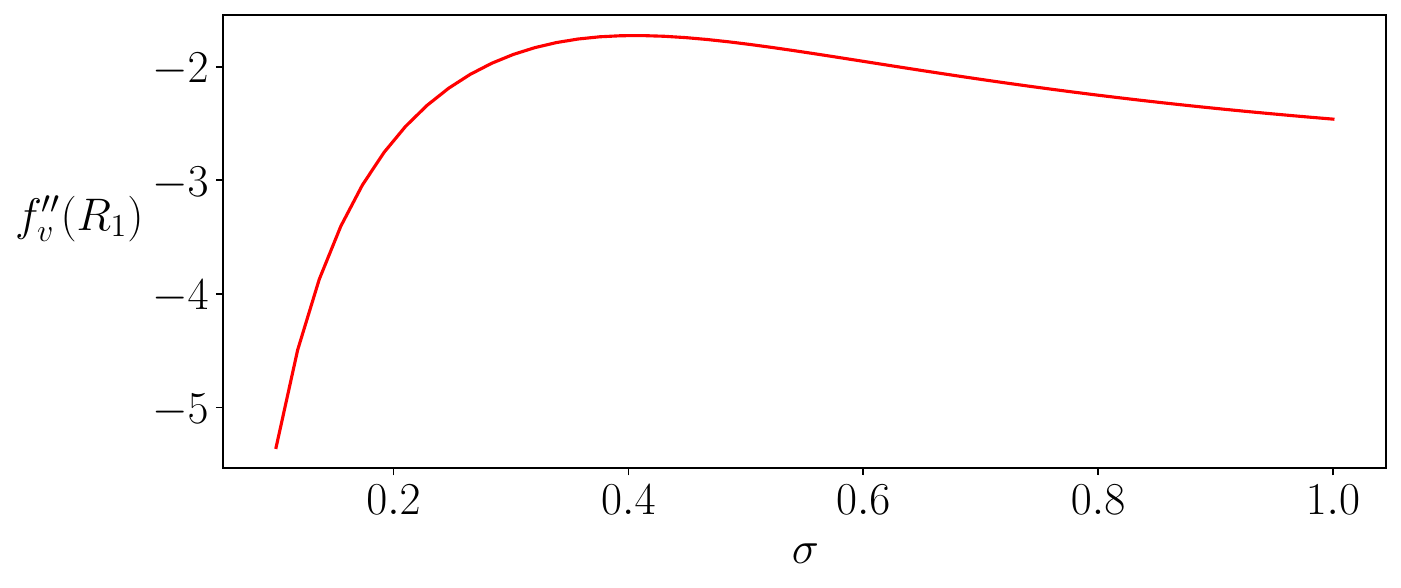}
	\caption{Graph of $f_v''(R_1)$ as a function of $\sigma$ when $N=1$, $R_1=1$
	and $\mathrm{sign}(a_1)=1$.}
	\label{normal_deriv}
\end{figure}

\paragraph{Numerical experiment ($\bm{N=1}$).} Here, we investigate the case
where $N=1$, $R_1=1$ and~${\mathrm{sign}(a_1)=1}$. \Cref{pre_certif} shows the
graph of $f_v$ for several values of~$\sigma$. This suggests that there
exists~$\sigma_0>0$ such that $\eta_v$ is a dual certificate~(and hence the one
with minimal norm) for every~$\sigma<\sigma_0$. It even seems
that~${\sigma_0\geq 0.75}$.
In \Cref{normal_deriv}, we numerically compute $f_v''(R_1)$ and notice it is
(strictly) negative, even when~${\eta_v\notin\partial\TV(0)}$. This suggests that
there exists~${\sigma_0>0}$ such that, for every~${\sigma\leq\sigma_0}$, the
non-degenerate source condition holds~(and, from our experiments, it seems
that~${\sigma_0\geq 0.75}$). Surprisingly, $\sigma\mapsto f_v''(R_1)$ does not
seem to be monotonous, even on~${[0,\sigma_0)}$.

\begin{figure}
	\centering
	\includegraphics[width=.99\textwidth]{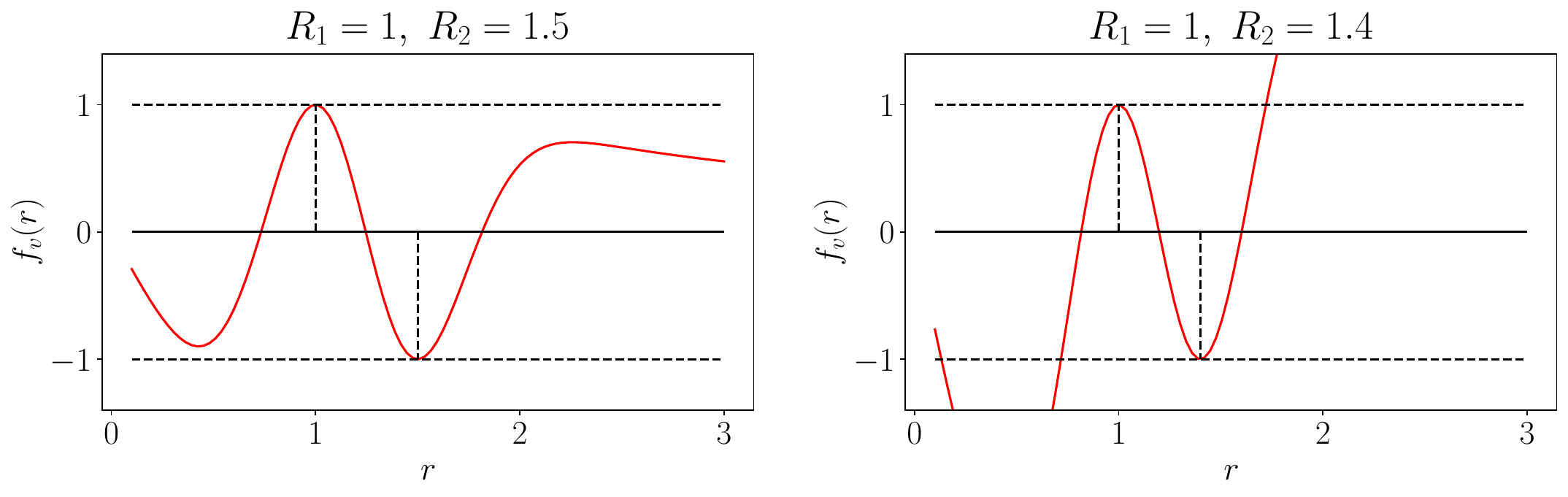}
	\caption{Graph of $f_v$ defined in \cref{f_v} when $\sigma=0.2$, $N=2$,
	$\mathrm{sign}(a_1)=-\mathrm{sign}(a_2)$, $R_1=1$, ${R_2=1.5}$~(left) and
	$R_2=1.4$ (right).}
	\label{pre_certif_opposite_signs}
\end{figure}

\begin{figure}
	\centering
	\includegraphics[width=.99\textwidth]{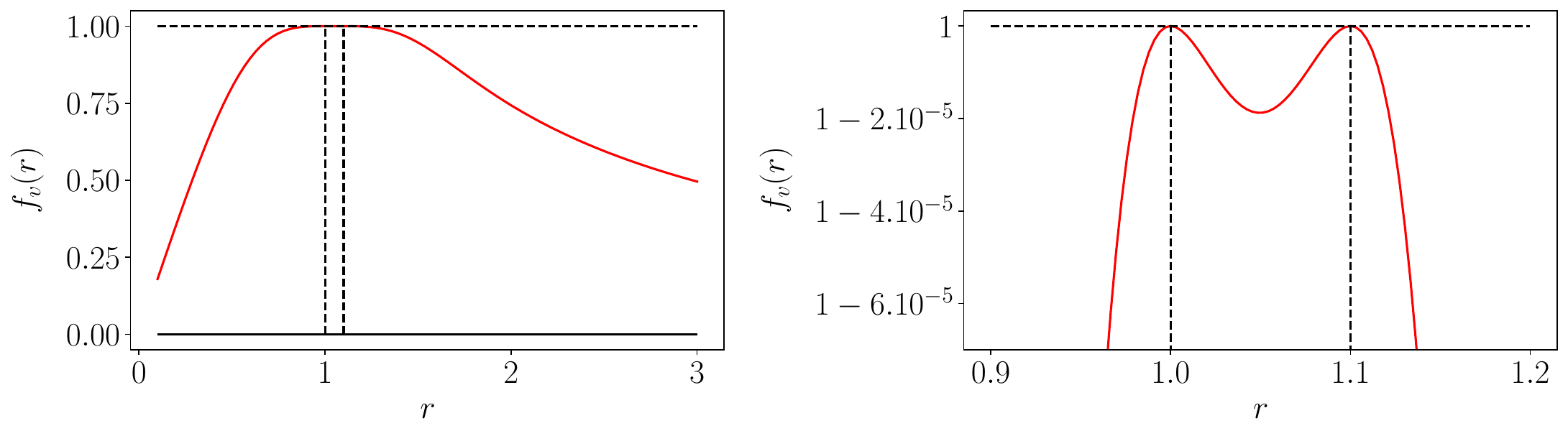}
	\caption{Graph of $f_v$ defined in \cref{f_v} when $\sigma=0.2$, $N=2$,
	$\mathrm{sign}(a_1)=\mathrm{sign}(a_2)$, $R_1=1$ and~${R_2=1.1}$ (left: global graph, right: zoom around $1$).}
	\label{pre_certif_equal_signs}
\end{figure}

\paragraph{Numerical experiments ($\bm{N\geq 2}$).} Now, we investigate the case
where $N=2$. Our experiments suggest the existence of two completely different
regimes. If $\mathrm{sign}(a_1)\neq \mathrm{sign}(a_2)$, then~$\eta_v$ is
non-degenerate only if $R_1$ and $R_2$ are not too close (see
\Cref{pre_certif_opposite_signs}). On the contrary, the case
where~$\mathrm{sign}(a_1)=\mathrm{sign}(a_2)$ seems to correspond to a real
super-resolution regime, as $\eta_v$ is non-degenerate even for arbitrarily
close $R_1$ and $R_2$ (see \Cref{pre_certif_equal_signs}). Still, we notice
that, in this last case, the quantities $f_v''(R_1)$ and $f_v''(R_2)$, which
control the stability of the recovery, go to $0$ as $R_1$ and $R_2$ get closer.

\paragraph{Beyond the radial case.}  In the general case, to numerically ensure that
$\eta_v\in\partial\TV(0)$, one can solve
\begin{equation}
	\underset{u\in\{\TV\leq 1\}}{\mathrm{sup}}~\int_{\RR^2}\eta_v u\,,
	\label{prob_max_face}
\end{equation}
which can be done by relying on standard discretization techniques.
Indeed, as underlined in~\Cref{subgrad_exposed}, we have that
$\eta_v\in\partial\TV(0)$ if and only if \cref{prob_max_face} has a value which
is no greater than~$1$. To ensure the non-degenerate source condition holds, one 
must also show that~${|\int_{E}\eta_v|<P(E)}$
for every simple set $E$ such that~$|E\triangle E_i|>0$ for every $i\in\{1,...,N\}$. This 
last property holds if and only if $\mathrm{Supp}(\Diff u)\subseteq \bigcup_{i=1}^N \partial E_i$ 
for every solution $u$ of \cref{prob_max_face}. It can therefore be checked by finding, among all 
solutions of \cref{prob_max_face}, the one such that $\mathrm{Supp}(\Diff u)$ is maximal. Finally, one should also verify the strict stability of $E_i$ as a solution to \customref{prescribed_curv}{$\mathcal{PC}(\mathrm{sign}(a_i)\eta_v)$} for every $i\in\{1,...,N\}$.

%% file: appendix.tex
\section{Proofs of \texorpdfstring{\Cref{sec_curv_pb}}{Section 4}}
\label{appendix_proof}
\subsection{Proof of \texorpdfstring{\Cref{conv_result}}{Proposition 4.1}}
We argue by contradiction and assume the existence of two sequences $(\eta_n)_{n\in\NN^*}$ and $(F_n)_{n\in\NN^*}$ such that
\begin{itemize}
	\item for all $n\in\NN^*$, $\eta_n\in\partial\TV(0)\cap\mathrm{C}^1(\RR^2)$,
	\item the sequence $(\eta_n)_{n\in\NN^*}$ converges to $\eta_0$ in $\LD$ and $\mathrm{C}^1(\RR^2)$,
	\item for all $n\in\NN^*$, the set $F_n$ is a non-empty solution to \customref{prescribed_curv}{$\mathcal{PC}(\eta_n)$} and cannot be written as a~$\mathrm{C}^2$-normal deformation of size at most $\epsilon$ of a non-empty solution to \customref{prescribed_curv}{$\mathcal{PC}(\eta_0)$}\,.
\end{itemize}
We hence have that $(F_n)_{n\in\NN^*}$ is bounded (by \Cref{lemma_bounds}) and that $F_n$ is a strong $(\Lambda,r_0)$-quasi minimizer of the perimeter~(in short
form~${F_n\in\mathcal{QM}(\Lambda,r_0)}$, see \cite[Section 21]{maggiSetsFinitePerimeter2012} and~\cite[Definition 4.7.3]{ambrosioCorsoIntroduttivoAlla2010}) for all~${n\in\NN^*}$, with $\Lambda=\mathrm{sup}\,\{\|\eta_n\|_{\infty},\,n\in\NN^*\}$ and~$r_0$ any
positive real number. Taking $r_0$ small enough to have $\Lambda\,r_0\leq
1$, from \cite[Propositions~21.13 and 21.14]{maggiSetsFinitePerimeter2012}
we get (up to the extraction of a not relabeled subsequence) that $(F_n)_{n\in\NN^*}$ converges in measure to a bounded set $E\in\mathcal{QM}(\Lambda,r_0)$, and that $(\partial F_n)_{n\in\NN^*}$
converges to $\partial E$ in the Hausdorff sense. From $|F_n\triangle E|\to 0$ we obtain that $E$ is a solution to \customref{prescribed_curv}{$\mathcal{PC}(\eta_0)$}. In addition, since~$F_n$ is non-empty for all $n$, using \Cref{lemma_bounds}, we get that $E$ is non-empty. The convergence of~$(\partial F_n)_{n\in\NN^*}$ towards $\partial E$ also yields
$$\forall r>0,~\exists n_0\in\NN,~\forall n\geq n_0,~\partial F_n\subset \bigcup\limits_{x\in\partial E}C(x,r,\nu_E(x))\,,$$
where $C(x,r,\nu_E(x))$ denotes the square of axis $\nu_E(x)$ and side $r$ centered at $x$, defined in \Cref{def_cylinder}.
From \cite[4.7.4]{ambrosioCorsoIntroduttivoAlla2010}, and arguing as in the proof of \cite[Theorem 26.6]{maggiSetsFinitePerimeter2012}, for every $x\in\partial E$ we obtain the existence
of $r>0$, of $n_0\in\NN$, of $u\in\mathrm{C}^{1,1}([-r,r])$ and of a sequence~${(u_n)_{n\geq n_0}}$ which is uniformly bounded in
$\mathrm{C}^{1,1}([-r,r])$, such that, in $C(x,r,\nu_E(x))$, the set $E$ is the hypograph of $u$ and, for every~${n\geq n_0}$, the set $F_n$ is the hypograph of $u_n$. Moreover, we have that $\|u_n-u\|_{\mathrm{C}^1([-r,r])}\to 0$.

Now, we also have that $u$ and $u_n$ (for $n\geq n_0$) respectively solve (in the sense of distributions) the following equations in $(-r,r)$:
\begin{equation}
	\begin{aligned}
		\frac{u''(z)}{\left(1+u'(z)^2\right)^{3/2}}&=H(z,u(z))\,,&\mathrm{with}~H(z,t)&\eqdef\eta_0(x+R_{\nu_E(x)}(z,t))\,, \\
		\frac{u_n''(z)}{\left(1+u_n'(z)^2\right)^{3/2}}&=H_n(z,u_n(z))\,,&\mathrm{with}~H_n(z,t)&\eqdef\eta_n(x+R_{\nu_E(x)}(z,t))\,.
	\end{aligned}
	\label{opt_pc}
\end{equation}
We hence immediately obtain that $u$ and $u_n$ belong to $\mathrm{C}^2([-r,r])$. Moreover, for every $z\in(-r,r)$ we have:
\begin{equation*}
	\begin{aligned}
		|u''_n(z)-u''(z)|&= \Big|H_n(z,u_n(z))\left(1+u'_n(z)^2\right)^{3/2}-H(z,u(z))\left(1+u'(z)^2\right)^{3/2}\Big|\\
		&\leq ~~~~(\|H_n-H\|_{\infty}+|H(z,u_n(z))-H(z,u(z))|)\left(1+u'_n(z)^2\right)^{3/2}\\&~~~~\,+\|H\|_{\infty}\left[\left(1+u'_n(z)^2\right)^{3/2}-\left(1+u'(z)^2\right)^{3/2}\right],
	\end{aligned}
\end{equation*}
from which we obtain that $\|u''_n-u''\|_{\infty}\to 0$.

Using these new results in combination with \cref{opt_pc}, we get that $u$ and $u_n$ belong to $\mathrm{C}^3([-r,r])$. Differentiating \cref{opt_pc}, we obtain, for every $z\in(-r,r)$:
\begin{equation*}
	\begin{aligned}
		u^{(3)}(z)&=~~~~\big[\partial_1 H(z,u(z))+u'(z)\,\partial_2 H(z,u(z))\big]\,(1+u'(z)^2)^{3/2}\\
		&~~~~\,+3\,H(z,u(z))\,u''(z)\,u'(z)\,(1+u'(z)^2)^{3/2}\,,\\
		u_n^{(3)}(z)&=~~~~\big[\partial_1 H_n(z,u_n(z))+u'_n(z)\,\partial_2 H_n(z,u_n(z))\big]\,(1+u'_n(z)^2)^{3/2}\\
		&~~~~\,+3\,H_n(z,u_n(z))\,u''_n(z)\,u'_n(z)\,(1+u'_n(z)^2)^{3/2}\,,
	\end{aligned}
\end{equation*}
from which we can finally show $\|u^{(3)}_n-u^{(3)}\|_{\infty}\to 0$.

Finally, using the compactness of $\partial E$, we obtain that $(F_n)_{n\geq 0}$ converges in $\mathrm{C}^3$ towards $E$, and \Cref{prop_conv_normal} allows to eventually write $F_n$ as a $\mathrm{C}^2$-normal deformation of~$E$, whose norm converges to zero. This yields a contradiction.

\subsection{Proofs of \texorpdfstring{\Cref{sec_stability}}{Section 4.1}}
\label{proof_stab_prob_pc}
To prove \Cref{lemma_cont_shape_hessian_3,lemma_conv_shape_hessian}, we need to compute~$\mathrm{j}''_{E}(\varphi)$ for $\varphi\in\mathrm{C}^1(\partial E)$ in a neighborhood of $0$. This may be done using \Cref{lemma_formula_shape_hessian} below. To state it, given a bounded set $E$ of class~$\mathrm{C}^2$ and $\varphi$ in a neighborhood of $0$
in $\mathrm{C}^1(\partial E)$, we introduce the mapping $f_\varphi=Id+\xi_\varphi$, with~$\xi_\varphi$ defined as in \Cref{lemma_extension}. If $\|\varphi\|_{\mathrm{C}^1(\partial E)}$ is
sufficiently small then $f_\varphi$ is a $\mathrm{C}^1$-diffeomorphism, and we denote its inverse by $g_\varphi$.
\begin{lemma}
	Let $E$ be a bounded set of class $\mathrm{C}^2$. Then for
	every $\varphi$ in a neighborhood of $0$ in~$\mathrm{C}^1(\partial E)$, and for every~${\psi\in\mathrm{H}^1(\partial E)}$, we have:
	\begin{equation}
		\mathrm{j}''_E(\varphi).(\psi,\psi)=j''_{E_\varphi}(0).(\xi_{\psi}\circ g_{\varphi}\cdot\nu_\varphi,\,\xi_{\psi}\circ g_{\varphi}\cdot\nu_\varphi)+\mathrm{j}'_{E_\varphi}(0).(Z_{\varphi,\psi})
		\label{formula_shape_hessian}
	\end{equation}
	where $\nu_\varphi$ is the unit outward normal to $E_\varphi$ and
	$$Z_{\varphi,\psi}=B_{\varphi}((\xi_\psi\circ g_\varphi)_{\tau_\varphi},(\xi_\psi\circ g_\varphi)_{\tau_\varphi})-2(\nabla_{\tau_\varphi}(\xi_\psi\circ g_\varphi\cdot\nu_{\varphi}))\cdot(\xi_\psi\circ g_\varphi)_{\tau_\varphi}\,,$$
	with $\zeta_{\tau_\varphi}$ and $\nabla_{\tau_\varphi}\zeta$ the tangential part and the tangential gradient of $\zeta$ with respect to $E_\varphi$, and~$B_{\varphi}$ the second fundamental form of $E_\varphi$.
	\label{lemma_formula_shape_hessian}
\end{lemma}
\begin{proof}
	To prove this result, we need to introduce $\mathcal{J}_E$\footnote{This mapping allows to study the behaviour of the objective in a neighborhood of $E$ with respect to general deformations, while $\mathrm{j}_E$ is only related to normal deformations.} defined by
	$$\! \begin{aligned}[t]
		\mathcal{J}_E : \mathrm{C}^1_b(\RR^2,\RR^2) &  \to \RR \\
		\xi   & \mapsto J((Id+\xi)(E))\,.
	\end{aligned}$$
	We denote by $\nu$ the outward unit normal to $E$ and $B$ its second fundamental form. We also denote~$\zeta_\tau$ and~$\nabla_\tau\zeta$ the tangential part and the tangential gradient of $\zeta$ with respect to $E$. The structure theorem~(see~e.g.~\cite[Theorem 5.9.2]{henrotShapeVariationOptimization2018} or \cite[Theorem 2.1]{dambrineStabilityShapeOptimization2019}) then yields, for every sufficiently smooth vector field $\zeta$:
	\begin{equation*}
		\begin{aligned}
			\mathcal{J}'_E(0).(\zeta)&=\mathrm{j}_E'(0).(\restriction{\zeta}{\partial E}\cdot\nu)\,,\\
			\mathcal{J}''_E(0).(\zeta,\zeta)&=\mathrm{j}_E''(0).(\restriction{\zeta}{\partial E}\cdot\nu,\restriction{\zeta}{\partial E}\cdot\nu)+\mathrm{j}_E'(0).(Z_\zeta)\,,
		\end{aligned}
	\end{equation*}
	where
	$$Z_\zeta\eqdef B(\zeta_\tau,\zeta_\tau)-2\,(\nabla_\tau(\zeta\cdot\nu))\cdot \zeta_\tau\,.$$
	Now, we first notice that, for every pair of vector fields $\xi,\zeta$ such that $Id+\xi$ is invertible, we have:
	\begin{equation*}
		(Id+\xi+\zeta)(E)=(Id+\zeta\circ(Id+\xi)^{-1})((Id+\xi)(E))\,.
	\end{equation*}
	Defining $F\eqdef(Id+\xi)(E)$ we hence obtain $\mathcal{J}_E(\xi+\zeta)=\mathcal{J}_{F}(\zeta\circ(Id+\xi)^{-1},\zeta\circ(Id+\xi)^{-1})$, which yields
	\begin{equation*}
		\mathcal{J}''_E(\xi).(\zeta,\zeta)=\mathcal{J}_F''(0).(\zeta\circ (Id+\xi)^{-1})\,.
	\end{equation*}
	Using this with $\xi=\xi_\varphi$ and $\zeta=\xi_\psi$, we get:
	\begin{equation*}
		\mathrm{j}_E''(\varphi).(\psi,\psi)=\mathcal{J}''_E(\xi_\varphi)(\xi_\psi,\xi_\psi)=\mathcal{J}''_{E_\varphi}(0).(\xi_\psi\circ g_\varphi,\xi_\psi\circ g_\varphi)\,,
	\end{equation*}
	and we finally obtain \cref{formula_shape_hessian} by applying the structure theorem.
\end{proof}

Most of the results below rely on the following technical lemma, whose first part is contained in~\cite[Lemma 4.7]{dambrineStabilityShapeOptimization2019}.
\begin{lemma}
	Let $E$ be a bounded $\mathrm{C}^2$ set. If $\|\varphi\|_{\mathrm{C}^1(\partial E)}\to 0$ we have:
	\begin{equation*}
		\begin{aligned}
			&(i)&&\|f_\varphi-Id\|_{\mathrm{C}^{1}(\partial E)}\to 0\,,
			&&\|\nu_\varphi\circ f_\varphi-\nu\|_{\mathrm{C}^{0}(\partial E)}\to 0\,,&(iii)\\
			&(ii)&&\|g_\varphi-Id\|_{\mathrm{C}^{1}(\partial E_\varphi)}\to 0\,,
			&&\|\mathrm{Jac}_{\tau}f_\varphi-1\|_{\mathrm{C}^{0}(\partial E)}\to 0\,.&(iv)\\
		\end{aligned}
	\end{equation*}
	If $\|\varphi\|_{\mathrm{C}^2(\partial E)}\to 0$ then we also have:
	\begin{equation*}
		\begin{aligned}
			&(v)&&\|H_\varphi\circ f_\varphi-H\|_{\mathrm{C}^{0}(\partial E)}\to 0\,,
			&&\|B_\varphi\circ f_\varphi-B\|_{\mathrm{C}^{0}(\partial E)}\to 0\,.&(vi)
		\end{aligned}
	\end{equation*}
	Moreover, the following holds:
	\begin{equation*}
		\begin{aligned}
			(a)~~&\underset{\|\varphi\|_{\mathrm{C}^1(\partial E)}\to 0}{\mathrm{lim}}~\underset{\psi\in\mathrm{L}^2(\partial E)\setminus\{0\}}{\mathrm{sup}}~\frac{\|(\xi_\psi\circ g_\varphi)_{\tau_\varphi}\|_{\mathrm{L}^2(\partial E_\varphi)}}{\|\psi\|_{\mathrm{L}^2(\partial E)}}=0\,,\\
			(b)~~&\underset{\|\varphi\|_{\mathrm{C}^1(\partial E)}\to 0}{\mathrm{lim}}~\underset{\psi\in\mathrm{H}^1(\partial E)\setminus\{0\}}{\mathrm{sup}}\frac{\left|\|\nabla_{\tau_\varphi}(\xi_\psi\circ g_\varphi\cdot\nu_\varphi)\|_{\mathrm{L}^2(\partial E_\varphi)}-\|\nabla_\tau\psi\|_{\mathrm{L}^2(\partial E)}\right|}{\|\psi\|_{\mathrm{H}^1(\partial E)}}=0\,,\\
			(c)~~&\underset{\|\varphi\|_{\mathrm{C}^2(\partial E)}\to 0}{\mathrm{lim}}~\underset{\psi\in\mathrm{H}^1(\partial E)\setminus\{0\}}{\mathrm{sup}}\frac{\|Z_{\varphi,\psi}\|_{\mathrm{L}^1(\partial E_\varphi)}}{\|\psi\|_{\mathrm{H}^1(\partial E)}^2}=0\,.
		\end{aligned}
	\end{equation*}
	\label{lemma_aux_c2}
\end{lemma}
\begin{proof}
	See \cite[Lemma 4.7]{dambrineStabilityShapeOptimization2019} for a proof of the results stated in the first part of the lemma. To prove $(a)$ we use the fact that
	\begin{equation*}
		\begin{aligned}
			\|(\xi_\psi\circ g_\varphi)_{\tau_\varphi}\|_{\mathrm{L}^2(\partial E_\varphi)}^2&=\int_{\partial E}(\nu\circ g_\varphi)^2_{\tau_\varphi}\circ f_\varphi~\mathrm{Jac}_{\tau}f_\varphi\,\psi^2\,d\mathcal{H}^1
			\\&\leq \|\mathrm{Jac}_{\tau}f_\varphi\|_{\mathrm{C}^0(\partial E)}\,\|(\nu\circ g_\varphi)_{\tau_\varphi}\circ f_\varphi\|_{\mathrm{C}^0(\partial E)}^2\,\|\psi\|_{\mathrm{L}^2(\partial E)}^2\\
			&=\|\mathrm{Jac}_{\tau}f_\varphi\|_{\mathrm{C}^0(\partial E)}\,\|\nu-(\nu\cdot \nu_\varphi\circ f_\varphi)\,\nu_\varphi\circ f_\varphi\|_{\mathrm{C}^0(\partial E)}^2\,\|\psi\|_{\mathrm{L}^2(\partial E)}^2\,.
		\end{aligned}
	\end{equation*}
	which, using $(i)$, $(iii)$ and $(iv)$, yields the result.
	
	To prove $(b)$, we notice that:
	\begin{equation*}
		\begin{aligned}
			\nabla_{\tau_{\varphi}}(\xi_\psi\circ g_{\varphi}\cdot\nu_{\varphi})&=\left[c^1_\varphi\,\psi\circ g_\varphi +c^2_\varphi\cdot\nabla_\tau
			\psi\circ g_\varphi\right]\tau_\varphi\,,
		\end{aligned}
	\end{equation*}
	with $\tau=\nu^\perp$, $\tau_\varphi=\nu_\varphi^\perp$\footnote{These two vectors are defined as the application of the rotation of angle $\pi/2$ to $\nu$ and $\nu_\varphi$.} and
	$$c^1_\varphi\eqdef \tau\circ g_\varphi\cdot\nu_\varphi~(J_{g_\varphi}\,\tau_\varphi)\cdot\tau\circ g_\varphi+\tau_\varphi\cdot\nu\circ g_\varphi\,,~~~~~~c^2_\varphi\eqdef\nu\circ g_\varphi\cdot\nu_\varphi~(J_{g_\varphi}\,\tau_\varphi)\,.$$
	We hence obtain $$|\nabla_{\tau_{\varphi}}(\xi_\psi\circ g_{\varphi}\cdot\nu_{\varphi})\circ f_\varphi\,\mathrm{Jac}_{\tau}f_\varphi-\nabla_\tau \psi|\leq c_\varphi\,(|\psi|+|\nabla_\tau\psi|)$$
	with~$c_\varphi$ independant of $\psi$. Moreover, using $(ii)$ and $(iii)$, we have:
	$$\underset{\|\varphi\|_{\mathrm{C}^1(\partial E)}\to 0}{\mathrm{lim}}~\|c_\varphi\|_{\mathrm{C}^{0}(\partial E)}\to 0\,.$$
	Denoting $\mathcal{A}\eqdef \left|\|\nabla_{\tau_\varphi}(\xi_\psi\circ g_\varphi\cdot\nu_\varphi)\|_{\mathrm{L}^2(\partial E_\varphi)}-\|\nabla_\tau\psi\|_{\mathrm{L}^2(\partial E)}\right|$, this finally yields
	\begin{equation*}
		\begin{aligned}
			\mathcal{A}&\leq\|\nabla_{\tau_{\varphi}}(\xi_\psi\circ g_{\varphi}\cdot\nu_{\varphi})\circ f_\varphi\,\mathrm{Jac}_{\tau}f_\varphi-\nabla_\tau \psi\|_{\mathrm{L}^2(\partial E)}\\
			&\leq \sqrt{2}~\|c_\varphi\|_{\mathrm{C}^{0}(\partial E)}\,\|\psi\|_{\mathrm{H}^1(\partial E)}\,.
		\end{aligned}
	\end{equation*}
	
	We now prove $(c)$. Since
	\begin{equation*}
		\begin{aligned}
			\|B_{\varphi}((\xi_\psi\circ g_\varphi)_{\tau_\varphi},(\xi_\psi\circ g_\varphi)_{\tau_\varphi})\|_{\mathrm{L}^1(\partial E_\varphi)}&\leq \|B_\varphi\|_{\mathrm{C}^{0}(\partial E_\varphi)}\,\|(\xi_\psi\circ g_\varphi)_{\tau_\varphi}\|_{\mathrm{L}^{2}(\partial E_\varphi)}^2
		\end{aligned}
	\end{equation*}
	and
	\begin{equation*}
		\begin{aligned}
			\mathcal{B}&\leq \|\nabla_{\tau_\varphi}(\xi_\psi\circ g_\varphi\cdot\nu_\varphi))\|_{\mathrm{L}^{2}(\partial E_\varphi)}\,\|(\xi_\psi\circ g_\varphi)_{\tau_\varphi}\|_{\mathrm{L}^{2}(\partial E_\varphi)}
		\end{aligned}
	\end{equation*}
	with $\mathcal{B}\eqdef \|(\nabla_{\tau_\varphi}(\xi_\psi\circ g_\varphi\cdot\nu_{\varphi}))\cdot(\xi_\psi\circ g_\varphi)_{\tau_\varphi}\|_{\mathrm{L}^1(\partial E_\varphi)}$, we get the result.
\end{proof}

Using the above result, we now prove the continuity of $\varphi\mapsto \mathrm{j}''_E(\varphi)$ by proving the continuity of the two terms appearing in its expression. We recall that, if $E$ is a (real) vector space, we denote by $\mathcal{Q}(E)$ the set of quadratic forms over $E$.
\begin{proposition}
	\label{lemma_cont_shape_hessian_1}
	If $E$ is a bounded $\mathrm{C}^2$ set and $\mathrm{p}_E:\varphi\mapsto P(E_\varphi)$, the mapping
	$$\! \begin{aligned}[t]
		\mathrm{p}''_E : \mathrm{C}^{2}(\partial E) &  \to \mathcal{Q}(\mathrm{H}^1(\partial E))  \\
		\varphi & \mapsto \mathrm{p}''_E(\varphi)
	\end{aligned}$$
	is continuous at $0$.
\end{proposition}
\begin{proof}
	Using \Cref{lemma_aux_c2}, for every $\varphi\in\mathrm{C}^2(\partial E)$ in a neighborhood of $0$ and $\psi\in\mathrm{H}^1(\partial E)$, we obtain:
	$$\mathrm{p}''_E(\varphi).(\psi,\psi)-\mathrm{p}''_E(0).(\psi,\psi)=\mathcal{A}+\mathrm{p}'_{E_\varphi}(0).(Z_{\varphi,\psi})\,,$$
	with $\mathcal{A}\eqdef \mathrm{p}''_{E_\varphi}(0).((\xi_\psi\circ g_{\varphi})\cdot\nu_{\varphi},(\xi_\psi\circ g_{\varphi})\cdot\nu_{\varphi})-\mathrm{p}''_E(0).(\psi,\psi)$. Now, we also have:
	\begin{equation*}
		\begin{aligned}
			\mathcal{A}=\|\nabla_{\tau_\varphi}(\xi_\psi\circ g_\varphi\cdot\nu_\varphi)\|_{\mathrm{L}^2(\partial E_\varphi)}^2-\|\nabla_\tau\psi\|_{\mathrm{L}^2(\partial E)}^2\,,
		\end{aligned}
	\end{equation*}
	and using \Cref{lemma_aux_c2} we obtain
	$$\underset{\|\varphi\|_{\mathrm{C}^2(\partial E)}\to 0}{\mathrm{lim}}~\underset{\psi\in \mathrm{H}^1(\partial E)\setminus\{0\}}{\mathrm{sup}}\frac{\left|\mathrm{p}''_{E_\varphi}(0).((\xi_\psi\circ g_{\varphi})\cdot\nu_{E_\varphi},(\xi_\psi\circ g_{\varphi})\cdot\nu_{E_\varphi})-\mathrm{p}''_E(0).(\psi,\psi)\right|}{\|\psi\|_{\mathrm{H}^1(\partial E)}^2}=0\,.$$
	Moreover
	\begin{equation*}
		\begin{aligned}
			|\mathrm{p}'_E(0).(Z_{\varphi,\psi})|&\leq \|H_\varphi\|_{\mathrm{L}^{\infty}(\partial E_\varphi)}\,\|Z_{\varphi,\psi}\|_{\mathrm{L}^1(\partial E_\varphi)}\,,
		\end{aligned}
	\end{equation*}
	and \Cref{lemma_aux_c2} allows to conclude.
\end{proof}

\begin{proposition}
	\label{lemma_cont_shape_hessian_2}
	If $E$ is a bounded $\mathrm{C}^2$ set, $\eta\in\mathrm{C}^1(\RR^2)$ and $\mathrm{g}_E:\varphi\mapsto \int_{E_\varphi}\eta$, the mapping
	$$\! \begin{aligned}[t]
		\mathrm{g}''_E : \mathrm{C}^{2}(\partial E) &  \to \mathcal{Q}(\mathrm{H}^1(\partial E))  \\
		\varphi & \mapsto \mathrm{g}''_E(\varphi)
	\end{aligned}$$
	is continuous at $0$.
\end{proposition}
\begin{proof}
	We proceed as in \Cref{lemma_cont_shape_hessian_1}. Defining
	$$\mathcal{A}\eqdef \mathrm{g}''_{E_\varphi}(0).((\xi_\psi\circ
	g_{\varphi})\cdot\nu_{E_\varphi},(\xi_\psi\circ g_{\varphi})\cdot\nu_{E_\varphi})$$
	we have:
	\begin{equation*}
		\begin{aligned}
			\mathcal{A}&=\int_{\partial
				E_{\varphi}}\left[H_{\varphi}\,\eta+\frac{\partial\eta}{\partial\nu_{\varphi}}\right]((\psi\,\nu)\circ
			g_{\varphi}\cdot\nu_{E_\varphi})^2\,d\mathcal{H}^1\\
			&=\int_{\partial E}\left[H_{\varphi}\,\eta+\frac{\partial\eta}{\partial\nu_{\varphi}}\right]\circ f_{\varphi}\,(\nu\cdot\nu_{\varphi}\circ f_{\varphi})^2\,\mathrm{Jac}_\tau f_\varphi\,\psi^2\,d\mathcal{H}^1.
		\end{aligned}
	\end{equation*}
	This yields:
	\begin{equation*}
		\begin{aligned}
			\frac{\left|\mathrm{g}''_{E_\varphi}(0).((\xi_\psi\circ g_{\varphi})\cdot\nu_{\varphi},(\xi_\psi\circ g_{\varphi})\cdot\nu_{\varphi})-
				\mathrm{g}''_E(0).(\psi,\psi)\right|}{\|\psi\|_{\mathrm{L}^2(\partial E)}^2}\leq c_\varphi\,,
		\end{aligned}
	\end{equation*}
	with $$c_\varphi\eqdef\left\|\left[H_{\varphi}\,\eta+\frac{\partial\eta}{\partial\nu_{\varphi}}\right]\circ
	f_{\varphi}\,(\nu\cdot\nu_{\varphi}\circ f_{\varphi})^2\,\mathrm{Jac}_{\tau}f_\varphi-\left[H\,\eta+\frac{\partial\eta}{\partial\nu}\right]\right\|_{\infty}.$$
	Using \Cref{lemma_aux_c2} we obtain
	$$\underset{\|\varphi\|_{\mathrm{C}^2(\partial E)\to 0}}{\mathrm{lim}}~\underset{\psi\in \mathrm{H}^1(\partial E)
		\setminus\{0\}}{\mathrm{sup}}\frac{\left|\mathrm{g}''_{E_\varphi}(0).((\xi_\psi\circ g_{\varphi})\cdot\nu_{\varphi},(\xi_\psi\circ g_{\varphi})\cdot\nu_{\varphi})-
		\mathrm{g}''_E(0).(\psi,\psi)\right|}{\|\psi\|_{\mathrm{H}^1(\partial E)}^2}=0\,.$$
	Moreover
	\begin{equation*}
		\begin{aligned}
			|\mathrm{g}'_{E_\varphi}(0).(Z_{\varphi,\psi})|&\leq \|\eta\|_{\infty}\,\|Z_{\varphi,\psi}\|_{\mathrm{L}^1(\partial E_\varphi)}\,,
		\end{aligned}
	\end{equation*}
	and using again \Cref{lemma_aux_c2} we finally obtain the result.
\end{proof}

Proof of \Cref{lemma_conv_shape_hessian}:
\begin{proof}
	Since $\left|(\mathrm{j}_E-\mathrm{j}_{0,E})''(\varphi).(\psi,\psi)\right|\leq c^1_\varphi+c^2_\varphi$ with
	\begin{equation*}
		\begin{aligned}
			c^1_\varphi&\eqdef\left|\int_{\partial E_\varphi}\left(H_\varphi\,(\eta-\eta_0)+\frac{\partial\,(\eta-\eta_0)}{\partial \nu_\varphi}\right)(\xi_\psi\circ g_\varphi\cdot \nu_\varphi)^2\,d\mathcal{H}^1\right|,\\
			c^2_\varphi&\eqdef\left|\int_{\partial E_\varphi}(\eta-\eta_0)\,Z_{\varphi,\psi}\,d\mathcal{H}^1\right|,
		\end{aligned}
	\end{equation*}
	the result readily follows from \Cref{lemma_aux_c2}.
\end{proof}